\documentclass{article}
\usepackage[english]{babel}
\usepackage[utf8]{inputenc}

\usepackage{bm}
\usepackage{amsmath}
\usepackage{graphicx}
\usepackage{xcolor}
\usepackage{multirow}
\usepackage{subcaption}
\usepackage{amsfonts}
\graphicspath{ {./figs/} }

\date{}

\title{Multiscale solver for multi-component reaction-diffusion systems in heterogeneous media}
\author{
Maria Vasilyeva, 
Alexey Sadovski
and D. Palaniappan \\
Department of Mathematics and Statistics\\ 
Texas A\&M University - Corpus Christi, Corpus Christi, Texas - 78412
   \thanks{{\tt maria.vasilyeva@tamucc.edu}, {\tt alexey.sadovski@tamucc.edu},{\tt devanayagam.palaniappan@tamucc.edu}}
}

\begin{document}

\maketitle

\begin{abstract}
Coupled nonlinear system of reaction-diffusion equations describing multi-component (species) interactions with heterogeneous coefficients is considered. Finite volume method based approximation for the space is used to construct semi-discrete form for the computation of numerical solutions. Two techniques for time approximations, namely, a fully implicit (\textit{FI}) and a semi-implicit (\textit{SI}) schemes are examined.  The fully implicit scheme is constructed using Newton's method and leads to the coupled system of equations on each nonlinear and time iterations which is computationally rather expensive. In order to minimize the latter hurdle, an efficient and fast multiscale solver is proposed for reaction-diffusion systems in heterogeneous media. To construct fast solver, we apply a semi-implicit scheme that leads to an uncoupled system for each individual component. 
Problems in heterogeneous domains require a very fine grid for accurate solutions of large systems of equations at each time step iteration. Here we present a multiscale model reduction technique to reduce the size of the discrete system.  Multiscale solver is  based on the uncoupled operator of the problem and constructed by the use of Generalized Multiscale Finite Element Method (GMsFEM). In GMsFEM we use a diffusion part of the operator and construct multiscale basis functions by solving spectral problems in each local domain associated with the coarse grid nodes.  We collect multiscale basis functions to construct a projection/prolongation matrix and generate reduced order model on the coarse grid for fast solution. Moreover, the prolongation operator is used to reconstruct a fine-scale solution and accurate approximation of the reaction part of the problem which then leads to a very accurate and computationally effective multiscale solver. We provide numerical results for two species competition test problems in two-dimensional domain with heterogeneous inclusions.
Our computed solutions predict strong competition between the two species and demonstrate that diffusion can dictate the dominance of one species over the other in heterogeneous environments.   
We investigate the influence of number of the multiscale basis functions to the method accuracy and ability to work with different values of the diffusion coefficients. 
\end{abstract}

\section{Introduction}
Reaction-diffusion (RD) system of equations often occur in biological and chemical modeling \cite{Keener,Brown} and have been found useful in other physical sciences and finance discipline as well. Modeling with RD systems is an efficient technique for describing processes relating to multi-species dynamics since it is flexible for more realistic environments. In many circumstances, it is common to assume that RD model system is spatially homogeneous so that existing mathematical tools can be utilized to explore various scenarios, especially in the case of two-species competition systems (see \cite{Parshad} and references therein). But studies, for instance \cite{Steinberg,Alonso}, have shown that spatial heterogeneity can affect spatial structure and dynamics significantly. Spatially distributed competitive species models can provide guidance to various control strategies \cite{Braverman}. Fairly general RD models describing multi-species competitions in heterogeneous domains are mathematically challenging due to additional complexity due to diffusion coefficient variations. Despite the importance, a systematic computationally efficient technique for the determination of numerical solutions is still missing which is a key focus of the present investigation.
     
%
In mathematical ecology, the reaction-diffusion models are used to describe population of the species \cite{okubo2001diffusion}.  Spatial-temporal Lotka - Volterra models, a variant of the RD system, are used to describe population of multiple interacting communities and take into account the influence of spatial dependence on various properties such as biodiversity and food web structure, among many others \cite{montagna2017modeling}. 
A system of partial differential equations (PDE) is used to describe population variation in time and space. The resulting mathematical model is described by a coupled system of time-dependent nonlinear reaction–diffusion equations accommodating the effects of spatial heterogeneity \cite{chen2011spatio}.  
In simulations of the subsurface flow and transport, the reactive transport models are used for the purpose of biochemical reactions \cite{keating2013co2, xu2011toughreact, moortgat2020higher, steefel2018approaches}.  
Actually, the reactive transport leads to a sharp reaction interface and requires specific numerical techniques for accurate approximation and effective solutions of the PDE system.
 
In the current manuscript we consider the coupled system of RD equations describing multi-species interactions/competitions in heterogeneous two-dimensional domain $\Omega \subset R^2$ given in section 2 (see \eqref{m}). 
%
Precisely, we consider the system of PDEs in a heterogeneous domain with multiple inclusions.  
We provide a construction strategy for the discrete system on a grid that resolves heterogeneity (inclusions) at the grid level (fine grid, $\mathcal{T}_h$) using semi-implicit time approximation scheme and finite volume method for space discretization.  The size of the resulting linear system that we solve at each time layer is to the number of fine grid cells.  To obtain accurate solution in heterogeneous media with inclusions, it is required to use a sufficiently fine grid which leads to computationally expensive large discrete system of equations.

Computation of solutions of problems in heterogeneous media is complicated due to the multiscale nature of the processes.  Numerical solutions of such problems are actually computationally expensive because the complexity involved in the mesh resolution. Indeed, the better mesh resolution is needed for accurate simulations with high resolution of the multiscale features on the grid level.  
The multiscale methods and homogenization techniques are widely used for such problems to reduce the size of the system by the construction of the coarse-scale approximation. The literature on multiscale methods is extensive and interest in the topic is rapidly growing. Many model reduction techniques, such as numerical homogenization, upscaling and multiscale methods \cite{allaire1992homogenization, efendiev2009multiscale, efendiev2013generalized, ReducedCon} have been proposed in the literature. 
Generalized multiscale finite element methods (GMsFEM) for problems in perforated heterogeneous domains were considered in the first author's earlier collaborative work \cite{chung2016generalized}.  In the cited reference, the authors presented the construction of coarse grid approximation based on the continuous Galerkin approach and described the generation of multiscale space for elliptic, elastic, and Stokes problems.  
In \cite{vasilyeva2019upscaling}, the upscaling method is given for problems in perforated domains with non-homogeneous boundary conditions on the perforations using Non-Local Multi-Continuum method (NLMC). 
In \cite{vasilyeva2021multiscale}, the authors provided a multiscale solver based on the discontinuous Galerkin GMsFEM for the numerical solutions of flow fields and reactive transport in thin domains. 
Online residual based  multiscale technique for problems in heterogeneous media has been demonstrated in \cite{chung2015residual, chung2016online}.  
Sequential homogenization techniques for the solution of reactive transport systems in porous media are discussed in  \cite{korneev2016sequential}. In particular, the hierarchy of effective equations that sequentially carry the smallest scale information through the intermediate scales up to the macroscale was highlighted in the latter cited work. 
Another application of the sequential upscaling based on GMsFEM for solution of the problems in heterogeneous and perforated domains is presented in \cite{chung2016reiterated}. 
In \cite{battiato2011hybrid}, hybrid models for the simulation of reactive transport in porous and fractured media are given. The  iterative hybrid numerical method described in \cite{battiato2011hybrid} links the pore and continuum scales.  The technique has been employed in the modeling of transport in a fracture with chemically reactive walls.
Recently, for effective solution of the multiscale problems machine learning techniques have been applied. 
Implementation of the convolutional neural networks for the prediction of effective medium properties is discussed in \cite{vasilyeva2021machine, vasilyeva2019convolutional}. The novel nonlinear upscaling technique for transport and flow problems in heterogeneous domains is given in \cite{vasilyeva2020learning}. 
In \cite{prasianakis2020neural}, the authors considered multiscale modeling of geochemical processes and discussed a neural network approach for reactive transport simulations.

The primary goal of this paper is the construction of a fast, efficient and accurate solver on the coarse grid $\mathcal{T}_H$ with $H >> h$, where $h$ is the mesh size. In our approach, the coarse grid approximation is based on the Generalized Multiscale Finite Element Method (GMsFEM).  As will be seen later, in GMsFEM, local multiscale basis functions are constructed in local domains linked to the coarse grid node.  Both offline and online stages are contained in the proposed multiscale algorithm to solve local spectral and reduced order problems. One of the advantages of the algorithm is the construction of the  basis functions based on the diffusion operators for each species separately.  We use fully-implicit and semi-implicit schemes for time approximation and compare their performance.  We employ finite volume based approximation to resolve spatial heterogeneity due to the presence of inclusions. As an application, we consider two species competition model in a heterogeneous domain and use our algorithm to generate accurate numerical solutions to demonstrate the power of the multiscale solver.

The paper is organized as follows. In Section 2, we describe the problem formulation with fine-scale approximation using finite volume method and provide details of fully implicit and semi-implicit schemes used for time approximations. We present the multiscale basis functions construction based on the diffusion operator and describe multiscale solver algorithm for reaction-diffusion systems in Section 3. Discussions of offline and online stages and a brief outline of coarse grid construction techniques are given in the same Section. Numerical results for two species competition model are presented in Section 4. The Section also includes graphical illustrations of fine grid and multiscale solver solutions, solution averages and error calculations. The paper ends with a conclusion in Section 5.

\section{Description of the mathematical model}
We consider the following reaction-diffusion (RD) system corresponding to $k$ competing species in a heterogeneous domain $\Omega$ written in the form
\begin{equation}
\label{m}
\begin{split}
\frac{\partial u^k}{\partial t} 
- \nabla \cdot (\varepsilon^k(x) \nabla u^k) = 
R^k(u^1 ..., u^L), 
\quad k = 1,..,L, 
\quad x \in \Omega, \quad t>0,
\end{split}
\end{equation}
with 
\[
R^k(u^1, ..., u^L) =  r^k(x)   u^k(1 - u^k)
 - \sum_{l \neq k} \alpha^{kl}(x) u^k u^l,
\]
where 
$u^k$ represents the population of the $k$-th species, 
$\varepsilon^k$ is the spatially dependent diffusion coefficient, 
$r^k$ and $\alpha^{kl}$  are the $k$-th population reproductive growth rate and the interaction coefficient due to competition ($u^l$ compete with $u^k$), respectively. Note that the growth rates and the competition coefficients also vary spatially in our model. The system of equations \eqref{m} is subject to the following initial condition 

\begin{equation}
\label{m-ic}
u^k = u^k_0, \quad x \in \Omega, \quad t = 0,
\end{equation}
and the Neumann/zero outflow boundary condition (i.e., no migration across the boundary)
\begin{equation}
\label{m-bc}
\nabla u^k \cdot n = 0, \quad x  \in \partial \Omega,
\end{equation}
where $n$ is the outer normal vector to the boundary $\partial \Omega$. The coupled system \eqref{m} together with \eqref{m-ic} and \eqref{m-bc} constitute a nonlinear initial-boundary value problem (IBVP) describing the competition of $k$ species in an inhomogeneous medium. 

Further, we consider the heterogeneous domain $\Omega$ to be composed of a distribution of two types of domains where the species react and diffuse. We let $\Omega = \Omega_m \cup \Omega_c$, where $\Omega_m$ is the background domain and $\Omega_c$ denotes the subdomain of the inclusions in the heterogeneous medium
In order to emphasize heterogeneity, we set:
\[
\varepsilon^k(x) = \left\{
\begin{matrix}
\varepsilon^k_m, & \in \Omega_m \\ 
\varepsilon^k_c, & \in \Omega_c
\end{matrix} \right., 
\quad 
r^k(x) = \left\{
\begin{matrix}
r^k_m, & \in \Omega_m \\ 
r^k_c & \in \Omega_c
\end{matrix}\right., 
\quad 
\alpha^{lk}(x) = \left\{
\begin{matrix}
\alpha^{lk}_m, & \in \Omega_m \\ 
\alpha^{lk}_c & \in \Omega_c
\end{matrix} \right..
\]
Here $\varepsilon^k_m$, $\varepsilon^k_c$, $r^k_m$, $r^k_c$, $\alpha^{lk}_m$ and $\alpha^{lk}_c$ are constants characterizing the  properties of the background and subdomains. 
We remark that, the shape of the inclusions does not affect our numerical algorithm proposed herein and therefore may be chosen arbitrary. In the following subsections, we describe space and time approximations for our numerical scheme. 

\subsection{Spatial approximation on the fine grid}

\begin{figure}[h!]
\centering
\includegraphics[width=0.7\linewidth]{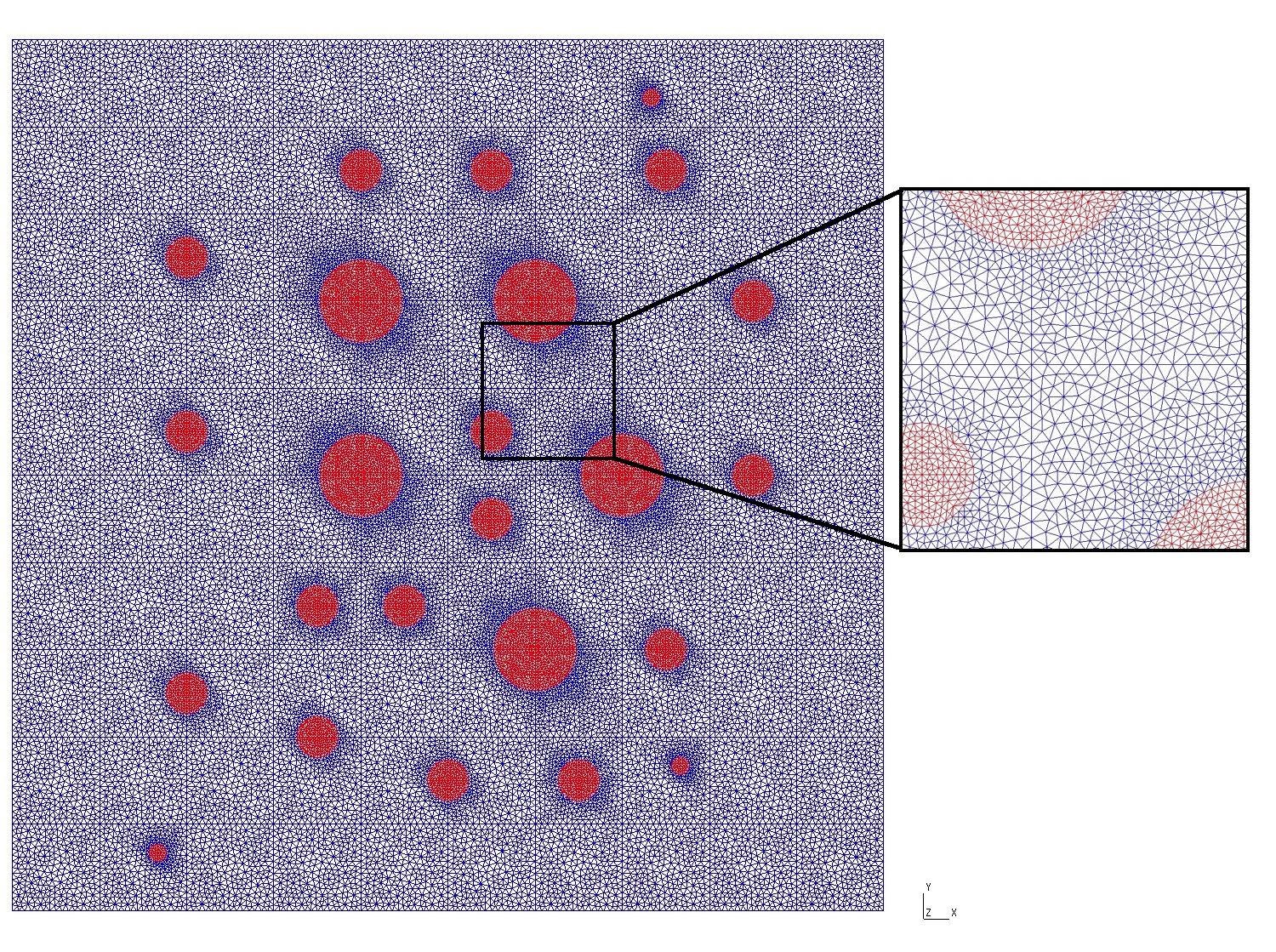}
\caption{
Fine grid that resolves heterogeneity on the grid level. Blue color: background subdomain $\Omega_m$. Red color: subdomain of the circle inclusions $\Omega_c$}
\label{domain}
\end{figure}

In order to determine the numerical solution of the IVBP set up in \eqref{m} together with the initial and boundary conditions \eqref{m-ic}-\eqref{m-bc}, we first construct the grid that resolves heterogeneity (due to the inclusions) at the grid level (fine grid, see Figure \ref{domain}). To this end, let $\mathcal{T}_h$ be the triangulation of the domain $\Omega$, with mesh size $h$, taken as
\[
\mathcal{T}_h = \cup_i K_i,
\]
where $K_i$ is the fine grid of the $i^{th}$ cell, with $i = 1,...,N$, and $N$ being the number of fine grid cells.  
Let $e_{ij}$ be the interface (facet) between the two cells $K_i$ and $K_j$, with $K_i \cap K_j \neq \O$.

In our numerical scheme, we employ a finite volume approximation by space. Thus, for each cell $K_i$, we have
\begin{equation}
\int_{K_i} \frac{ \partial u^k }{\partial t} dx
+  \sum_{j} \int_{e_{ij}} \varepsilon^k \nabla u^k \cdot n_{ij} \ ds
= \int_{K_i} R^k (u^1 ..., u^L) \ dx.
\end{equation}
 
Now, let $u_i^k$ be the cell average value of the function $u^k$ on cell $K_i$ written as
\[
\frac{1}{|K_i|}\int_{K_i} u^k dx = u_i^k,
\]
where $|K_i|$ denotes the cell volume. For the diffusion operator, we use a regular two-point flux approximation (TPFA)
\[
\int_{e_{ij}} \varepsilon^k \nabla u^k \cdot n_{ij} \ ds \approx T^k_{ij} ( u^k_i - u^k_j ).
 \quad 
 T^k_{ij} = \varepsilon^k_{ij} \ |e_{ij}| / d_{ij}, 
\]
Here $d_{ij}$ is the distance between the two cell center points $x_i$ and $x_j$,  
$|e_{ij}|$ is the length of the interface between to cells $K_i$ and $K_j$, and  
$\varepsilon^k_{ij}$ is the harmonic average given by $\varepsilon^k_{ij} = 2/(1/\varepsilon^k_{i} +1/ \varepsilon^k_{j})$.

Therefore, for the spatial approximation, we have  
\begin{equation}
\frac{ \partial u_i^k }{\partial t} |K_i| 
+  \sum_{j} T^k_{ij} ( u^k_i - u^k_j ) 
= R^k_i (u_i^1 ..., u_i^L) |K_i|, 
\end{equation}
with 
\[
R^k_i(u_i^1 ..., u_i^L) =  r^k_i   u_i^k (1 - u_i^k)
 - \sum_{l \neq k} \alpha^{kl}_i u_i^k  u_i^l.
\]

\subsection{Time approximation} 
For the time approximation, we propose the following two schemes, namely, a fully implicit (FI) scheme giving a coupled system and a semi-implicit (SI) scheme leading to an uncoupled system of equations \cite{vabishchevich2013additive, vabishchevich2012explicit, owolabi2014higher}.

Let $u_i^k = u^k(x_i, t_n)$ and $\check{u}_i^k = u^k(x_i, t_{n-1})$, where $t_n = n \tau$,  $n=1,2, ...$, and $\tau > 0$, be the fixed time step size. 
We first apply backward Euler's approximation for time derivative and obtain a fully implicit (\textit{FI}) scheme
\begin{equation}
\label{nd}
\frac{ u^k_i - \check{u}^k_i}{\tau} |K_i| 
+  \sum_{j} T^k_{ij} ( u^k_i - u^k_j ) 
= R^k_i (u_i^1 ..., u_i^L) |K_i|, 
\end{equation}
with 
\[
R^k_i(u_i^1 ..., u_i^L) =  
r^k_i u_i^k (1 - u_i^k) |K_i|
 - \sum_{l \neq k} \alpha^{kl}_i u_i^k  u_i^l |K_i|.
\]
\noindent
We utilize Newton's method to solve the nonlinear system of equations. Let $s$ be the nonlinear iteration number and
\[
u_i^{s+1} = u_i^s + \delta u_i, \quad s=0,1,2...
\]
Applying linearization (for the reaction term) 
\[
R^k(u^{s+1}) \approx R^k(u^s) + \sum_j (R^k(u^s))'_j \delta u^j, 
\]\[
(R^k(u))'_j = \left\{ \begin{matrix}
r^k (1 - 2 u^k) - \sum_{l \neq k} \alpha^{kl} u^l  & j = k, \\
- \alpha^{kj} u^k  & j \neq k, 
\end{matrix} \right.
\]
one obtains the following system of linear equations for each nonlinear iteration 
\begin{equation}
\label{fi}
\frac{ \delta u^k_i}{\tau} |K_i| 
+  \sum_{j} T^k_{ij} (  \delta u^k_i -  \delta u^k_j ) 
- r^k (1 - 2 u^{k,s})  \delta u^k_i  |K_i|
+ \sum_{l \neq k} \alpha^{kl} (u^{l,s} \delta u^k_i + u^{k,s}_i \delta u^l_i) |K_i|
= - F^{k,s}_i,  
\end{equation}
with
\[
F^{k,s}_i = \frac{ u^{k,s}_i - \check{u}^k_i}{\tau} |K_i| 
+  \sum_{j} T^k_{ij} ( u^{k,s}_i - u^{k,s}_j ) 
- r^k_i u_i^{k,s} (1 - u_i^{k,s}) |K_i|
+ \sum_{l \neq k} \alpha^{kl}_i u_i^{k,s}  u_i^{l,s} |K_i|.
\]
Here, we perform iterations till $|| \delta u^k ||_{L_2} < \epsilon_{nl}$ or until we reach a maximum number of the nonlinear iterations. We use the solution from previous time layer as the initial condition for the successive nonlinear iterations. 
We update the matrix as well as the right-hand side of the linear system \eqref{fi} in each time step and nonlinear iteration since both  depend on the current/preceding solution $u^{k,s}$.  We remark that updating the matrix and the right-hand side vector is computationally expensive and may lead to the large computational cost. Additionally, the convergence of the nonlinear iterations depends on the time step size and so larger size time steps require more iterations to converge. However, it should be pointed out that, Newton's method is the most accurate method for the numerical solutions of the nonlinear systems. The illustrated scheme is fully implicit and leading a way to the  solution of the coupled system of equations. 

Next, we provide a computationally effective uncoupling scheme. In this construction, we use a semi-implicit scheme (\textit{SI}) for the time approximation given by
\begin{equation}
\label{si1}
\frac{ u^k_i - \check{u}^k_i}{\tau} |K_i| +
\sum_{j} T^k_{ij} ( u^k_i - u^k_j ) 
= R^k_i (\check{u}_i^1 ..., \check{u}_i^L) |K_i|, 
\end{equation}
with 
\[
R^k_i(\check{u}_i^1 ..., \check{u}_i^L) =  
r^k_i \check{u}_i^k (1 - \check{u}_i^k) |K_i|
- \sum_{l \neq k} \alpha^{kl}_i \check{u}_i^k \check{u}_i^l |K_i|.
\]
The reaction term here is approximated by using the solution from previous time layer. In the present scheme set-up, the matrix corresponding to the linear equations is fixed (does not change with time) and further the system gets decoupled for each component (species). So, one can solve the equations for each species independently by just updating the vector on the right-hand-side. 


It will be shown later (in the numerical section) that the uncoupling technique resulting from the (\textit{SI}) scheme is  accurate and computationally effective. Our main objective of the present work is the construction and investigation of the multiscale method for generating numerical solutions of the reaction-diffusion systems on the coarse grid.  To build a computationally effective model reduction technique for problems that require space approximation in heterogeneous media, we use a multiscale finite element method and construct multiscale space for coarse grid approximation based on uncoupling time approximation procedure. 

The system \eqref{si1} can be conveniently re-written in the following matrix form
\begin{equation}
M  \frac{ u^k - \check{u}^k}{\tau}  + 
A^k u^k=   R^k,
\end{equation}
with 
\[
A^k = \{a^k_{ij}\}, \quad 
a^k_{ij} = 
\left\{\begin{matrix}
\sum_j T^k_{ij} & i = j, \\ 
-T^k_{ij} & i \neq j
\end{matrix}\right. , 
\quad 
M = \{m_{ij}\}, \quad 
m_{ij} = 
\left\{\begin{matrix}
|K_i| & i = j, \\ 
0 & i \neq j
\end{matrix}\right. , 
\]
and $R^k = \{R^k_i|K_i|\}$.
The size of the system for each species ($DOF_h$) equals the number of fine grid cells $N$. Since the number of unknowns $DOF_h$ is large, the determination of accurate numerical solutions in heterogeneous media is challenging and computationally expensive as well. As stated earlier, the main purpose of this paper is to develop a reduced order model on the coarse grid $\mathcal{T}_H$ with $H >> h$ ($H$ and $h$ here represent the fine and coarse mesh sizes). We employ a multiscale finite element method on $\mathcal{T}_H$ to construct fast and accurate multiscale solver with $DOF_H >> DOF_h$ as discussed in the next section.

\section{Multiscale model reduction using GMsFEM} 

\begin{figure}[h!]
\centering
\includegraphics[width=0.377\linewidth]{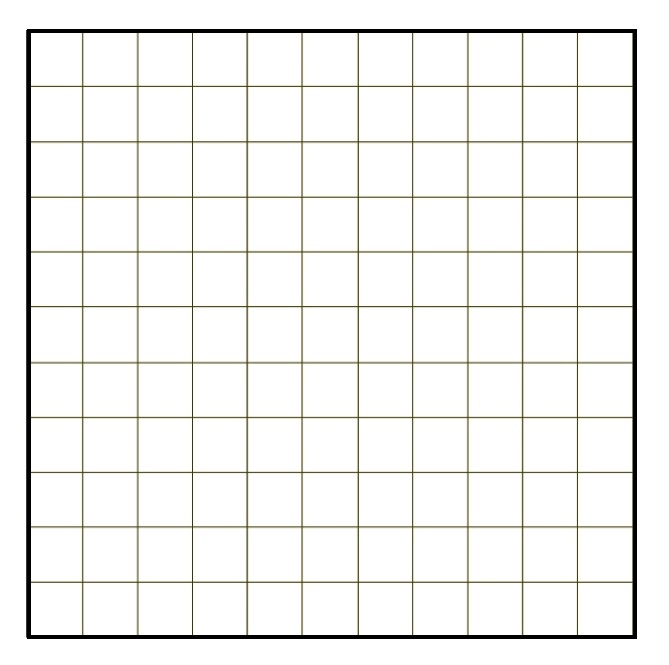}
\includegraphics[width=0.5\linewidth]{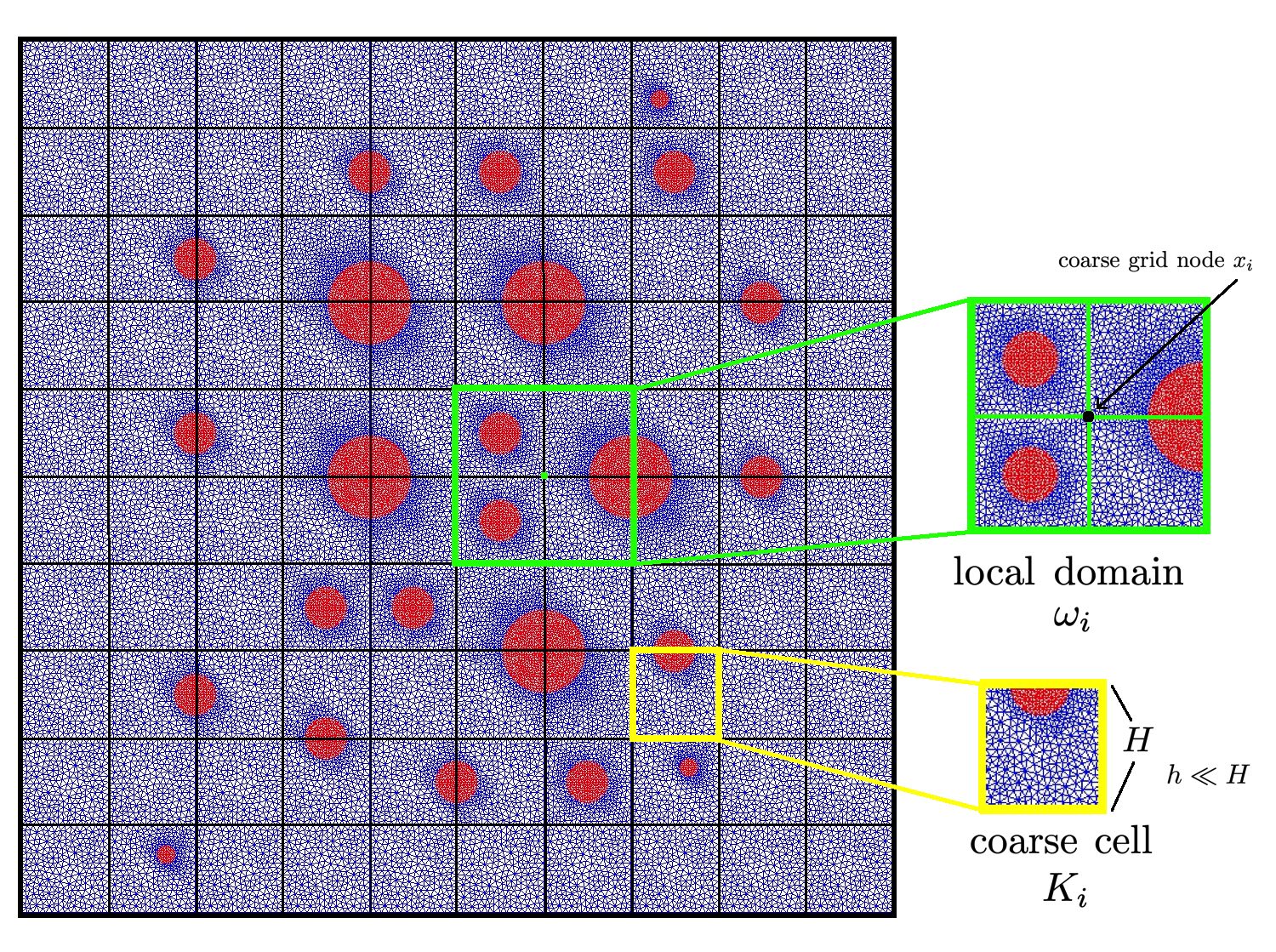}
\caption{Illustration of the $10 \times 10$ coarse grid, $\mathcal{T}_H$  (left).  
Fine grid for heterogeneous media conformed with coarse grid cell interfaces and local domain $\omega_i$ (right)}
\label{meshc}
\end{figure}

Let $\mathcal{T}_H$ be the coarse grid with coarse cell denoted by $K_i$ (see Figure \ref{meshc}). 
The key idea of the multiscale method is the solution of the problem on the coarse grid $\mathcal{T}_H$, rather than the solution on the fine grid $\mathcal{T}_h$ with $h << H$. 
For the reduction of the discrete system size, we use a Generalized Multiscale Finite Element Method (GMsFEM). 
In GMsFEM, we construct local multiscale basis functions in local domain $\omega_i$. Here $\omega_i$ is the domain related to the coarse grid node $x_i$ and designed as a combination of the several coarse cells containing the corresponding coarse grid node. We adopt conforming fine triangulation approach in our present work. In general, one can use a mesh partitioning to define coarse cells alternatively.

The GMsFEM algorithm consists of offline and online stages.  
In the offline setting, we construct local multiscale basis functions via the solution of local spectral problems in each local domain $\omega_i$. On the other hand, in the online stage, we solve a reduced order problem on the coarse grid. Specifically, the operations performed in these two stages of our algorithm are enumerated below. 

\begin{enumerate}
\item Coarse and fine grids construction and generation of local domains (\textit{Offline stage}).
\item Solution of the local spectral problems and construction of the multiscale basis functions (\textit{Offline stage}).
\item For each time iteration (\textit{Online stage}):
\begin{enumerate}
\item Update and solve the coarse grid system.
\item Project solution into fine grid.
\end{enumerate}
\end{enumerate}


We now proceed to describe our proposed multiscale solver design. We begin with the construction of the multiscale basis functions in each local domain. For the reaction-diffusion systems \eqref{m}, we construct basis functions based on the diffusion operators for each species separately. We then use them to construct (project) coarse scale system with smaller number of unknowns ($DOF_H << DOF_h$). Next we solve coarse grid reduced order model and reconstruct fine-scale solution for accurate reaction term approximation. We update the right-hand side (that is, the reaction term) based on the fine grid solution information and project it back to the coarse grid. Note that the matrix associated with the linear system, that we solve at each time step based on semi-implicit time approximation, doesn't change and therefore can be constructed at one time in the offline stage. Moreover, the multiscale basis functions do not depend on the reaction term and so can be used for multiparameteric simulations with various values for the reaction coefficients. Some mathematical details pertaining to multiscale basis functions and coarse grid system are provided in the next subsections.

\begin{figure}[h!]
\centering
\includegraphics[width=0.4\linewidth]{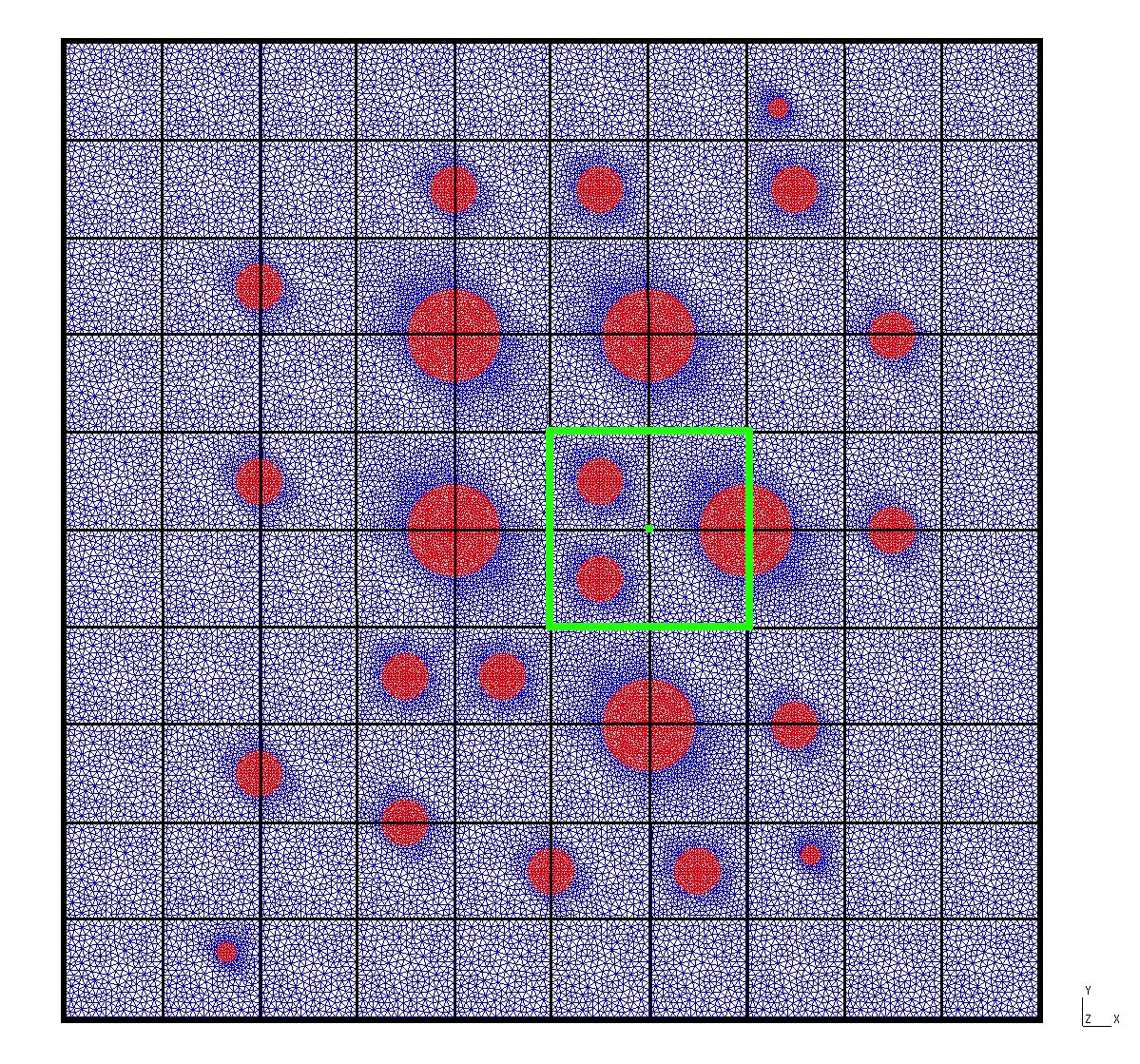}\\
\includegraphics[width=0.11\linewidth]{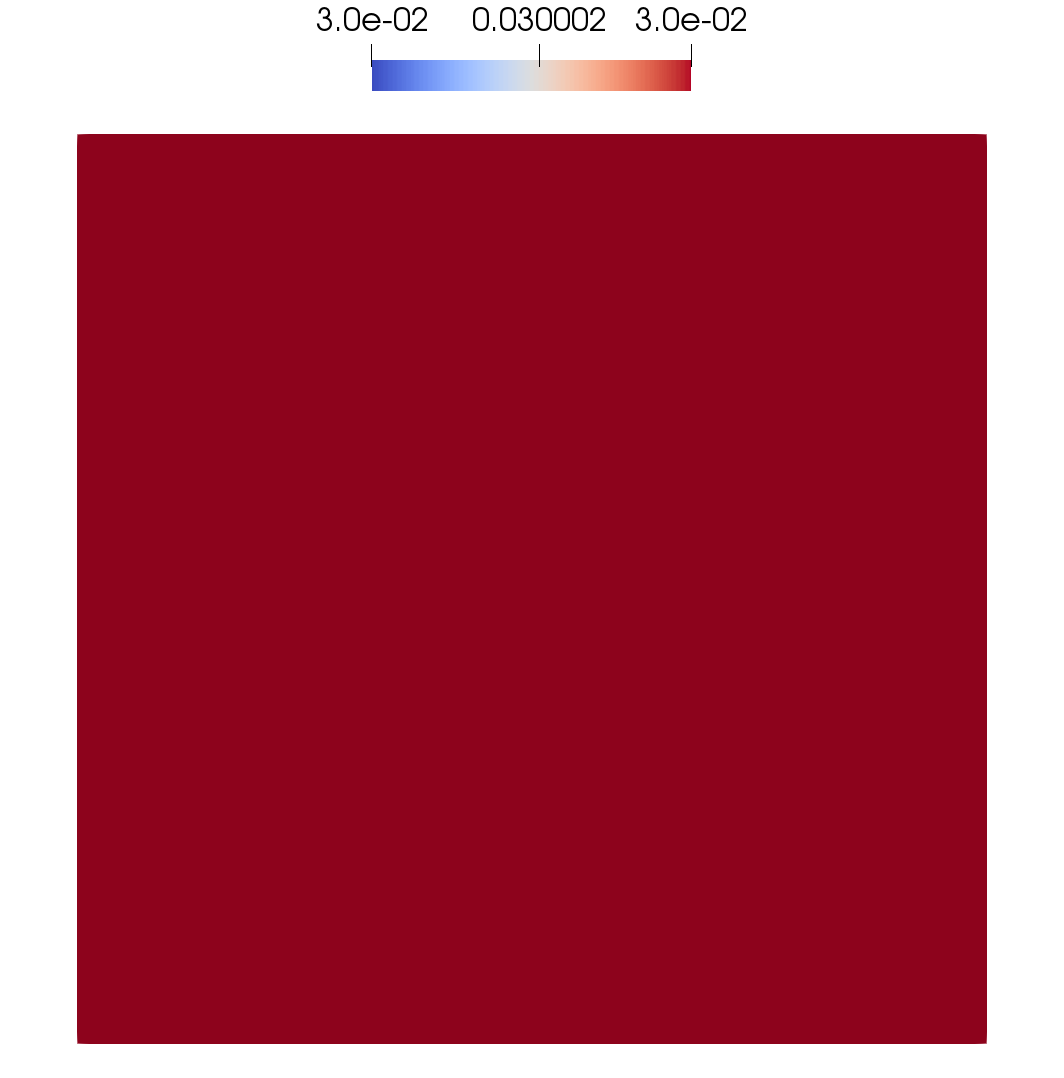}
\includegraphics[width=0.11\linewidth]{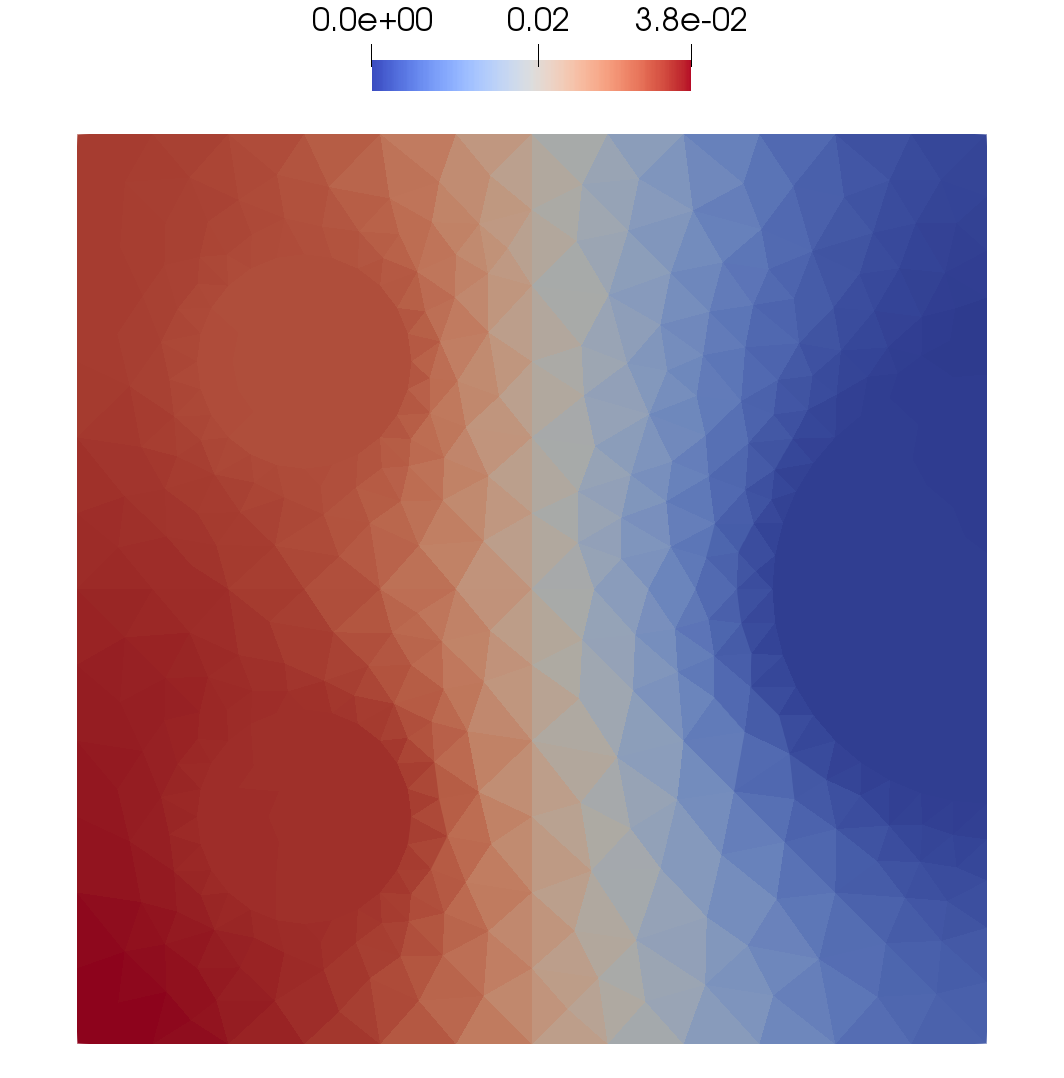}
\includegraphics[width=0.11\linewidth]{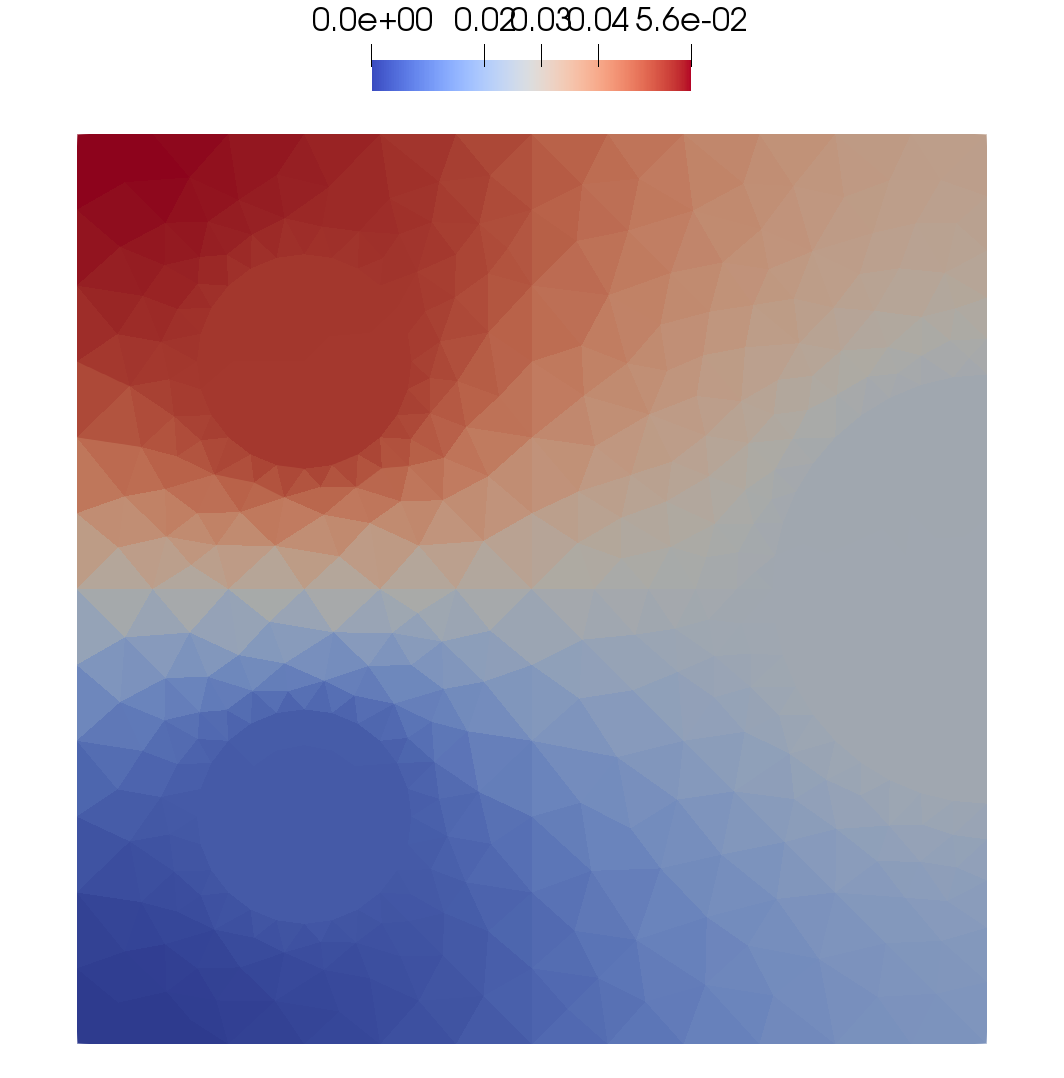}
\includegraphics[width=0.11\linewidth]{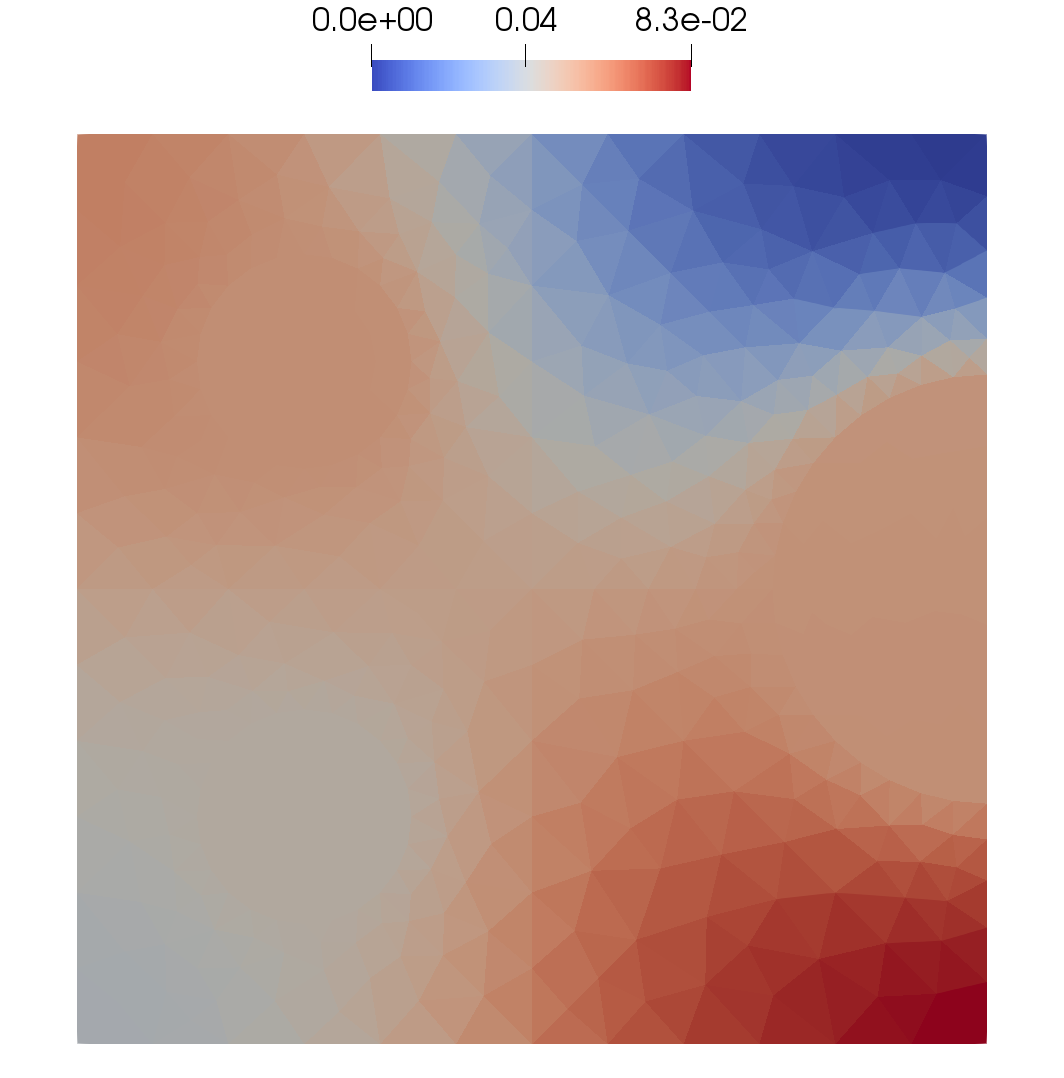}
\includegraphics[width=0.11\linewidth]{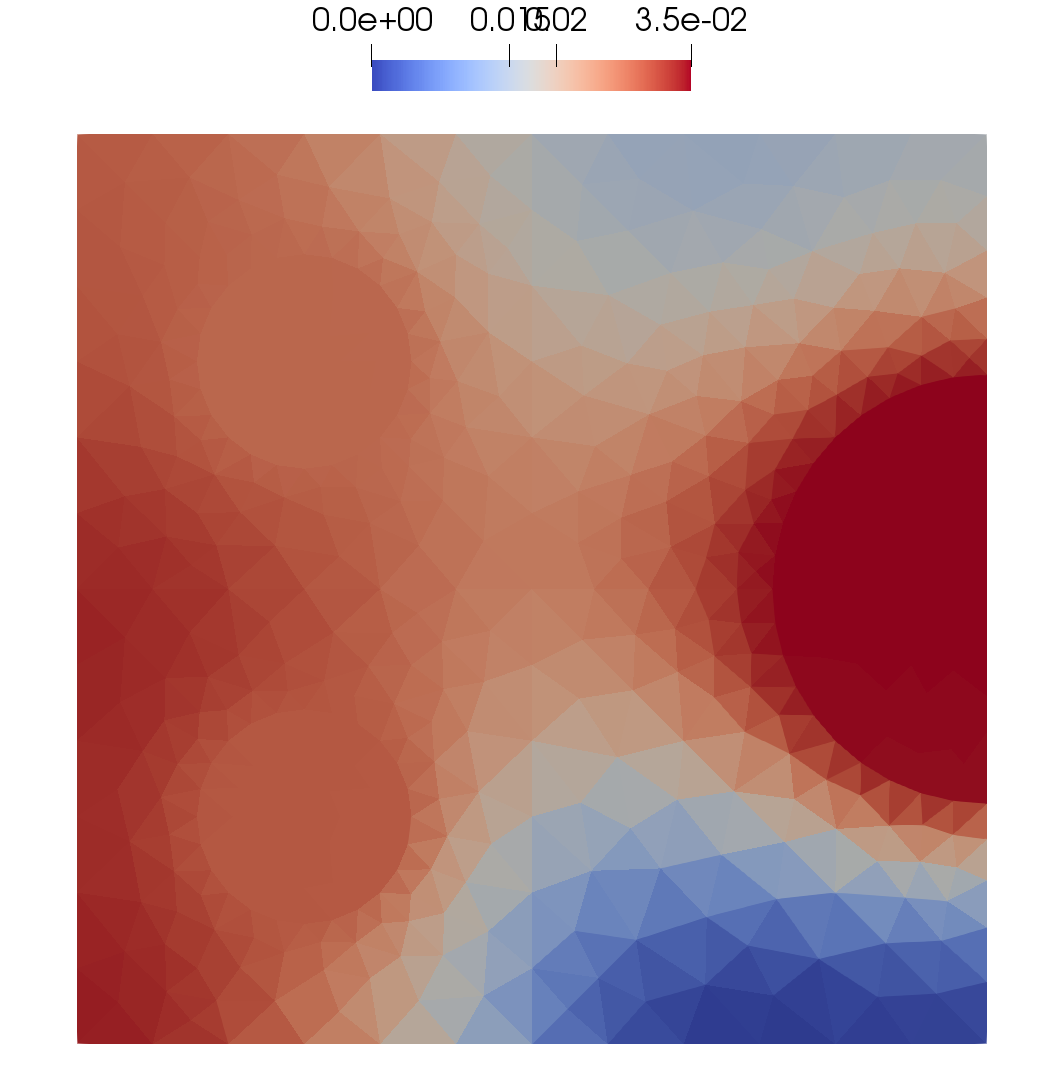}
\includegraphics[width=0.11\linewidth]{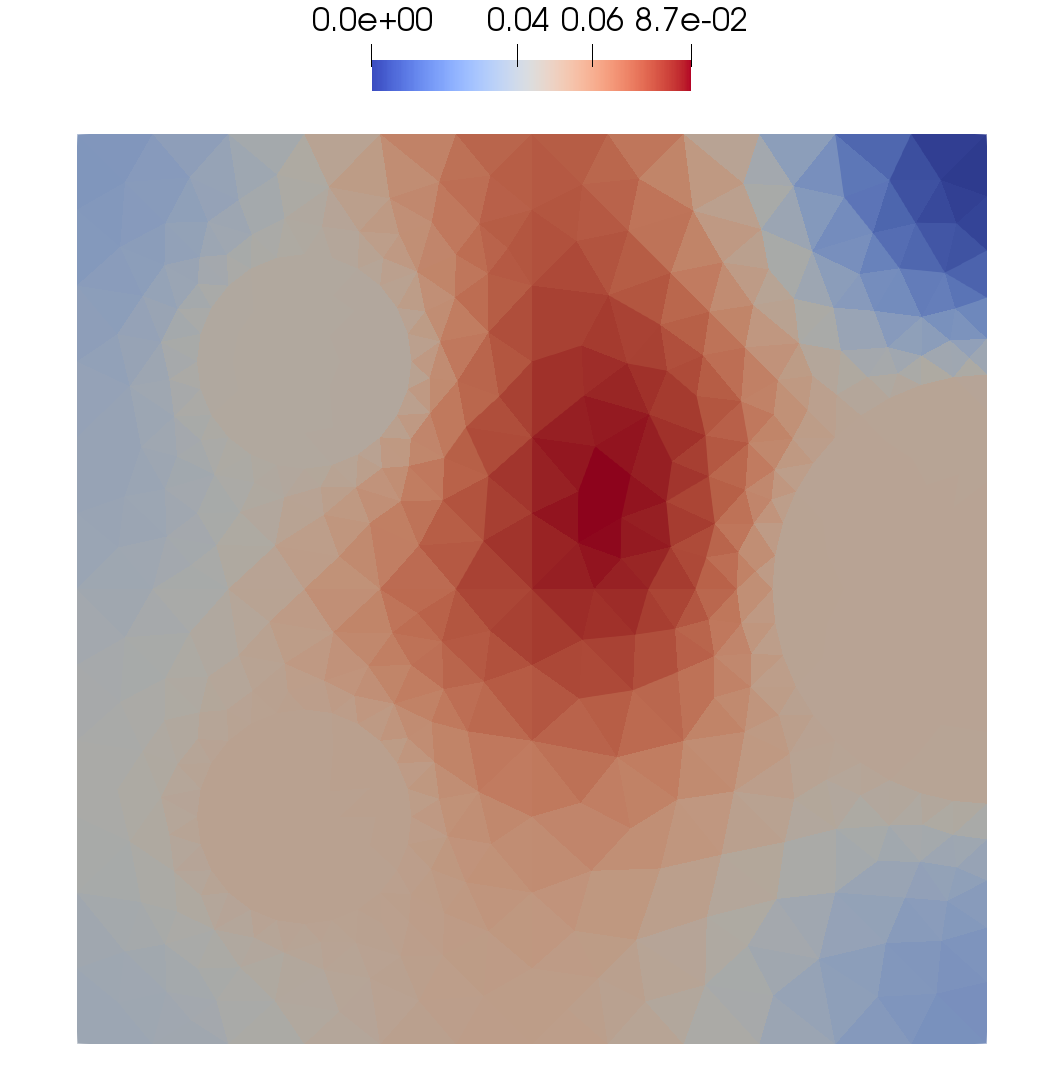}
\includegraphics[width=0.11\linewidth]{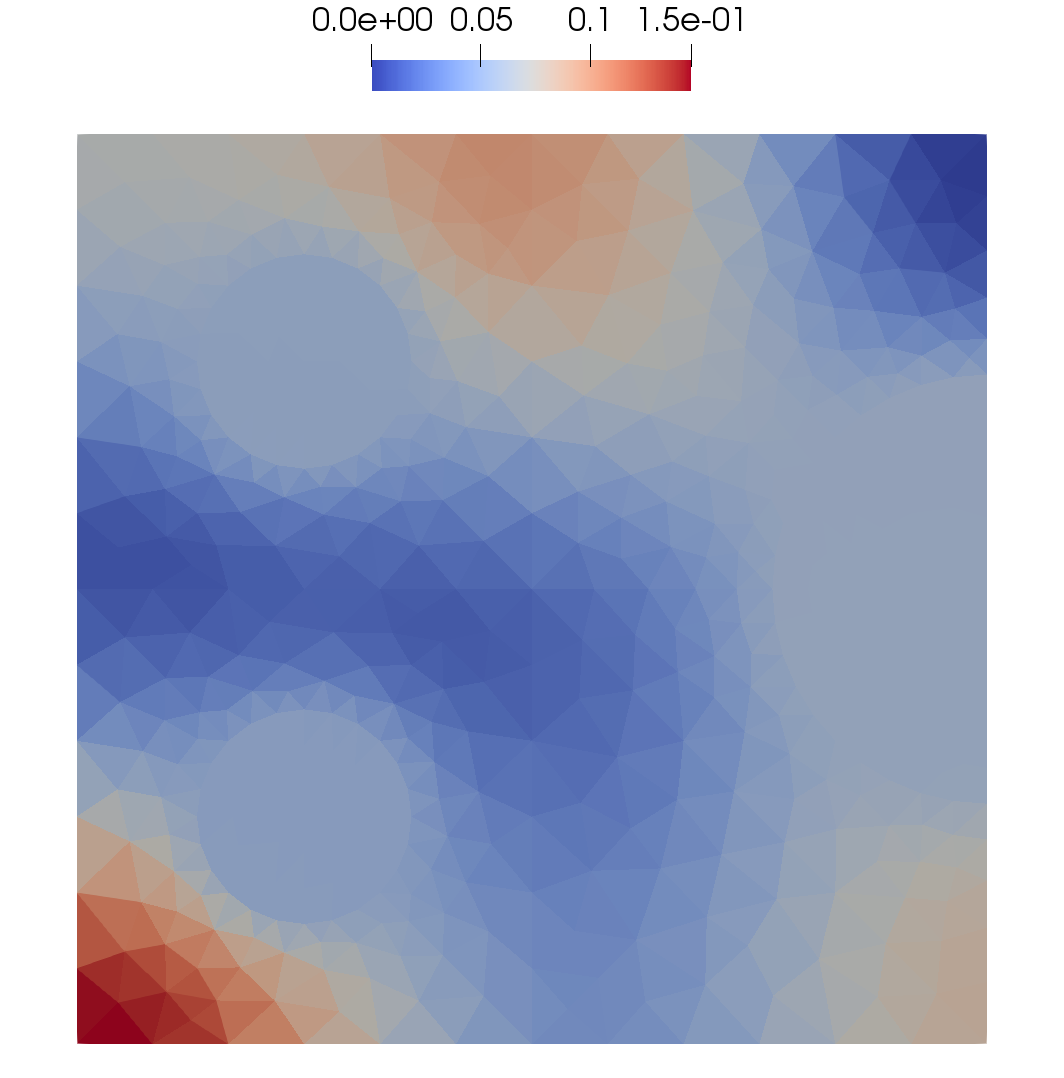}
\includegraphics[width=0.11\linewidth]{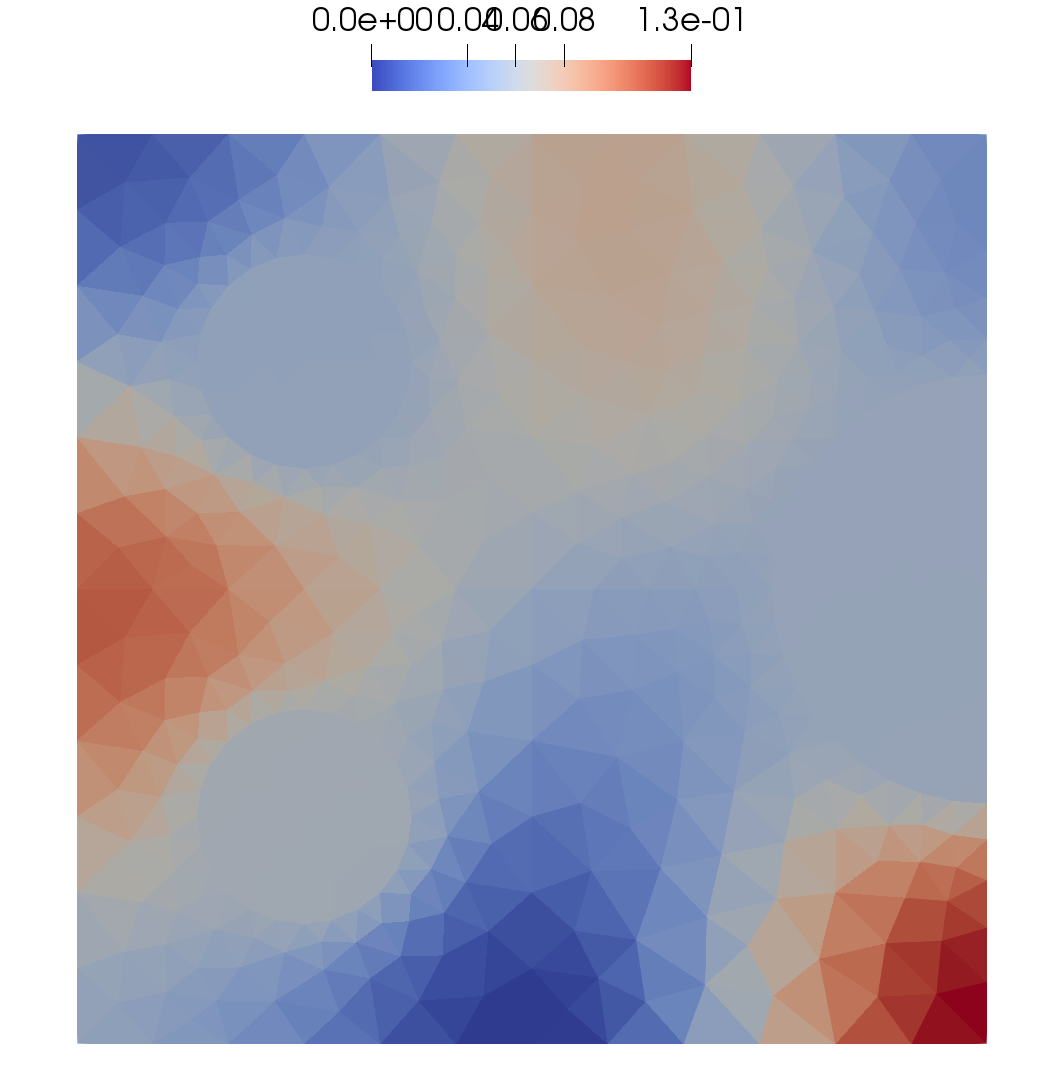}\\
\includegraphics[width=0.11\linewidth]{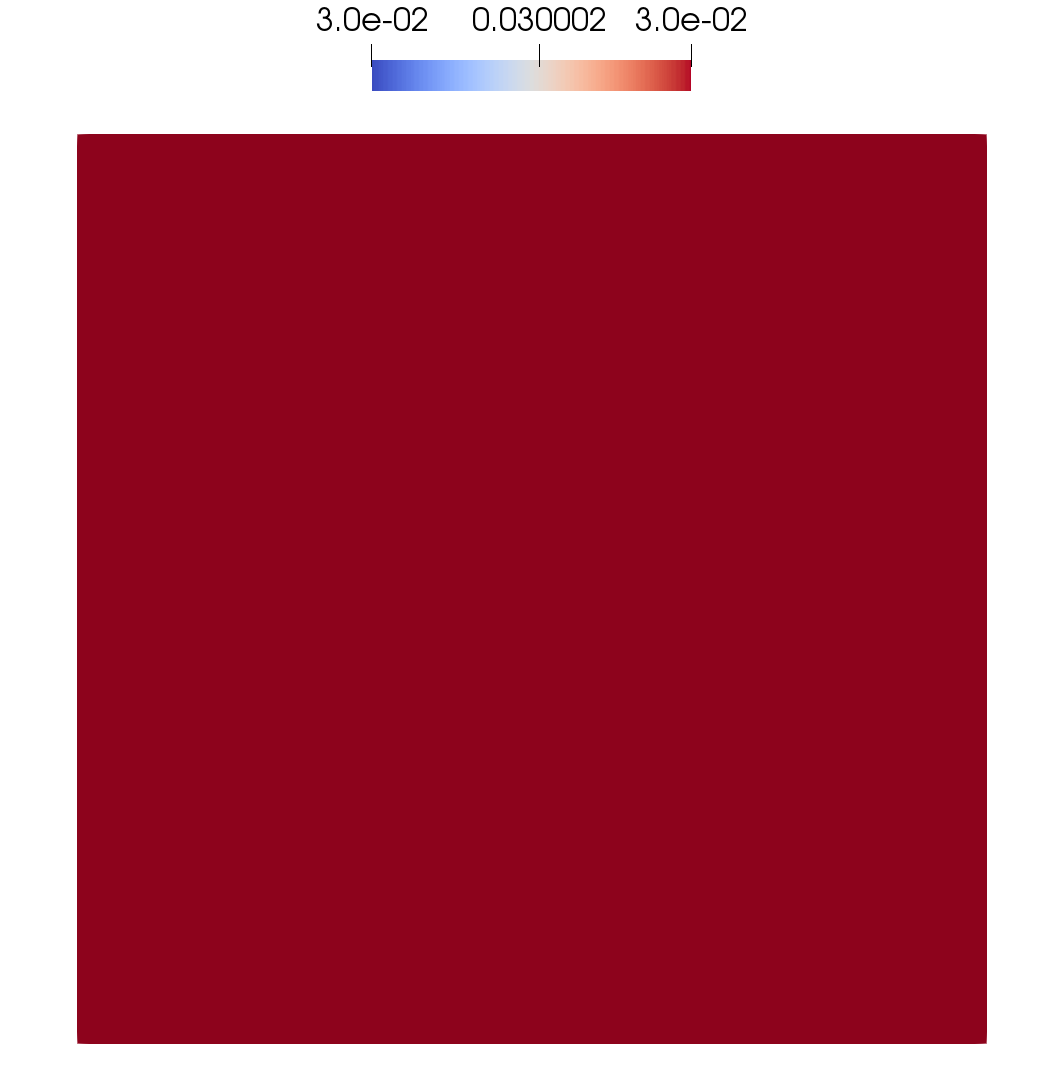}
\includegraphics[width=0.11\linewidth]{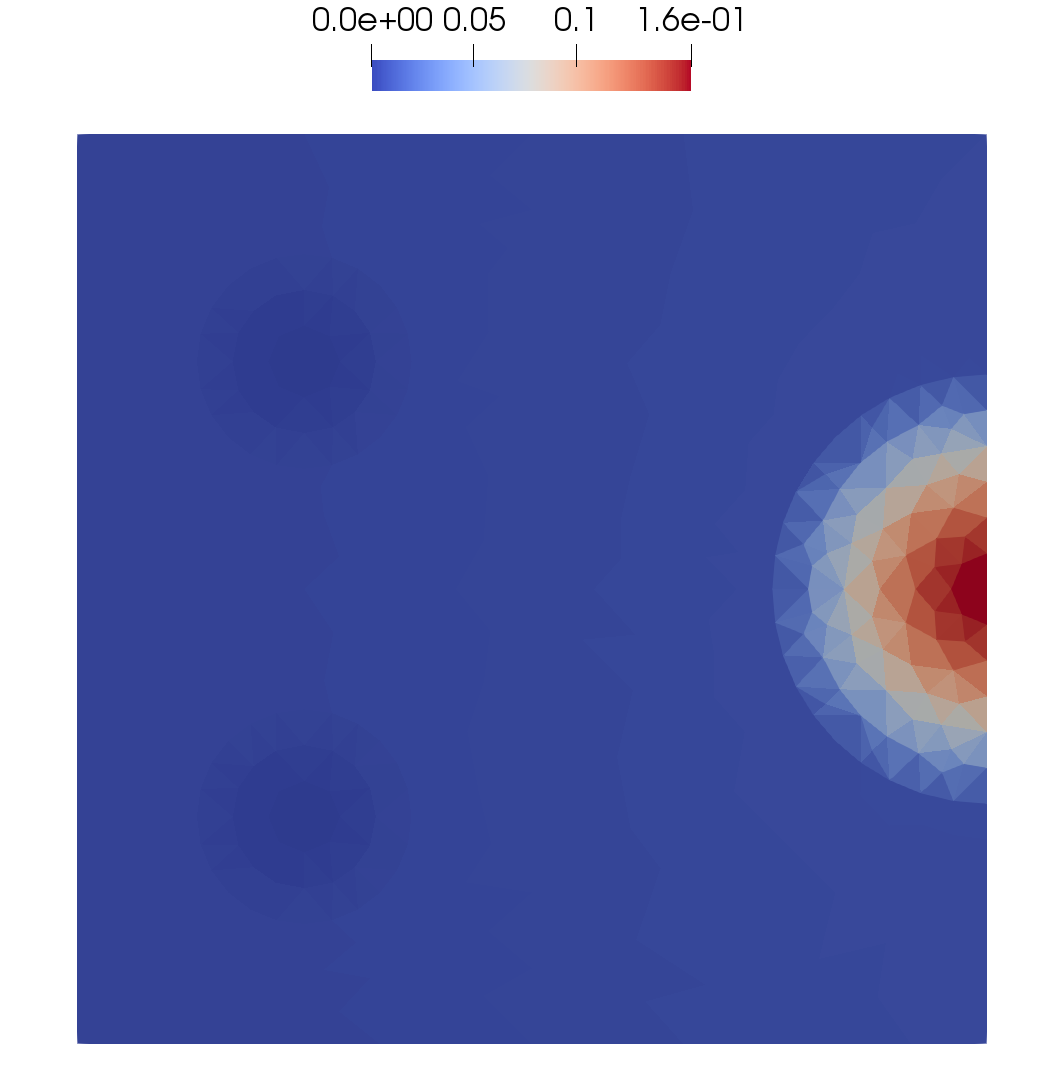}
\includegraphics[width=0.11\linewidth]{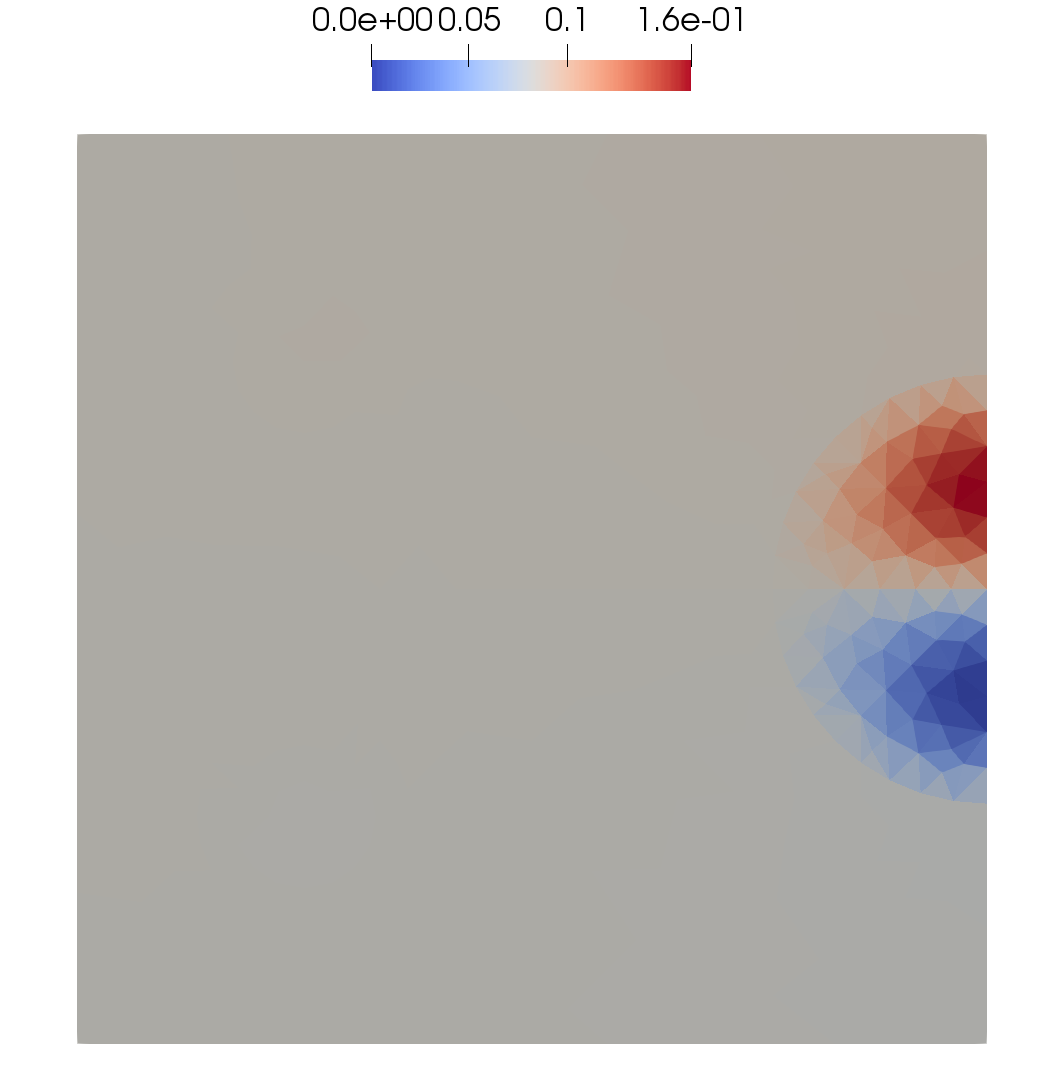}
\includegraphics[width=0.11\linewidth]{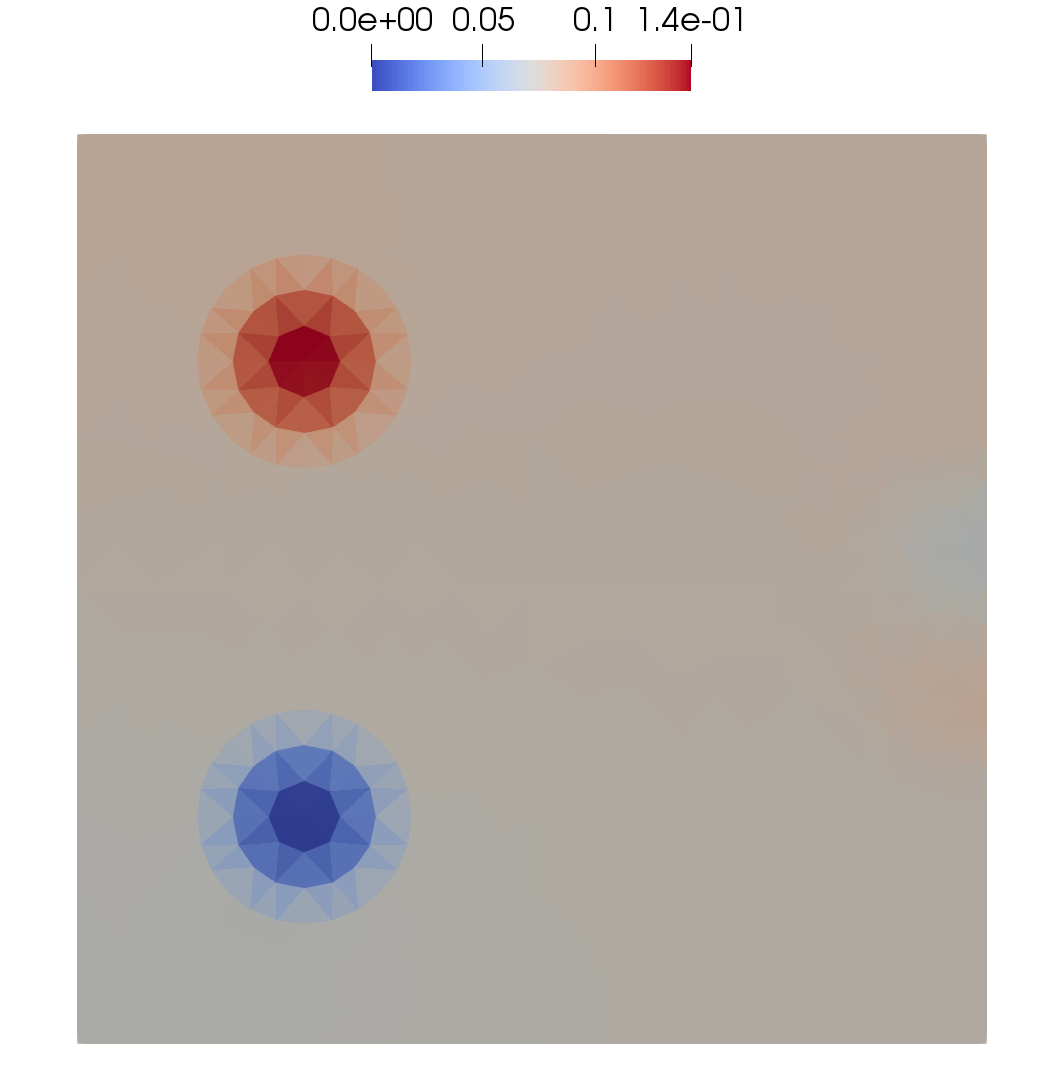}
\includegraphics[width=0.11\linewidth]{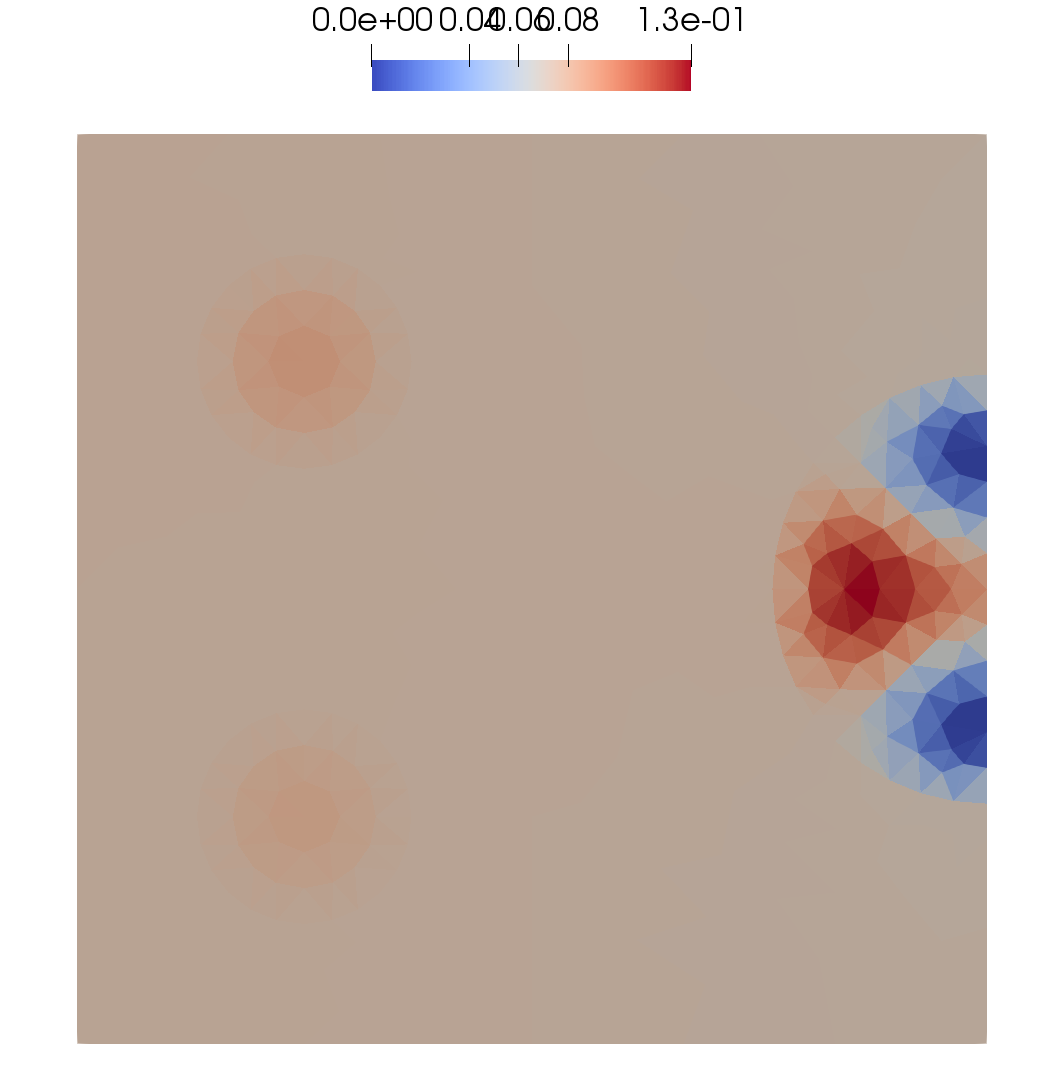}
\includegraphics[width=0.11\linewidth]{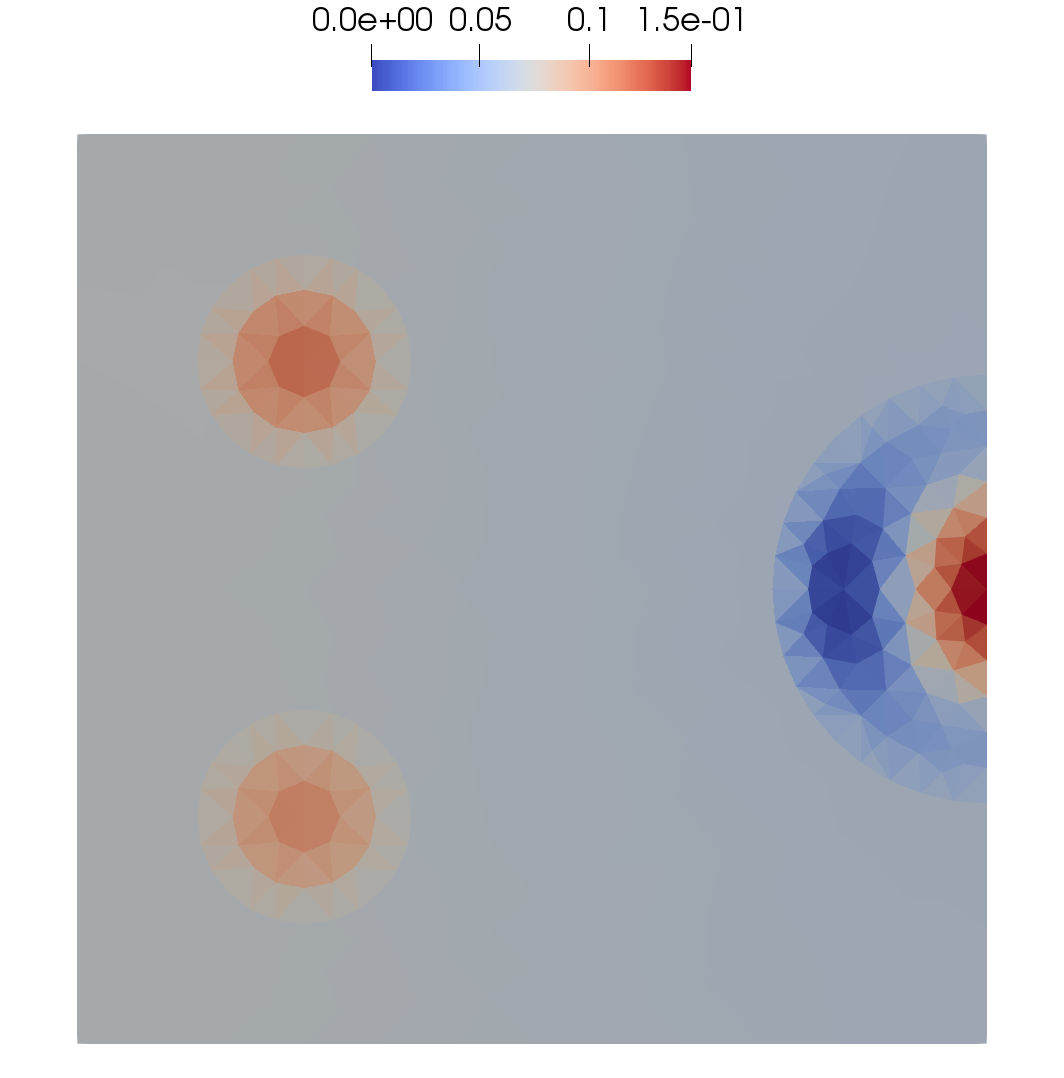}
\includegraphics[width=0.11\linewidth]{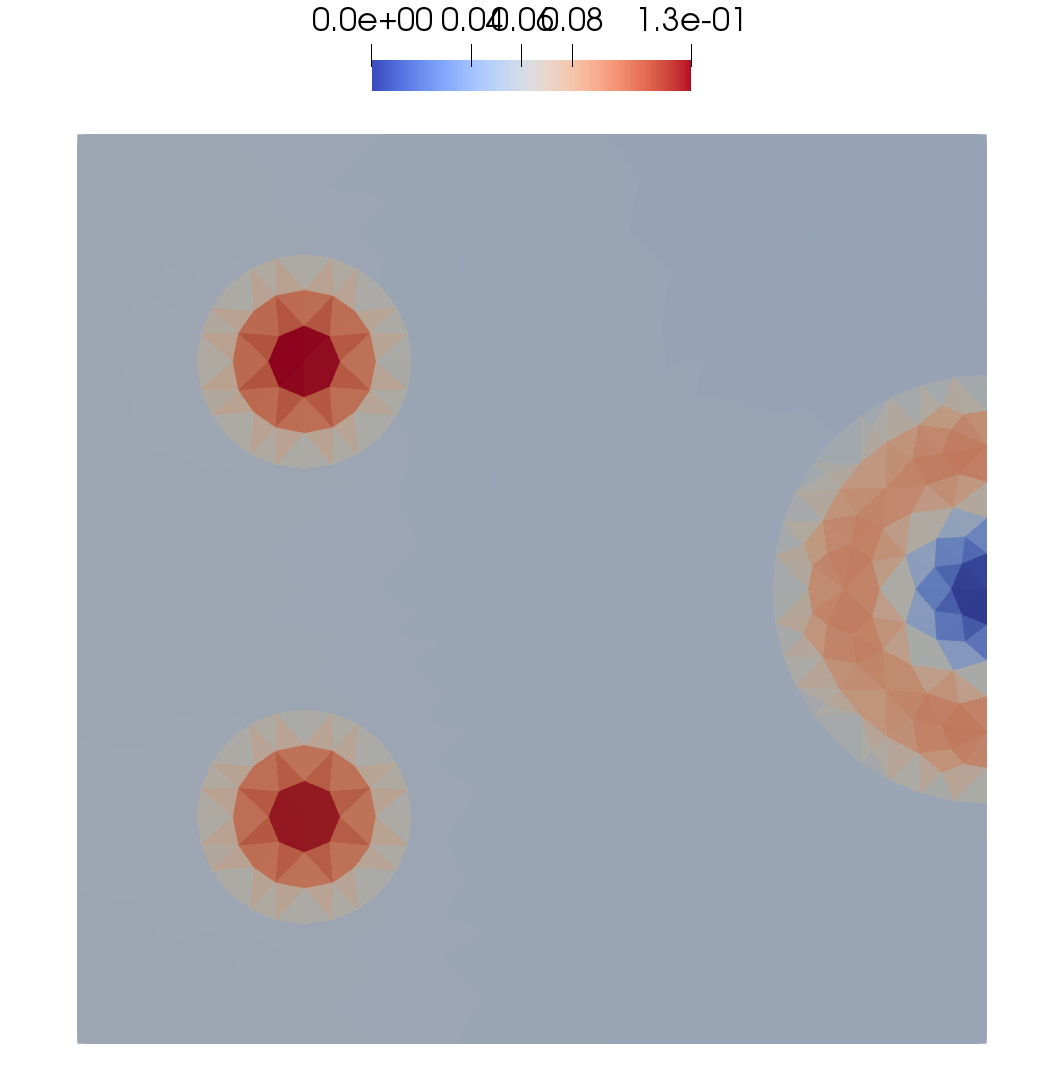}
\includegraphics[width=0.11\linewidth]{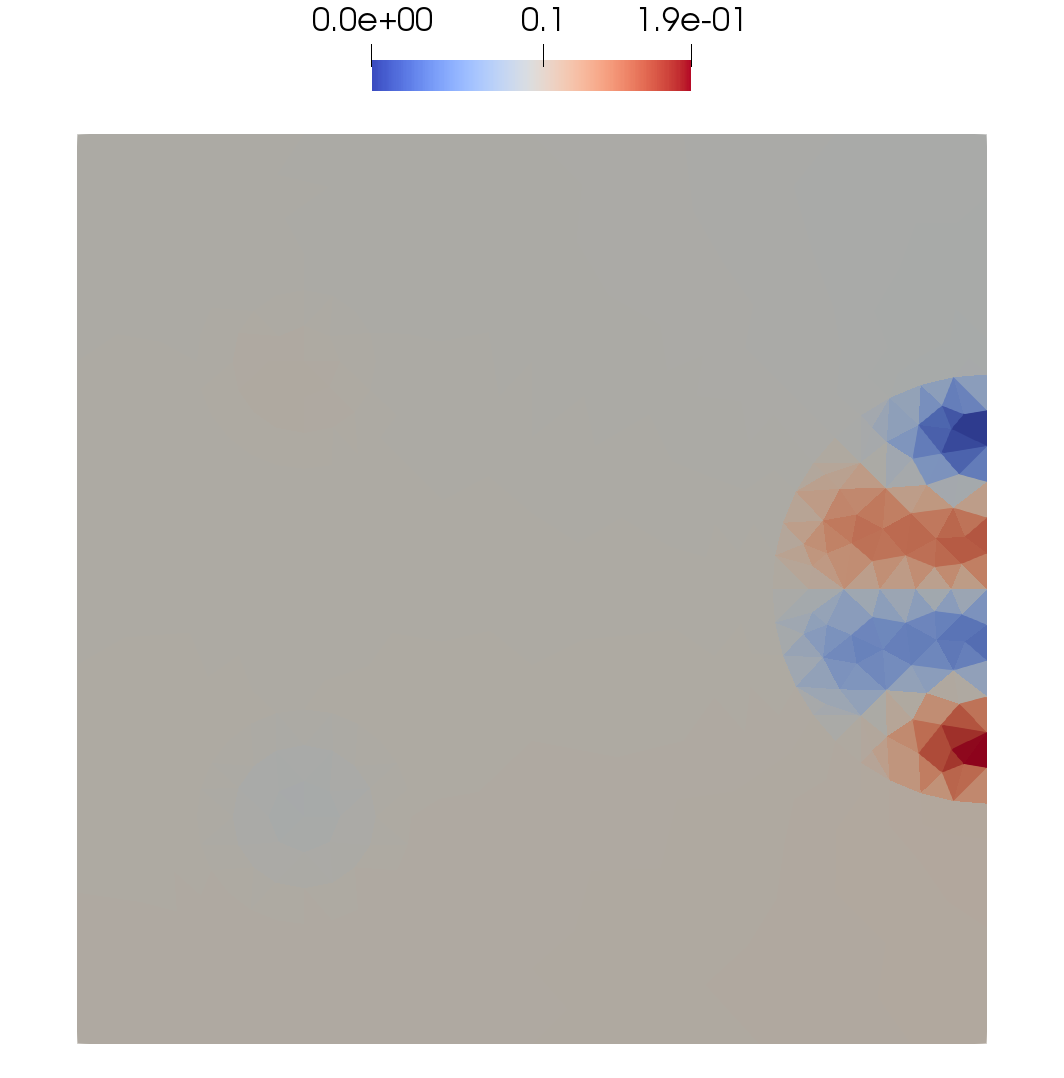}
\caption{Local domain and eight eigenvectors corresponding to the first eight smallest eigenvalues for $u^1$ (first row) and $u^2$ (second row)}
\label{basis}
\end{figure}

\subsection{Multiscale basis functions generation}

We solve the following eigenvalue problem in each local domain $\omega_i$
\begin{equation} 
\label{ms-sp}
A^{k, \omega_i} \psi^{k, \omega_i}_l = 
\lambda^l  \psi^{k, \omega_i}_l,
\end{equation}
where $A^{k, \omega_i}$ is the local diffusion matrix for component $k$
\[
A^{k, \omega_i} = \{a^k_{ij}\}, \quad 
a^k_{ij} = 
\left\{\begin{matrix}
\sum_j T^k_{ij} & i = j, \\ 
-T^k_{ij} & i \neq j
\end{matrix}\right. , 
\]
and $i,j = 1,...,N^{\omega_i}_f$, $N^{\omega_i}_f$ is the number of coarse grid cells in local domain $\omega_i$. 
To generate multiscale basis functions, we choose eigenvectors corresponding to the smallest $M_i$ eigenvalue and create a projection matrix for each component $k = 1,...,L$
\begin{equation} 
\label{ms-r}
P^k = \left[ 
\chi^1 \psi^{k, \omega_1}_1, \ldots, \chi^1 \psi^{k, \omega_1}_{M_1}
\ldots
\chi^{N_v} \psi^{k, \omega_{N_v}}_1, \ldots, \chi^{N_v} \psi^{k, \omega_{N_v}}_{M_{N_v}}
\right]^T, 
\end{equation}
where $\chi^i$ is the linear partition of unity functions and $N_v$ is the number of the local domains (number of coarse grid nodes). 
Illustration of the local domain and the relevant eigenvectors are presented in figure \ref{basis}. The first eight eigenvectors are shown in this figure representing inclusions that have larger diffusion coefficients (the first component) and the ones with smaller diffusion (the second component).

\subsection{Coarse grid system}

As mentioned before, we use the projection matrices for each component to obtain the reduced order model. 
On the coarse grid we have the following system in matrix form
\begin{equation}
M^k_H  \frac{ u^k _H- \check{u}^k_H}{\tau}  + 
A_H^k u^k_H=   R_H^k,
\end{equation}
with 
\[
A^k_H = P^k A^k (P^k)^T, 
\quad 
M^k_H = P^k M(P^k)^T, 
\quad 
R^k_H = P^k R^k
\]
where $u^k _H$ is the solution on the coarse grid. We use $u^k _{ms} = (P^k)^T u^k _H$ and compute the fine-scale solution.
The size of the resulting equation is $DOF_H = \sum_{i=1}^{N_v} M_i$, with $M_i$ denoting the number of local multiscale basis functions in $\omega_i$. In our numerical approach, we take the same number of basis functions in each local domain $\omega_i$ ($M_i = M$). Thus, we have $DOF_H = M \cdot N_v$.  Note that on the fine grid, the number of unknowns is $DOF_h = N$, where $N$ is the number of fine grid cells.
The convergence of the proposed method depends on the number of local basis functions and coarse grid size. We remark that a larger number of the multiscale basis functions may increase the size of system and consequently, take longer time to achieve accurate numerical solutions.

\section{Numerical results}

In this section, we numerically investigate the reaction-diffusion two-species competition model in heterogeneous media based on the proposed multiscale method. Setting $L=2$ in equation  \eqref{m}, we get
\[
\begin{split}
\frac{\partial u^1}{\partial t} 
- \nabla \cdot (\varepsilon^1(x) \nabla u^1) & 
=  r^1(x)  u^1 (1 - u^1)  - \alpha^{12}(x) u^1 u^2, 
\quad x \in \Omega, \\
\frac{\partial u^2}{\partial t} 
- \nabla \cdot (\varepsilon^2(x) \nabla u^2) & 
=  r^2(x)  u^2 (1 - u^2)  - \alpha^{21}(x) u^1 u^2, 
\quad x \in \Omega, 
\end{split}
\]
where $\Omega = \Omega_m \cup \Omega_c$ (see Figure \ref{domain}). In our test problems, we consider both small and regular diffusion scenarios. Thus, we set

\begin{itemize}
\item[] \textit{Small diffusion (a)}
\[
\varepsilon^1 = \left\{
\begin{matrix}
10^{-4}, & x \in \Omega_m \\ 
10^{-2}, & x \in \Omega_c
\end{matrix}
\right., 
\quad 
\varepsilon^2 = \left\{
\begin{matrix}
10^{-2}, & x \in \Omega_m \\ 
10^{-4}, & x \in \Omega_c
\end{matrix}
\right.
\]
\item[] \textit{Regular diffusion (b)}
\[
\varepsilon^1 = \left\{
\begin{matrix}
10^{-3}, & x \in \Omega_m \\ 
10^{-1}, & x \in \Omega_c
\end{matrix}
\right., 
\quad 
\varepsilon^2 = \left\{
\begin{matrix}
10^{-1}, & x \in \Omega_m \\ 
10^{-3}, & x \in \Omega_c
\end{matrix}
\right.
\]
\end{itemize}
Note that the regular diffusion is related to the same scale as of the birth rate, and the smaller rate means having 10 times smaller diffusion coefficients. We specifically consider the following two test problems:
\begin{itemize}
\item \textbf{Test 1}
\[
r^1 = \left\{
\begin{matrix}
0.15, & x \in \Omega_m \\ 
0.1 & x \in \Omega_c
\end{matrix} \right., 
\quad
r^2 = \left\{
\begin{matrix}
0.1, & x \in \Omega_m \\ 
0.15 & x \in \Omega_c
\end{matrix} \right., 
\]\[
\alpha^{12} = \left\{
\begin{matrix}
0.055, & x \in \Omega_m \\ 
0.05 & x \in \Omega_c
\end{matrix} \right., 
\quad 
\alpha^{21} = \left\{
\begin{matrix}
0.05, & x \in \Omega_m \\ 
0.055 & x \in \Omega_c
\end{matrix} \right., 
\]
and $t_{max} = 50$ and 100 time iterations. 
\begin{itemize}
\item[] \textit{Test 1a} (small diffusion).
\item[] \textit{Test 1b} (regular diffusion).
\end{itemize}
\item \textbf{Test 2}
\[
r^1 = \left\{
\begin{matrix}
0.15, & x \in \Omega_m \\ 
0.1 & x \in \Omega_c
\end{matrix} \right., 
\quad
r^2 = \left\{
\begin{matrix}
0.1, & x \in \Omega_m \\ 
0.15 & x \in \Omega_c
\end{matrix} \right., 
\]\[
\alpha^{12} = \left\{
\begin{matrix}
0.15, & x \in \Omega_m \\ 
0.01 & x \in \Omega_c
\end{matrix} \right., 
\quad 
\alpha^{21} = \left\{
\begin{matrix}
0.01, & x \in \Omega_m \\ 
0.075 & x \in \Omega_c
\end{matrix} \right., 
\]
and $t_{max} = 150$ and 100 time iterations. 
\begin{itemize}
\item[] \textit{Test 2a} (small diffusion).
\item[] \textit{Test 2b} (regular diffusion).
\end{itemize}
\end{itemize}
The fine grid contains 69948 cells and we have $DOF_h  = 69948$ for each species/component in the uncoupled  \textit{SI} scheme.  
But in the coupled \textit{FI} scheme,  the size of system is $DOF_h  = L \cdot 69948$, where $L=2$ (two-species). 
In both test problems, we perform simulation for 100 time steps with $u^1_0 = u^2_0 = 0.5$ as the initial conditions.

For the construction of the computational geometry with inclusions and grid generation, we use Gmsh \cite{geuzaine2009gmsh}. 
We adopt PETSc library \cite{balay2019petsc} for numerical implementation.  In order to work with unstructured grid elements, we use FEniCS library \cite{logg2012automated, logg2012efficient} and Paraview \cite{ahrens2005paraview} for the visualization of results.

\begin{figure}[h!]
\centering
\begin{subfigure}{0.48\textwidth}
\includegraphics[width=1\linewidth]{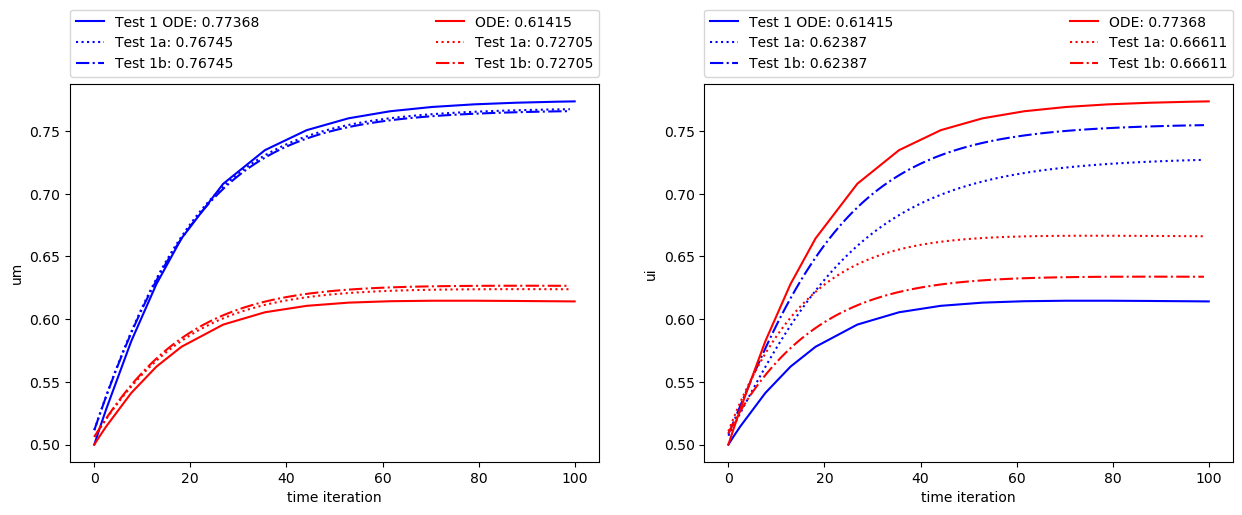}
\caption{Test 1. $\bar{u}^k_m$ (left) and $\bar{u}^k_c$ (right)}
\end{subfigure}
\,\,\,\,
\begin{subfigure}{0.48\textwidth}
\includegraphics[width=1\linewidth]{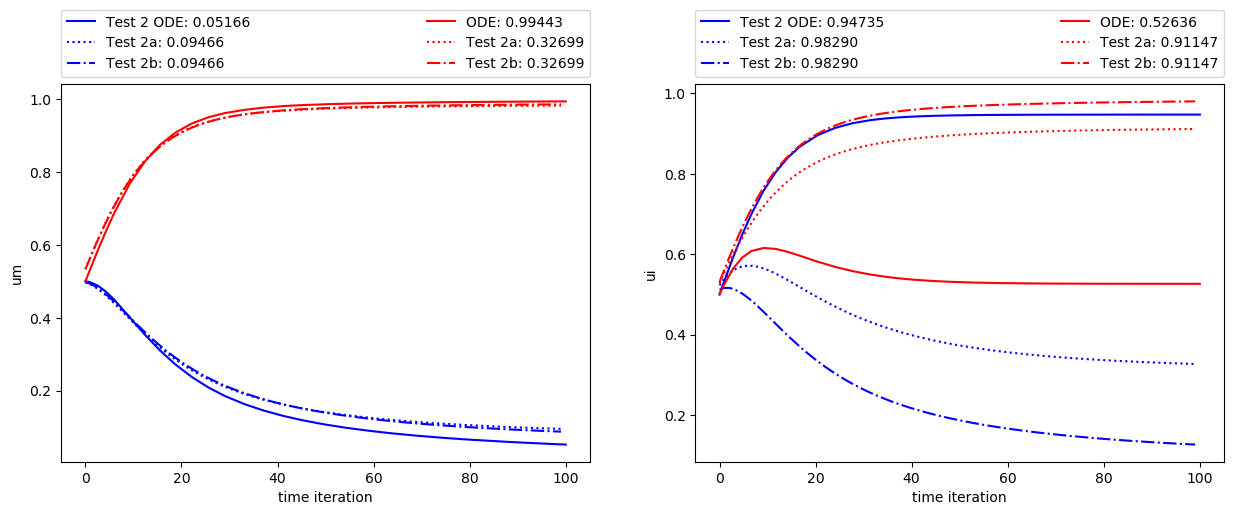}
\caption{Test 2. $\bar{u}^k_m$ (left) and $\bar{u}^k_c$ (right)}
\end{subfigure}
\caption{Effect of diffusion. Dynamic average solution $\bar{u}^1$ (red color) and $\bar{u}^2$ (blue color) with no diffusion (ODE) and solutions on fine grid for small and regular diffusion. 
Left: average in background domain $\Omega_m$, $\bar{u}^k_m$.  
Right: average in subdomain $\Omega_c$, $\bar{u}^k_c$} 
\label{sol-d}
\end{figure}

\begin{figure}[h!]
\centering
\begin{subfigure}{0.22\textwidth}
\includegraphics[width=1\linewidth]{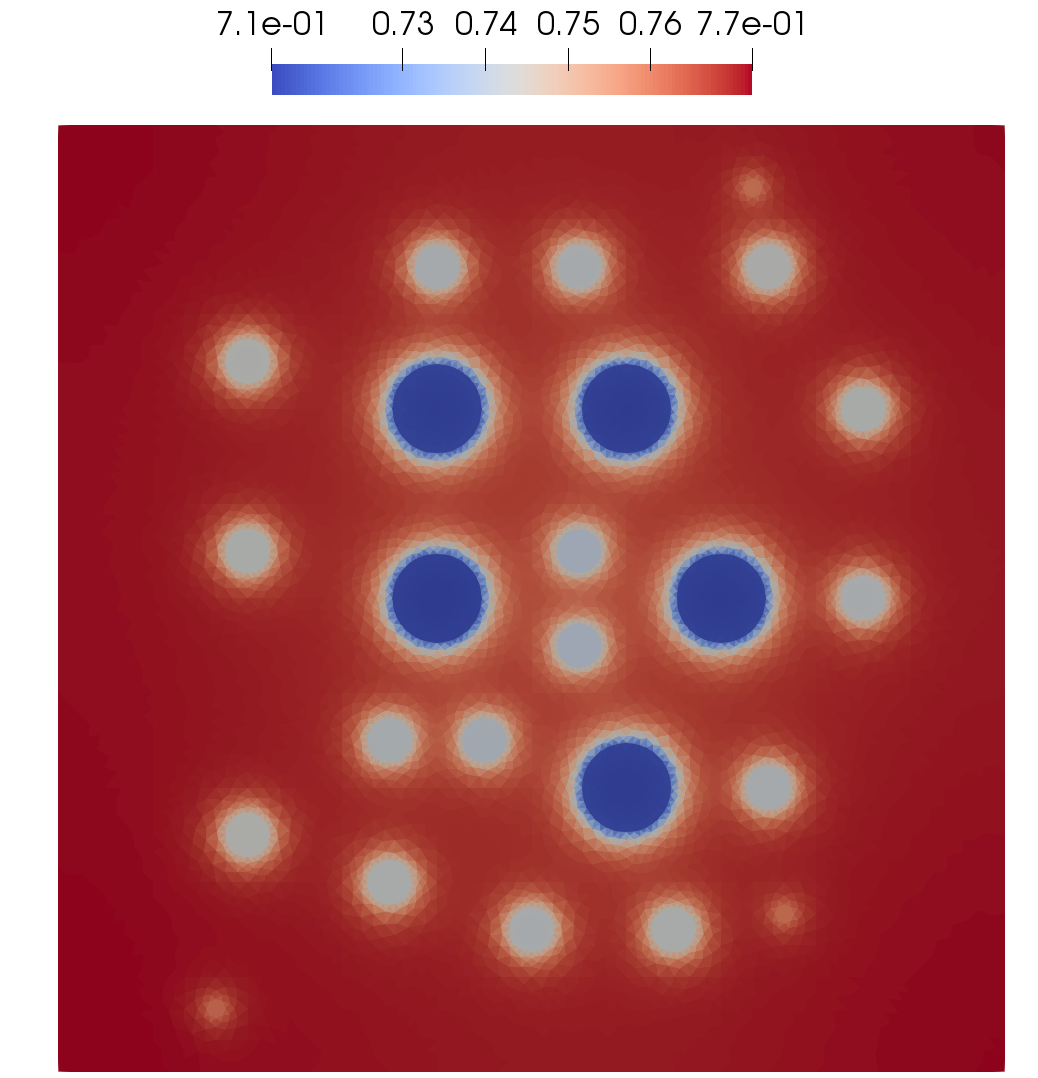}
\includegraphics[width=1\linewidth]{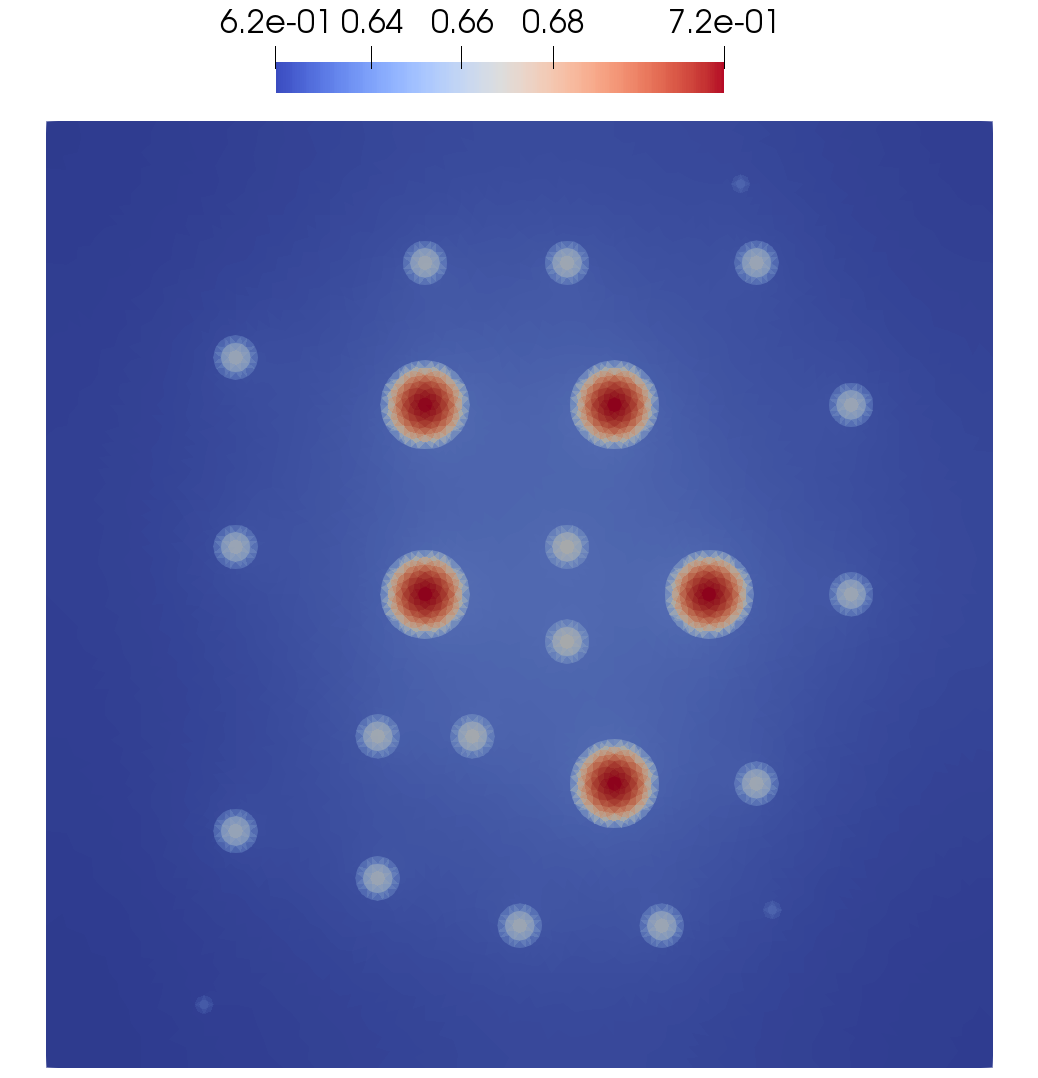}
\caption{Test 1a}
\end{subfigure}
\,\,\,\,
\begin{subfigure}{0.22\textwidth}
\includegraphics[width=1\linewidth]{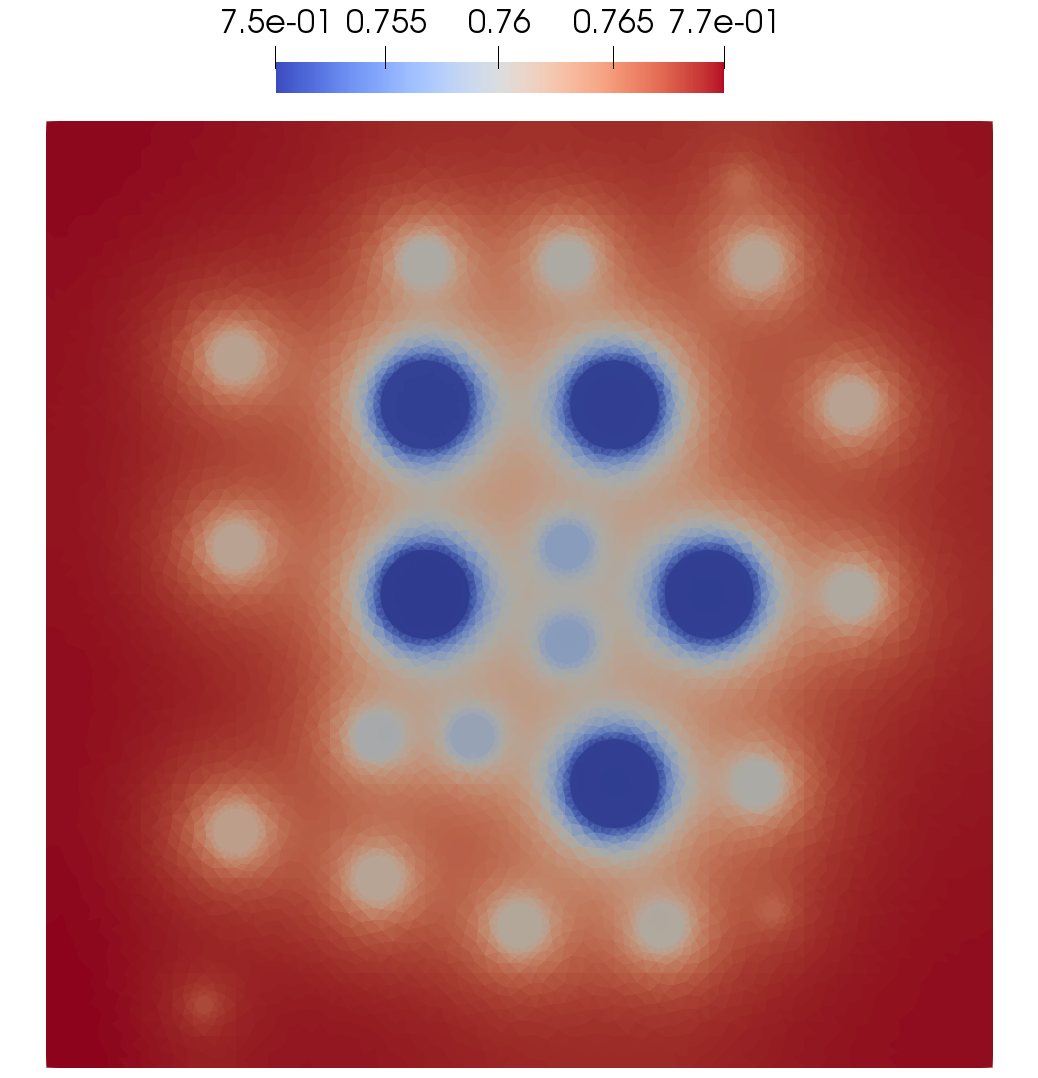}
\includegraphics[width=1\linewidth]{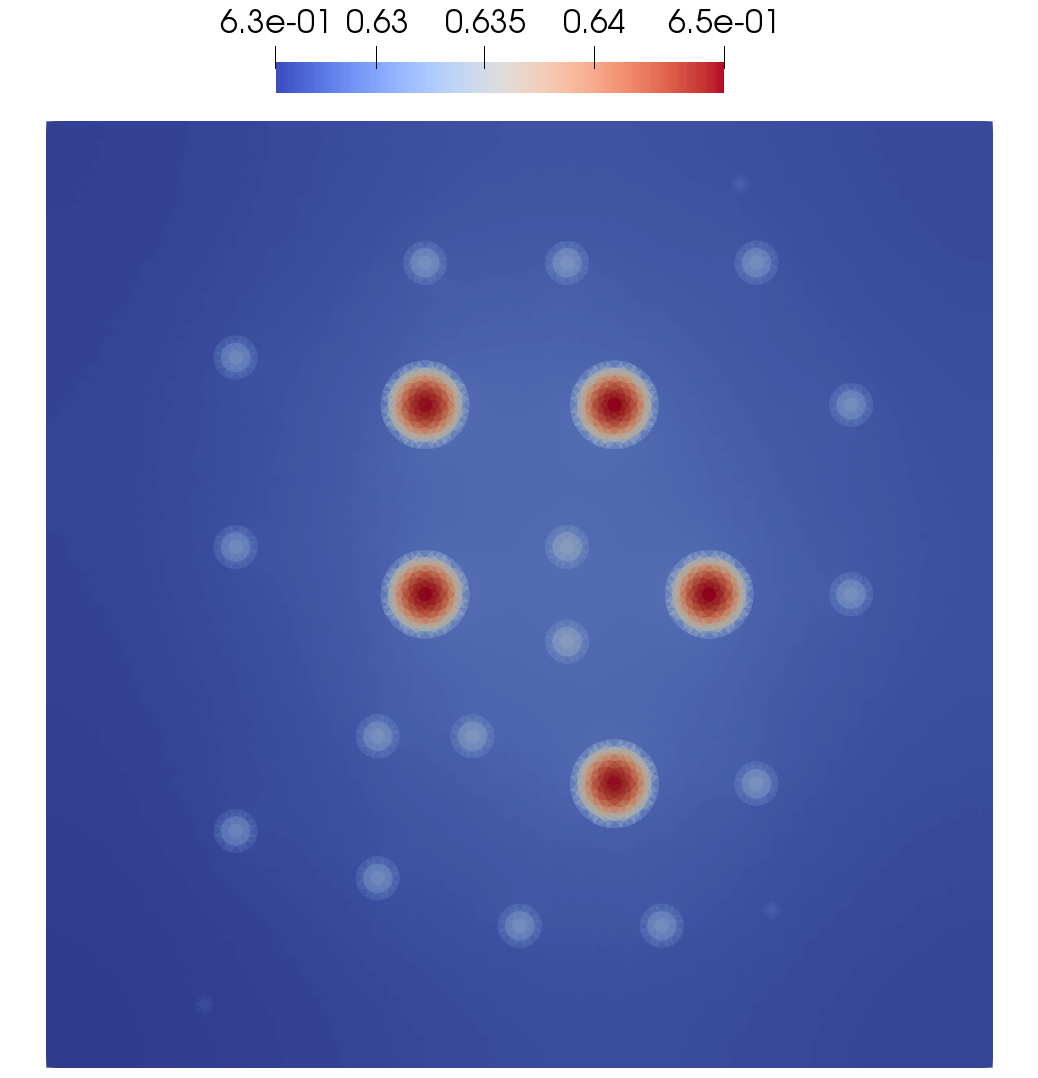}
\caption{Test 1b}
\end{subfigure}
\,\,\,\,
\begin{subfigure}{0.22\textwidth}
\includegraphics[width=1\linewidth]{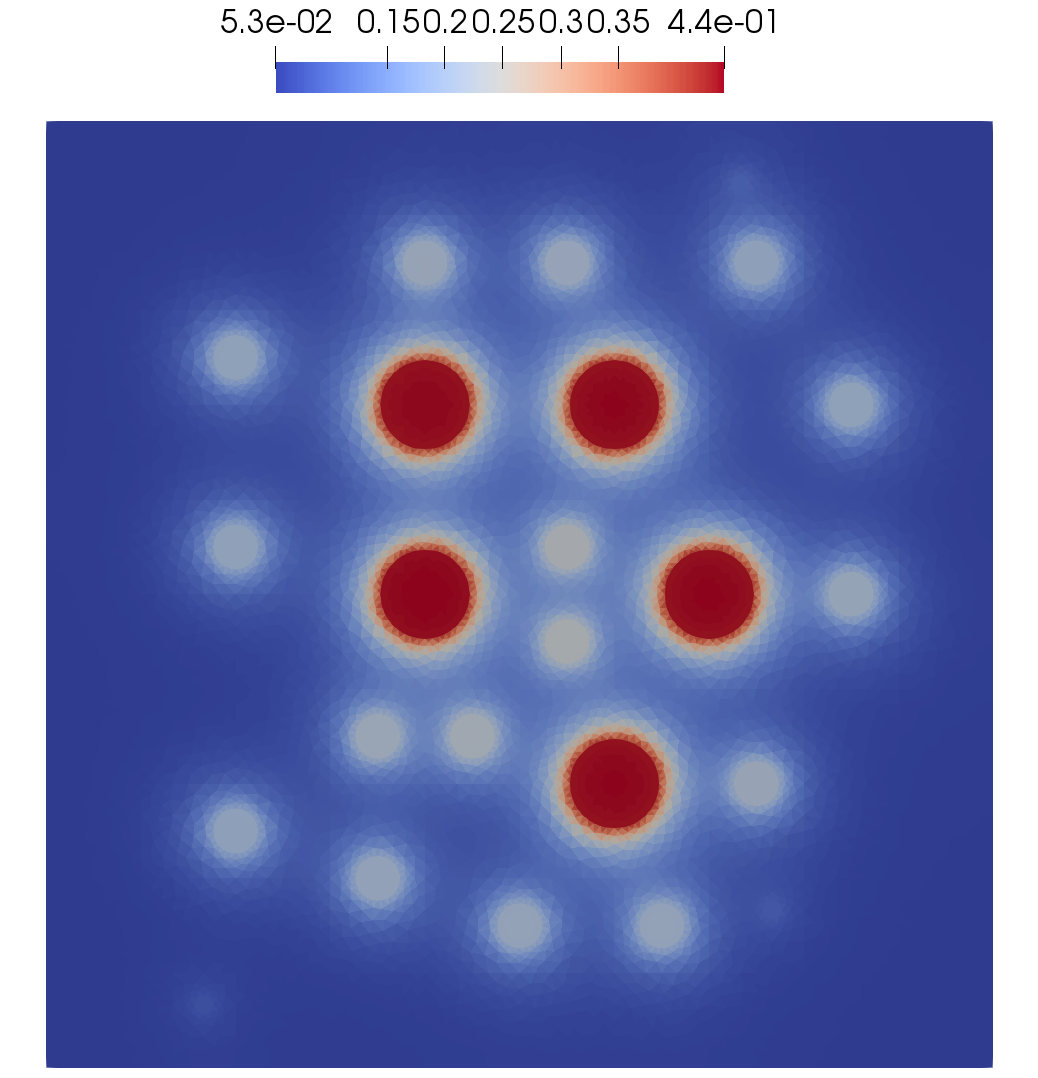}
\includegraphics[width=1\linewidth]{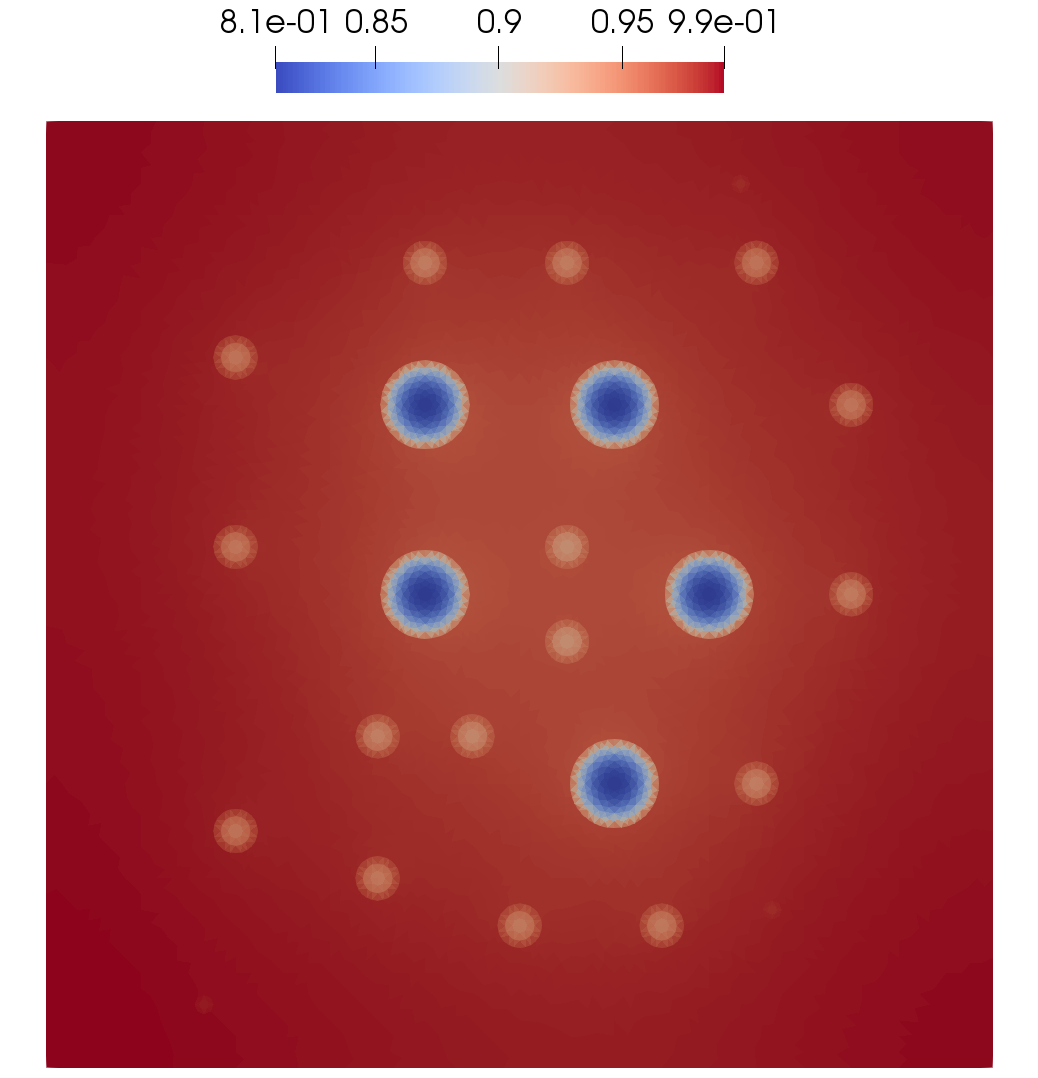}
\caption{Test 2a}
\end{subfigure}
\,\,\,\,
\begin{subfigure}{0.22\textwidth}
\includegraphics[width=1\linewidth]{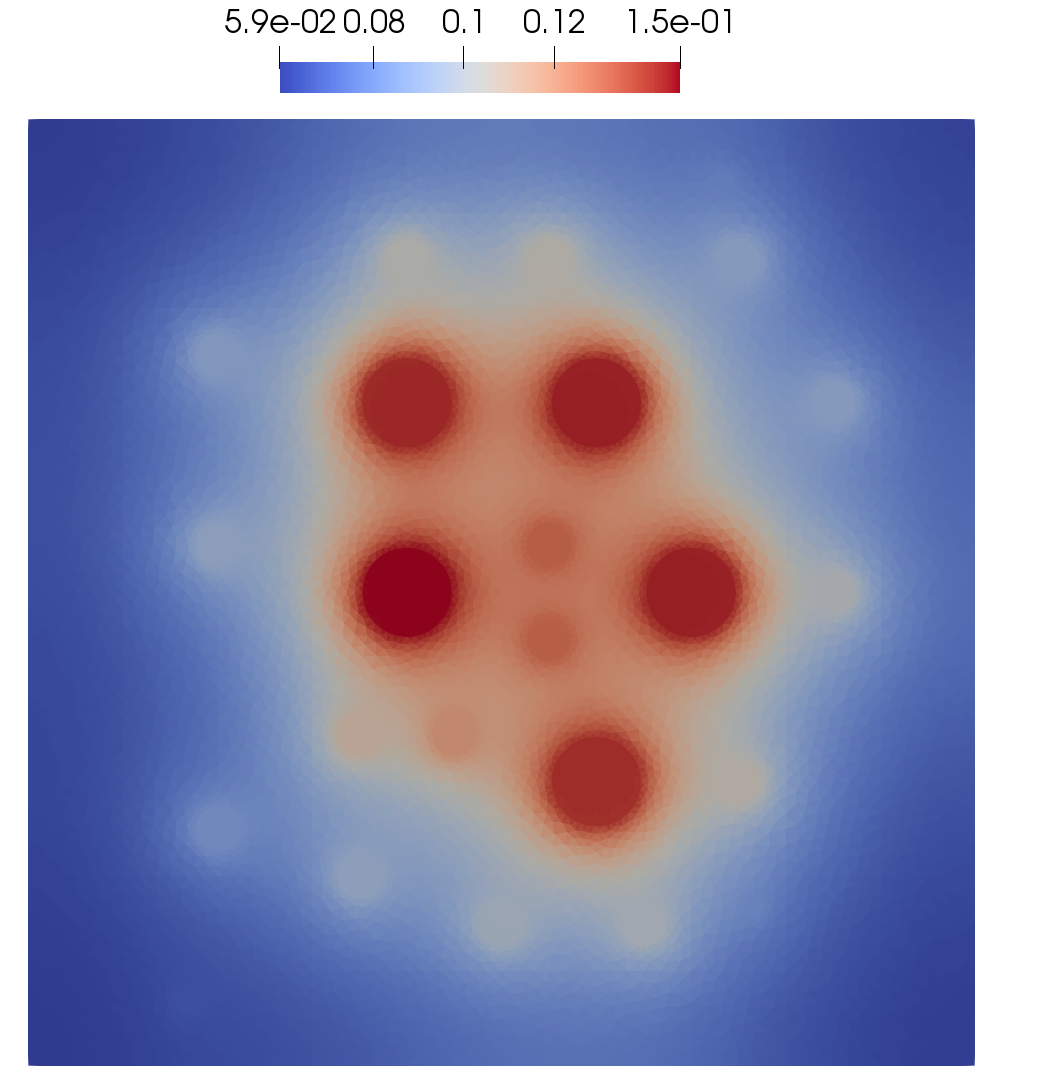}
\includegraphics[width=1\linewidth]{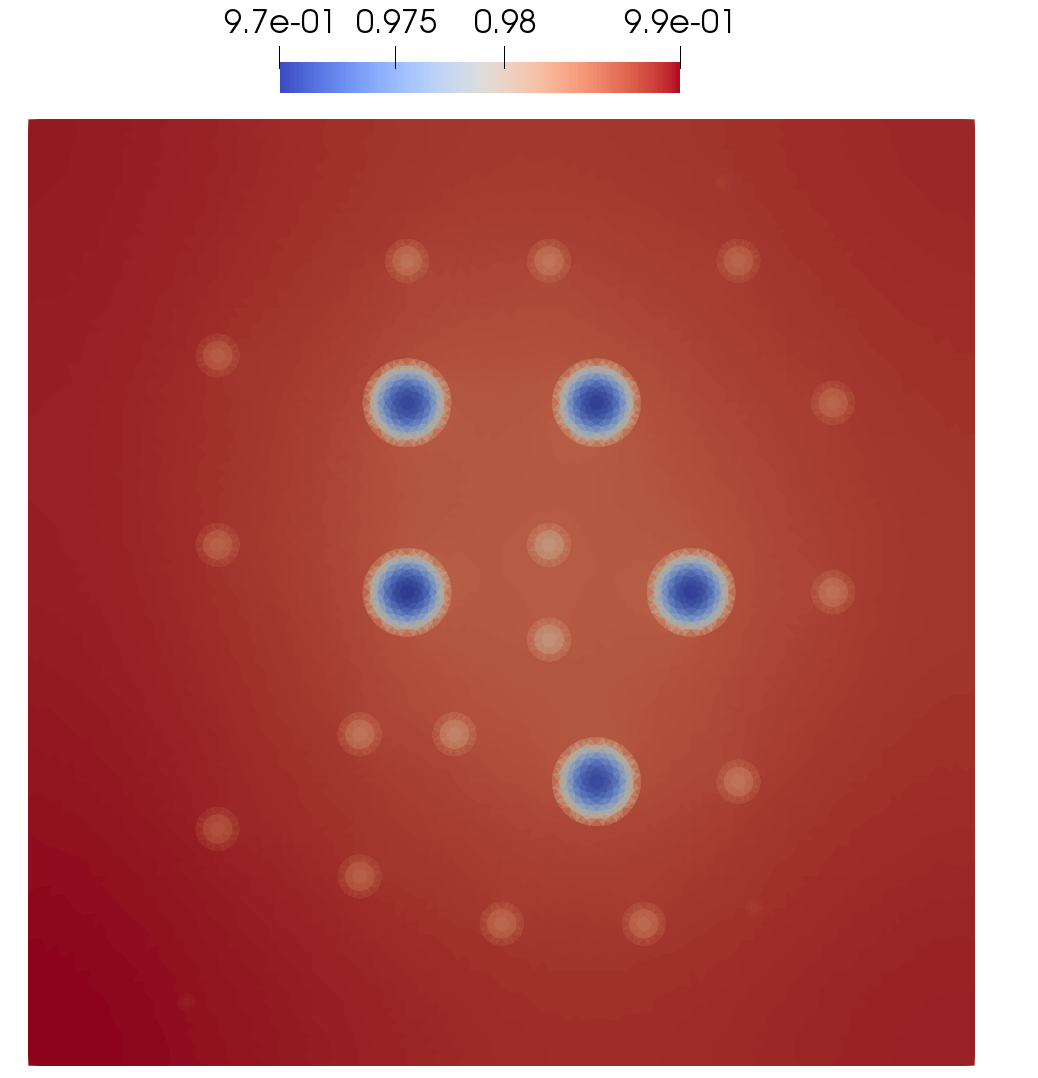}
\caption{Test 2b}
\end{subfigure}
\caption{Fine grid solution at the final time for Test 1a (small diffusion), Test 1b (regular diffusion), Test 2a (small diffusion) and Test 2b (regular diffusion).  
Firsr row: $u^1$. Second row: $u^2$ }
\label{sol-u}
\end{figure}

The solution average values of the solutions for Test 1 and Test 2 are plotted in Figure \ref{sol-d}. The plots display the cases with no diffusion (Test 1 (ODE) and 2 (ODE)), small diffusion (Test 1a and 2a) and regular diffusion (Test 1b and 2b). The average solution is computed using the formulas given by
\[
\bar{u}^k_m = \frac{1}{|\Omega_m|}\int_{\Omega_m} u^k \ dx, 
\quad 
\bar{u}^k_c = \frac{1}{|\Omega_c|}\int_{\Omega_c} u^k \ dx, 
\quad 
k = 1,2,
\]
where $|\Omega_m|$ and $|\Omega_c|$ are the volume of the domains $\Omega_m$ and $\Omega_c$, respectively. The displayed values show the solutions at the final time. 
To illustrate the impact of diffusion on the test problems solutions, we utilize the solution on the fine grid with no diffusion as well as with small and regular diffusion.  Note that we employ LSODE solver to solve the ODE system. 
We notice a greater influence of diffusion on the solution in inclusion subdomains. However effect of the diffusion on the solution average in the background domain is smaller. The plots in Figure \ref{sol-d} demonstrates that both species survive and coexist in the background media for Test 1 run, but the population of the second species dominates the first. In Test 2, we see the opposite trend, that is, the  population of the second species is very small and first species prevail in the background media. The influence of the diffusion in inclusions is very large and thus can change control the choice of dominance of the species. For instance, the first species dominates over second in Test 1 in the case of no diffusion, but with diffusion we see that second species is in command over the first. A similar behavior is observed in the inclusions for Test 2. Hence, one can say that there is a strong competition between the species with diffusion. Alternatively, diffusion can dictate the dominance of one species over the other.

The fine grid solutions at the final time for Test 1a (small diffusion), Test 1b (regular diffusion), Test 2a (small diffusion) and Test 2b (regular diffusion) are portrayed in Figure \ref{sol-u}.  The population distributions for the first $u^1$ and the second $u^2$ species are shown in row 1 and 2, respectively.  The spatial influence due to the diffusion on the numerical solutions is apparent in this figure.
 
Next, we turn our attention on the numerical comparison of time approximation techniques on the fine grid to explain the computational effectiveness of the uncoupled scheme.

\subsection{Solutions using the coupled and uncoupled schemes for time approximation}

We solve problems on the fine grid using coupled fully implicit (\textit{FI}) scheme and uncoupled semi-implicit time approximation scheme (\textit{SI}).

\begin{figure}[h!]
\centering
\begin{subfigure}{0.48\textwidth}
\includegraphics[width=1\linewidth]{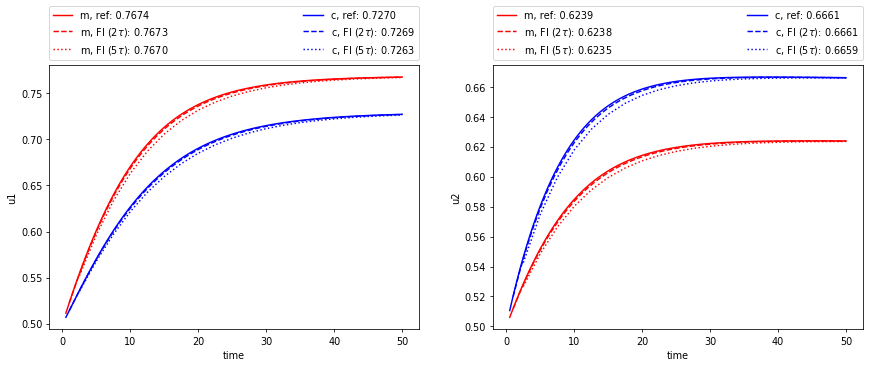}
\includegraphics[width=1\linewidth]{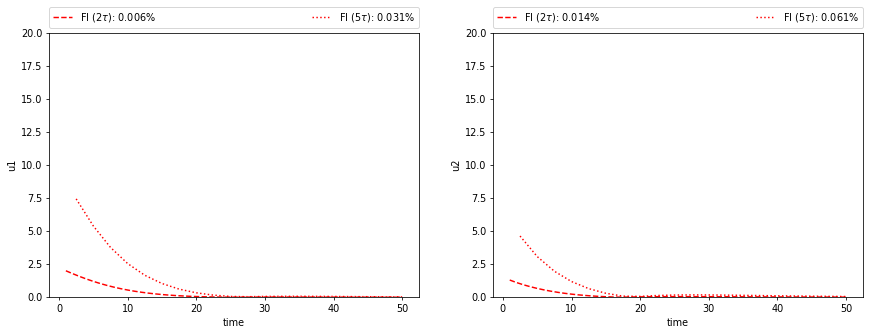}
\caption{Test 1a \textit{FI}}
\end{subfigure}
\,\,\,\,
\begin{subfigure}{0.48\textwidth}
\includegraphics[width=1\linewidth]{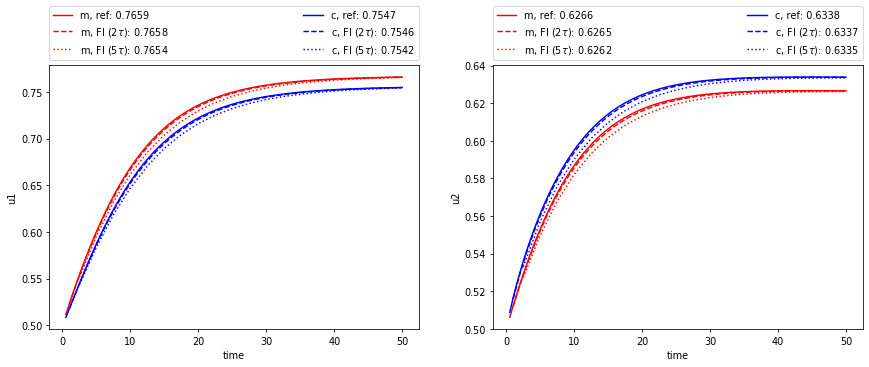}
\includegraphics[width=1\linewidth]{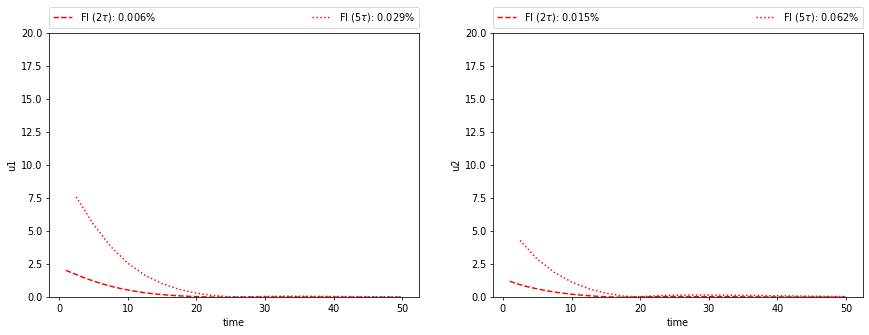}
\caption{Test 1b \textit{FI}}
\end{subfigure}\\
\begin{subfigure}{0.48\textwidth}
\includegraphics[width=1\linewidth]{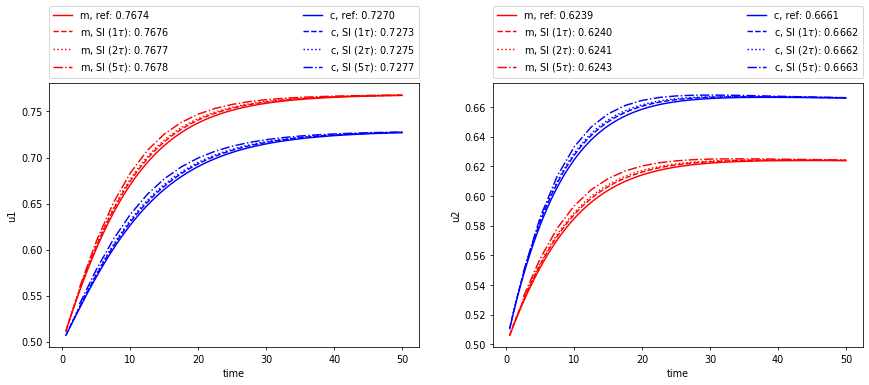}
\includegraphics[width=1\linewidth]{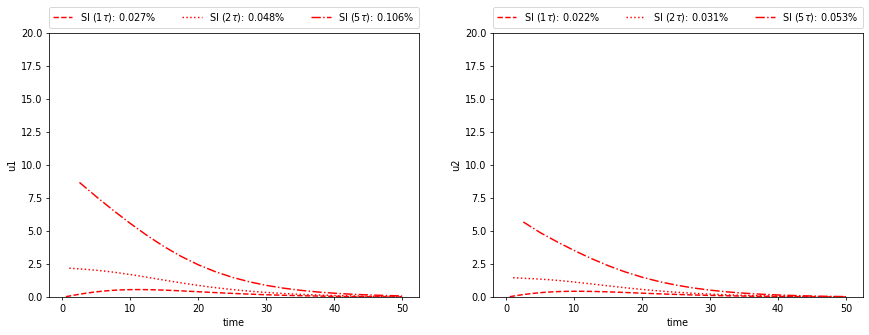}
\caption{Test 1a \textit{SI}}
\end{subfigure}
\,\,\,\,
\begin{subfigure}{0.48\textwidth}
\includegraphics[width=1\linewidth]{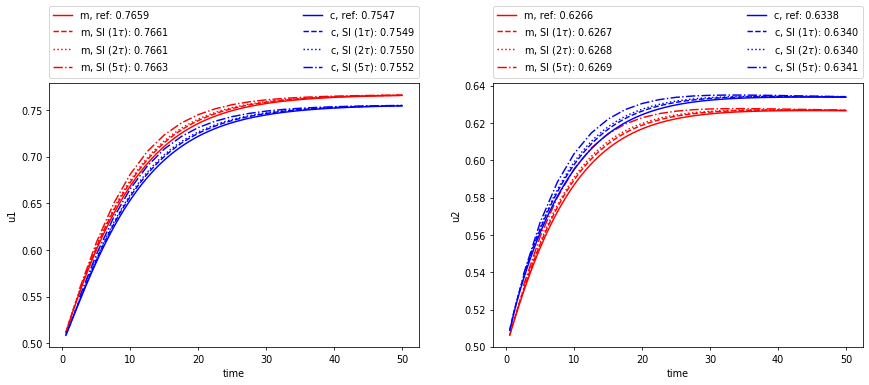}
\includegraphics[width=1\linewidth]{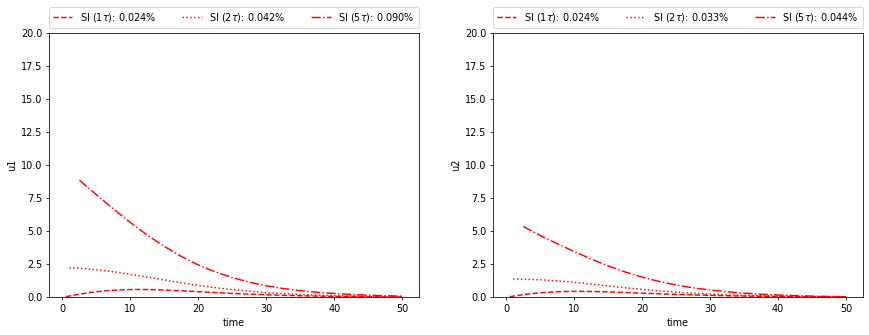}
\caption{Test 1b \textit{SI}}
\end{subfigure}
\caption{Test 1. Dynamic average solutions (first row) and errors (second row) for different time step size, $\tau = 1/M$, $2 \tau = 2/M$, and $10 \tau = 10/M$ ($M=100$) } 
\label{sol-fi-1}
\end{figure}

\begin{table}[h!]
\centering
\begin{tabular}{|c|cccc|cccc|}
\hline
&
 \multicolumn{4}{|c|}{Test 1a}&
 \multicolumn{4}{|c|}{Test 1b}\\ 
  & $N$ 
& $e_1$ (\%)  & $e_2$( \%) & time(sec) 
& $N$ 
& $e_1$ (\%)  & $e_2$ (\%) & time(sec)  \\
\hline
\multicolumn{9}{|c|}{\textit{FI}}\\
\hline
reference 	& 157 	& - 	& - 	& 1971.70 	& 156	& - 	& - 	& 2199.84 \\
$2 \tau$ 	& 88 	& 0.006 	& 0.014 	& 1098.32 	& 87		& 0.006 	& 0.015 	& 1251.94 \\
$5 \tau$ 	& 43 	& 0.031 	& 0.061 	& 542.58 	& 43		& 0.029 	& 0.062 	& 655.12 \\
\hline
\multicolumn{9}{|c|}{\textit{SI}}\\
\hline
$\tau$ 		& 100 	& 0.027 	& 0.022 	& 84.67 		& 100	& 0.024 	& 0.024 	& 229.95 \\
$2 \tau$ 	& 50 	& 0.048 	& 0.031 	& 54.98 		& 50		& 0.042 	& 0.033 	& 175.81 \\
$5 \tau$ 	& 20 	& 0.106 	& 0.053 	& 35.15 		& 20		& 0.090 	& 0.044 	& 104.50 \\
\hline
\end{tabular}
\caption{Test 1. Solution time with number of total iterations and errors (\%) at final time for first and second species.  
\textit{FI} (coupled) and \textit{SI} (uncoupled) schemes}
\label{table-fisi-1}
\end{table}

\begin{figure}[h!]
\centering
\begin{subfigure}{0.48\textwidth}
\includegraphics[width=1\linewidth]{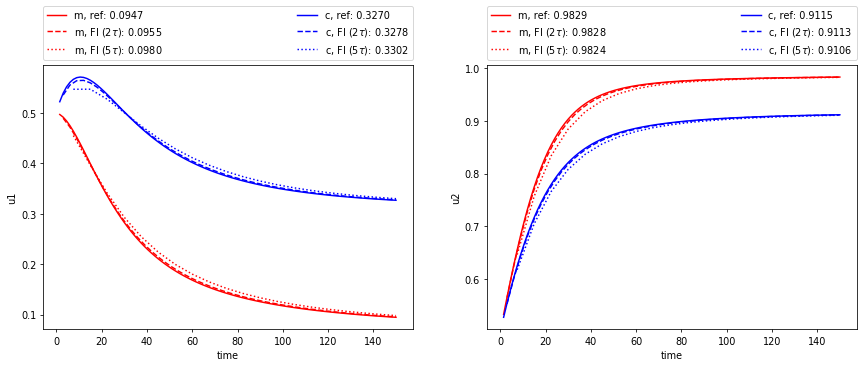}
\includegraphics[width=1\linewidth]{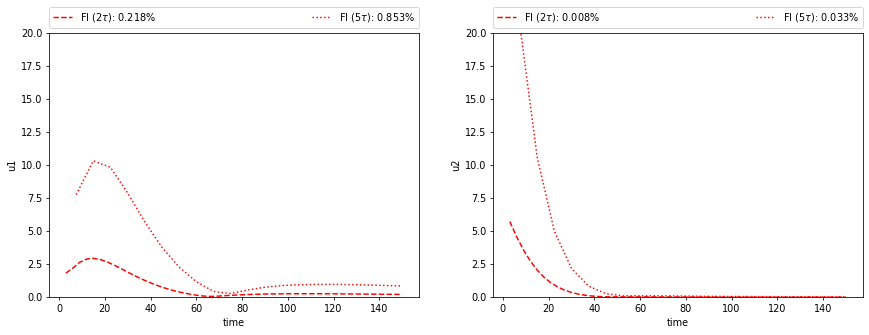}
\caption{Test 2a \textit{FI}}
\end{subfigure}
\,\,\,\,
\begin{subfigure}{0.48\textwidth}
\includegraphics[width=1\linewidth]{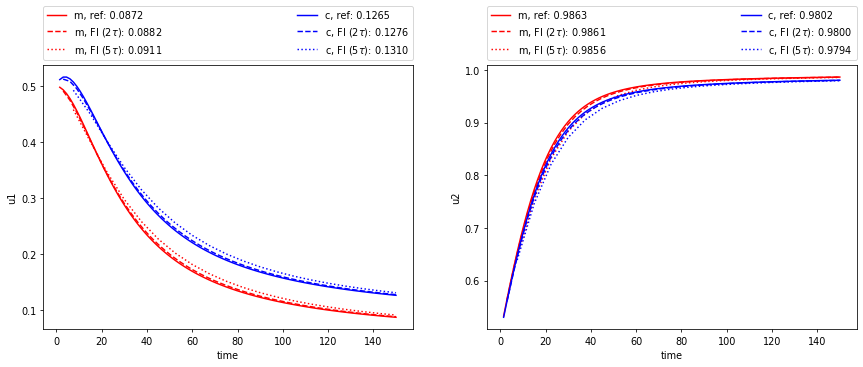}
\includegraphics[width=1\linewidth]{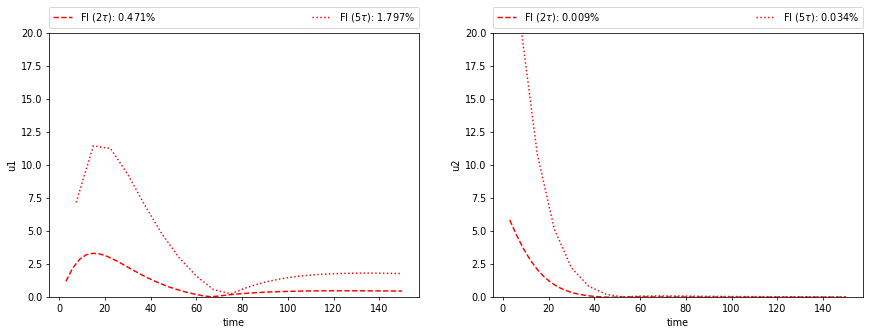}
\caption{Test 2b \textit{FI}}
\end{subfigure}\\
\begin{subfigure}{0.48\textwidth}
\includegraphics[width=1\linewidth]{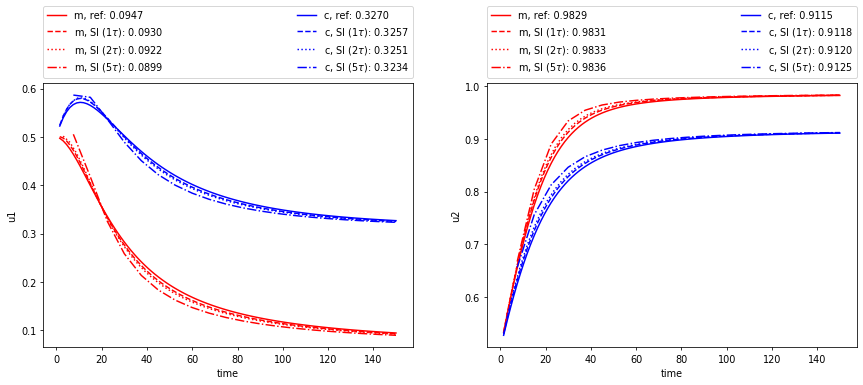}
\includegraphics[width=1\linewidth]{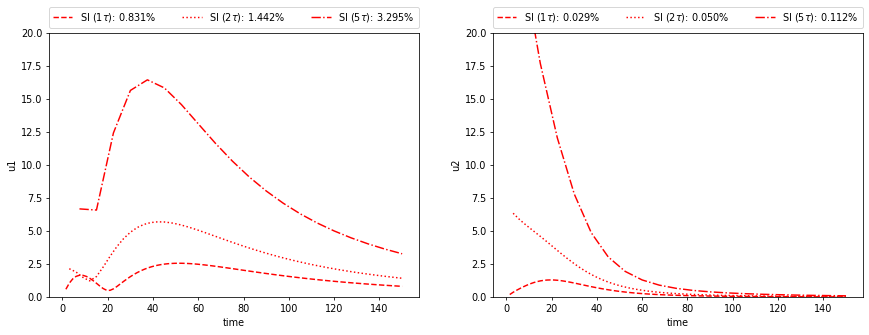}
\caption{Test 2a \textit{SI}}
\end{subfigure}
\,\,\,\,
\begin{subfigure}{0.48\textwidth}
\includegraphics[width=1\linewidth]{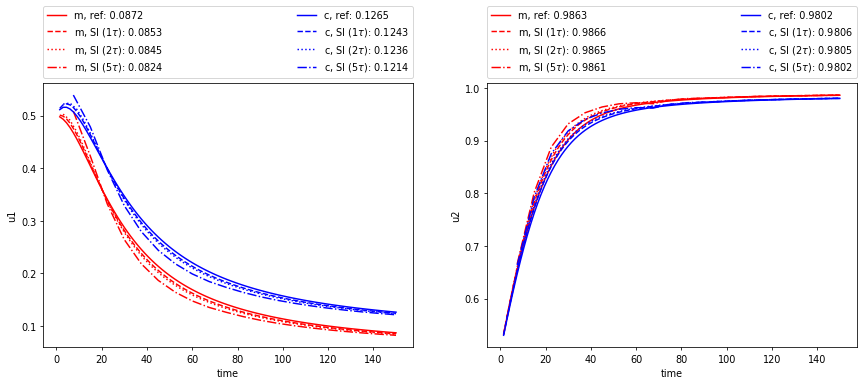}
\includegraphics[width=1\linewidth]{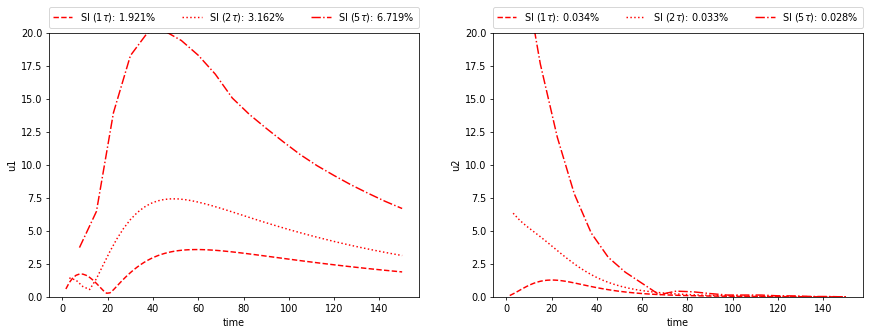}
\caption{Test 2b \textit{SI}}
\end{subfigure}
\caption{Test 2. Solution averages (first row) and errors (second row) for different time step size, $\tau = 1/M$, $2 \tau = 2/M$, and $10 \tau = 10/M$ ($M=100$)}
\label{sol-fi-2}
\end{figure}

\begin{table}[h!]
\centering
\begin{tabular}{|c|cccc|cccc|}
\hline
&
 \multicolumn{4}{|c|}{Test 2a}&
 \multicolumn{4}{|c|}{Test 2b}\\ 
 & $N$ 
& $e_1$ (\%)  & $e_2$( \%) & time(sec) 
& $N$ 
& $e_1$ (\%)  & $e_2$ (\%) & time(sec)  \\
\hline
\multicolumn{9}{|c|}{\textit{FI}}\\
\hline
reference 	& 200 	& - 	& - 	& 2313.36 	& 200	& - 	& - 	& 2684.84 \\
$2 \tau$ 	& 105 	& 0.218 	& 0.008 	& 1302.20 	& 105	& 0.471 	& 0.009 	& 1560.85 \\
$5 \tau$ 	& 47 	& 0.853 	& 0.033 	& 632.55 	& 66		& 1.797 	& 0.034 	& 934.96 \\
\hline
\multicolumn{9}{|c|}{\textit{SI}}\\
\hline
$\tau$ 		& 100 	& 0.831 	& 0.029 	& 121.18 	& 100	& 1.921 	& 0.034 	& 405.20 \\
$2 \tau$ 	& 50 	& 1.442 	& 0.050 	& 133.41 	& 50		& 3.162 	& 0.033 	& 247.99 \\
$5 \tau$ 	& 20 	& 3.295 	& 0.112 	& 59.91 		& 20		& 6.719 	& 0.028 	& 129.65 \\
\hline
\end{tabular}
\caption{Test 2. Solution time with number of total iterations and errors (\%) at final time for first and second species.  
\textit{FI} (coupled) and \textit{SI} (uncoupled) schemes}
\label{table-fisi-2}
\end{table}

\begin{figure}[h!]
\centering
\begin{subfigure}{0.48\textwidth}
\includegraphics[width=1\linewidth]{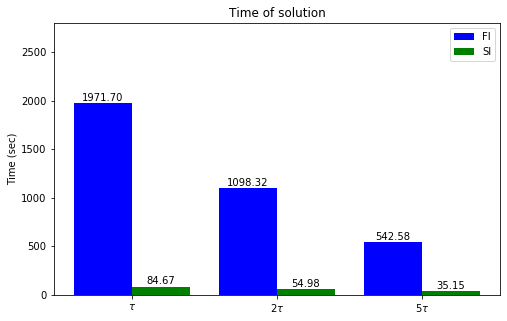}
\caption{Test 1a}
\end{subfigure}
\,\,\,\,
\begin{subfigure}{0.48\textwidth}
\includegraphics[width=1\linewidth]{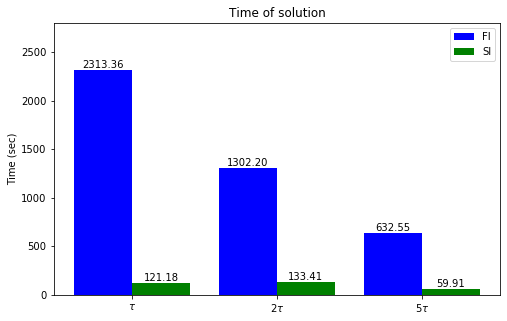}
\caption{Test 2a}
\end{subfigure}
\\
\begin{subfigure}{0.48\textwidth}
\includegraphics[width=1\linewidth]{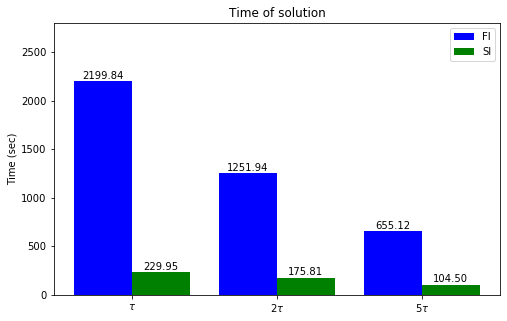}
\caption{Test 1b}
\end{subfigure}
\,\,\,\,
\begin{subfigure}{0.48\textwidth}
\includegraphics[width=1\linewidth]{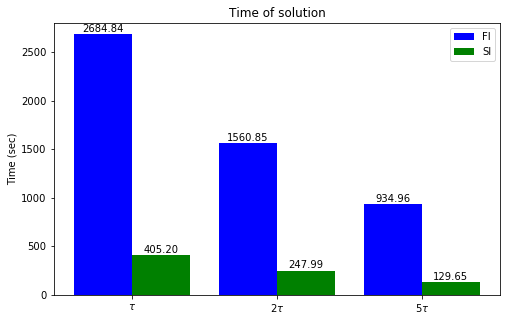}
\caption{Test 2b}
\end{subfigure}
\caption{Simulations time for \textit{FI} (coupled) and \textit{SI} (uncoupled) schemes} 
\label{sol-fisi-t}
\end{figure}

The average solutions for the first and second species in subdomains $\Omega_m$ and $\Omega_c$ are presented in Figures \ref{sol-fi-1} and \ref{sol-fi-2} for Test 1 and Test 2, respectively. The results using the \textit{FI} scheme are shown in  (a) and (b),  and solutions obtained via \textit{SI} scheme are exhibited in (c) and (d) of these separate figure sets. 
The numerical solutions are given using three time steps sizes, viz.,  $\tau = 1/N_t$, $2 \tau = 2/N_t$, and $10 \tau = 10/N_t$ with $N_t=100$ (the total number of time step iterations). It can be seen that the larger time stepping size influences the solution noticeably.  We observe that a larger difference occurs at the beginning as well as in the middle of the simulations. However the final solutions approach almost the same value. In Figures \ref{sol-fi-1} and \ref{sol-fi-2}, we have also depicted the error dynamics for the first and second species.  As a reference solution in our error calculations, we have used the solution using \textit{FI} scheme with smallest time step size $\tau = 1/N_t$. 
The errors are then calculated using the following formulas: 
\[
e_k = \left( \frac{\int_{\Omega} (u^k_{ref} - u^k)^2 \ dx}{\int_{\Omega} (u^k_{ref})^2 \ dx} \right)^{1/2}, 
\quad k = 1,2 
\]
where $u^k_{ref}$ represents the reference solution. 

From Figure \ref{sol-fi-1} for Test 1,  we observe that one can obtain good results using an uncoupled \textit{SI} scheme with larger time stepping size.  Also, from Figure \ref{sol-fi-2}  for Test 2, we see that better numerical results can be achieved using uncoupled strategy. But the latter has a higher influence on the time stepping size, especially at the beginning and middle time instances during the simulations. Furthermore, smaller diffusion coefficients lead to smaller errors. 

In Tables \ref{table-fisi-1}  and \ref{table-fisi-2}, we record the solution time along with the number of total iterations and errors (\%) at the final time step for first and second species. Results are provided for both \textit{FI} (coupled) and \textit{SI} (uncoupled) schemes.  The total number of time steps are: 100 for $\tau$, 50 for $2\tau$, 20 for $5\tau$ and 10 for $10\tau$.  In the coupled \textit{FI} scheme, the nonlinear iterations effect the total number of iterations which leads to the longer simulation time. Size of the system that are solved in each time step and nonlinear iterations have a greater impact on the solution time.  Also, in the coupled \textit{FI} scheme,  the size of system is $DOF_h  = L \cdot 69948$, where $L=2$ (two-species). 
An illustration of the computational efficiency of the uncoupled scheme is furnished in Figure  \ref{sol-fisi-t}.  As seen from the values, one can obtain a huge simulation time reduction even with the use of a small time stepping size.  For example,  the solution time for coupled \textit{FI} scheme is 30-40 minutes  for Test 1a and 2a with small diffusion, but using uncoupled \textit{SI} scheme we can obtain accurate solution by running the simulation for about 1-2 minutes. We remark that for larger diffusion the simulation time is slightly larger in both schemes.  Specifically, we obtain solutions in 35-45 minutes using \textit{FI} scheme and 3-7 minutes via \textit{SI} scheme for Test 1b and 2b. 

Below, we discuss the possibility reducing the solution time further by applying a multiscale order reduction technique for space approximation and using the uncoupled \textit{SI} scheme.

\subsection{Multiscale method for approximation by space}

We now present numerical investigation using the proposed multiscale solver for the solution of the problem on the coarse grid.  We use an uncoupled \textit{SI} scheme with fixed time step size $\tau = 1/N_t$ ($N_t = 100$). The coarse grid is fixed and contains 100 coarse cells ($10 \times 10$ coarse grid). 
We solve the problem on the fine grid and use the solution as a reference solution to calculate errors of the multiscale solver.  

\begin{figure}[h!]
\centering
\begin{subfigure}{1\textwidth}
\includegraphics[width=0.16\linewidth]{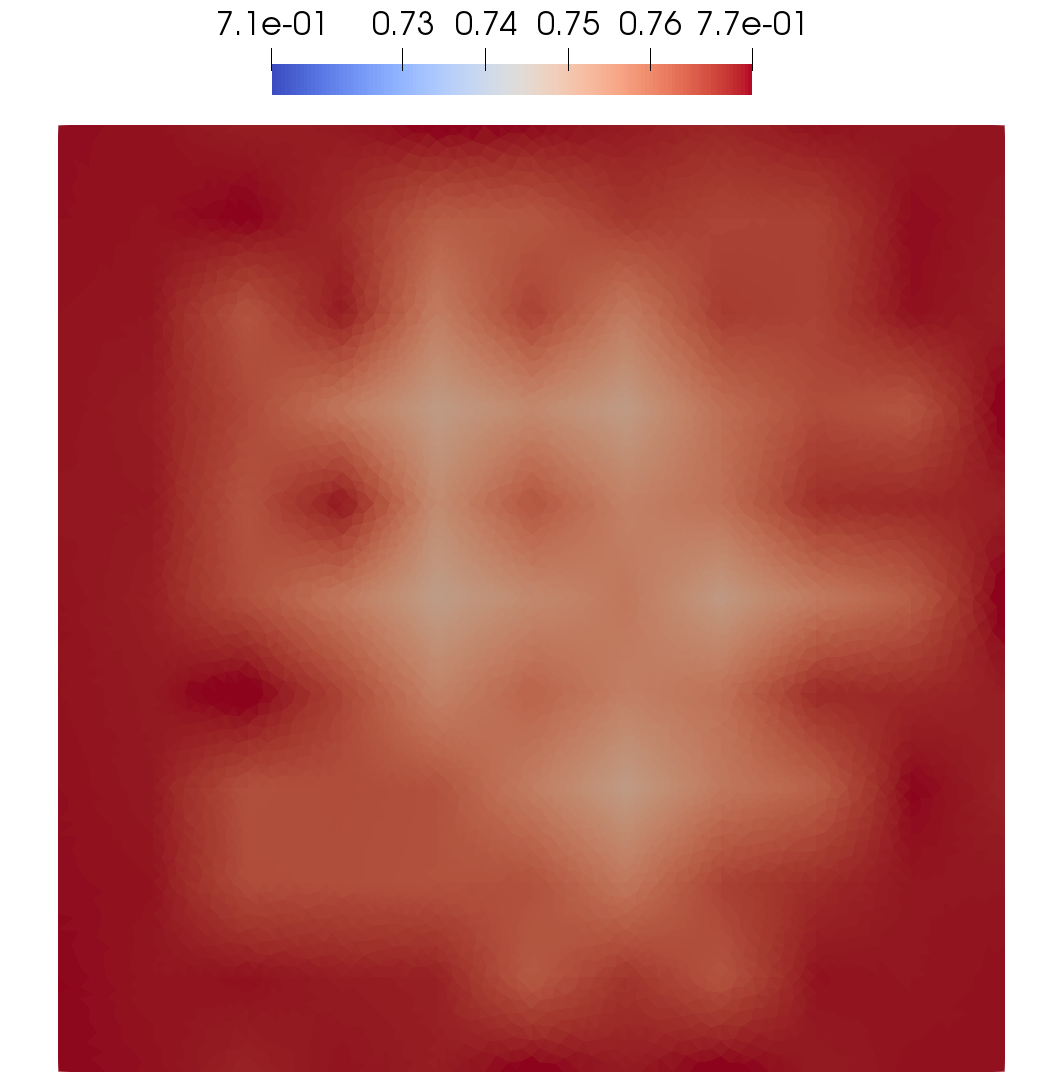}
\includegraphics[width=0.16\linewidth]{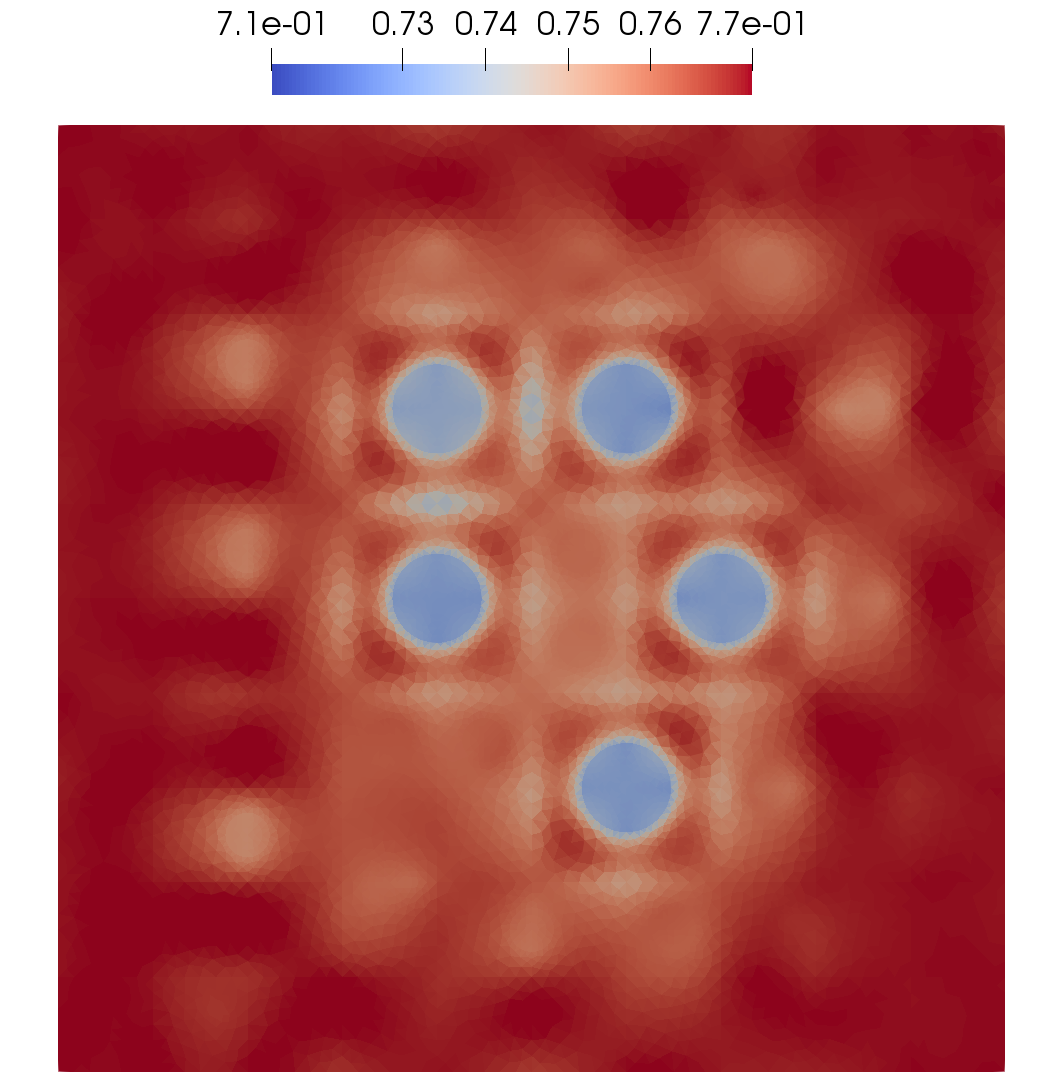}
\includegraphics[width=0.16\linewidth]{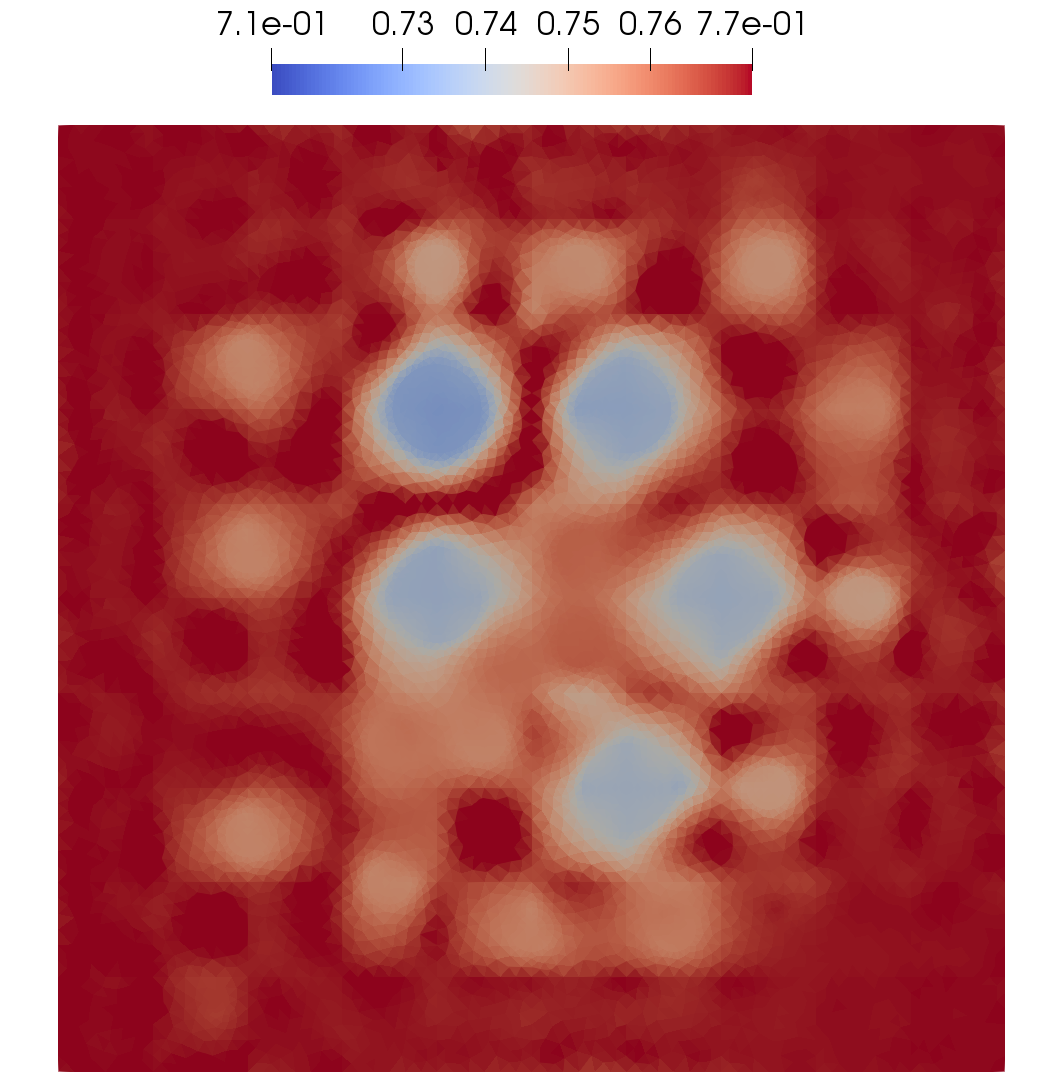}
\includegraphics[width=0.16\linewidth]{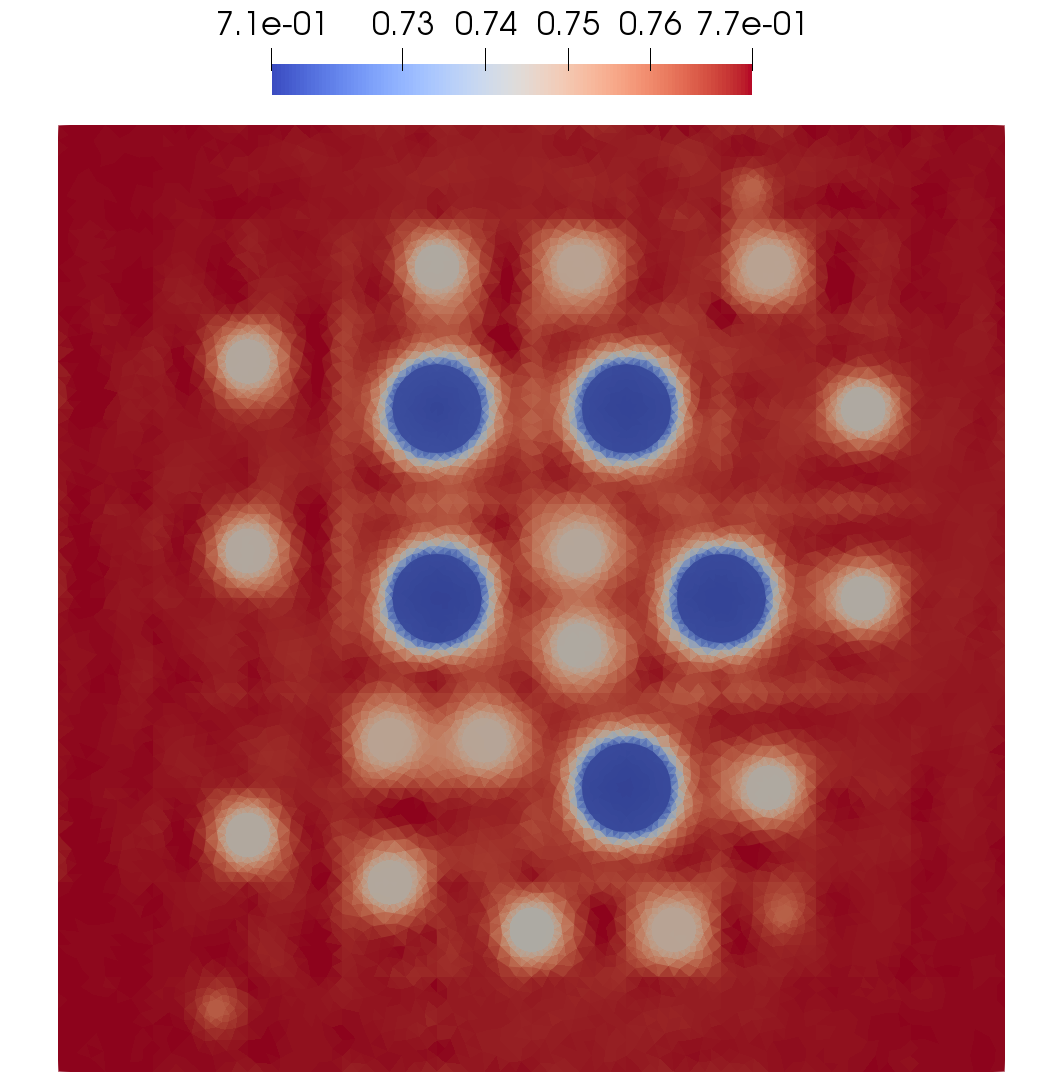}
\includegraphics[width=0.16\linewidth]{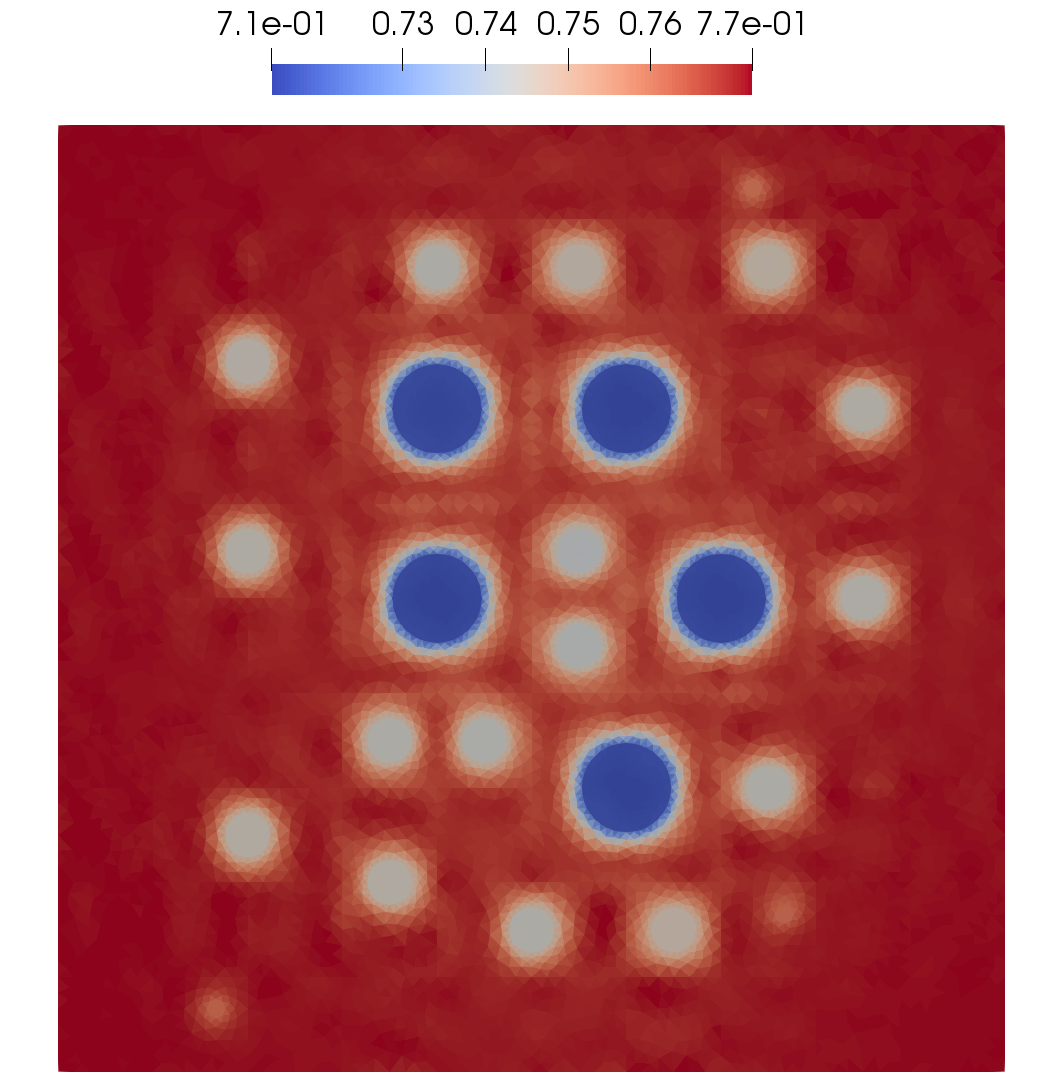}
\includegraphics[width=0.16\linewidth]{uf1}\\
\includegraphics[width=0.16\linewidth]{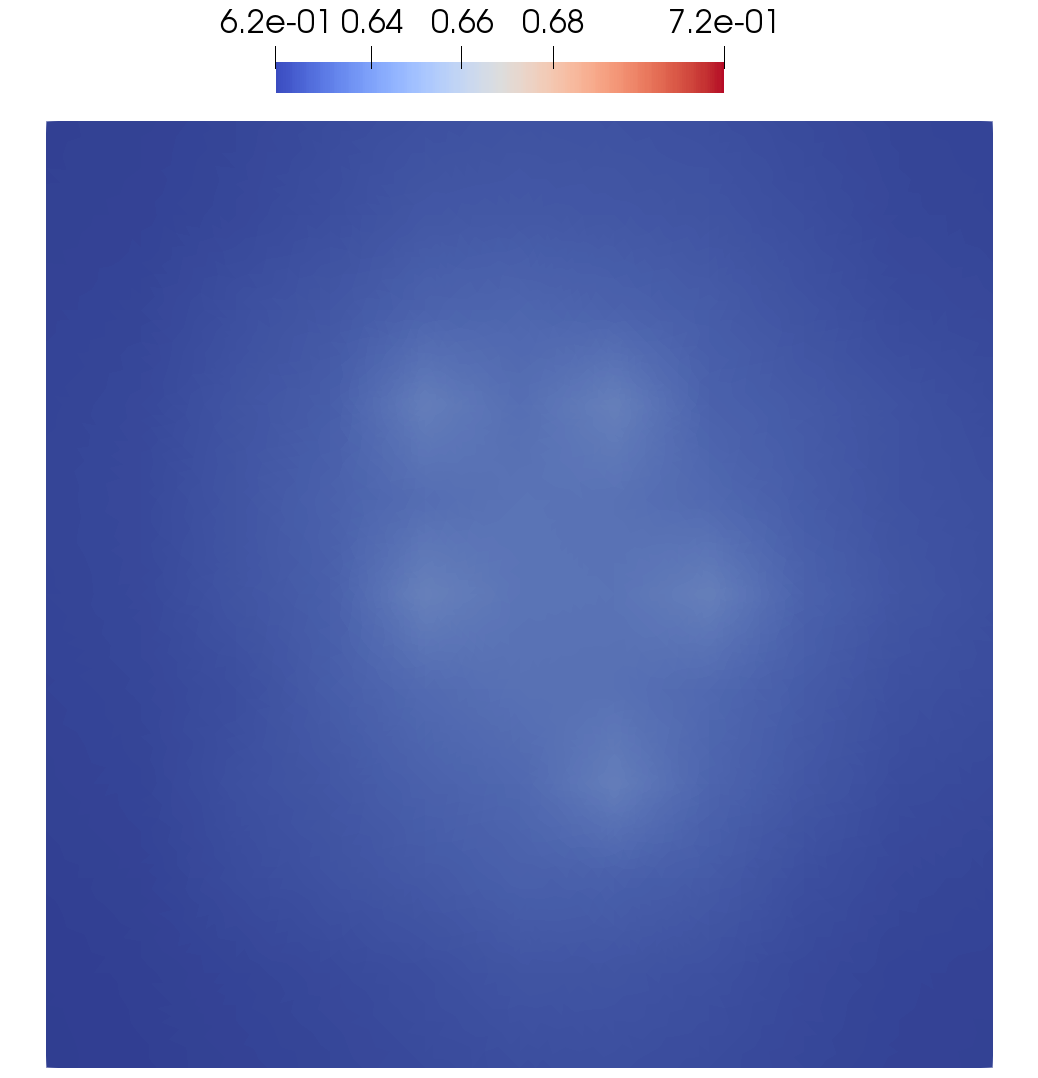}
\includegraphics[width=0.16\linewidth]{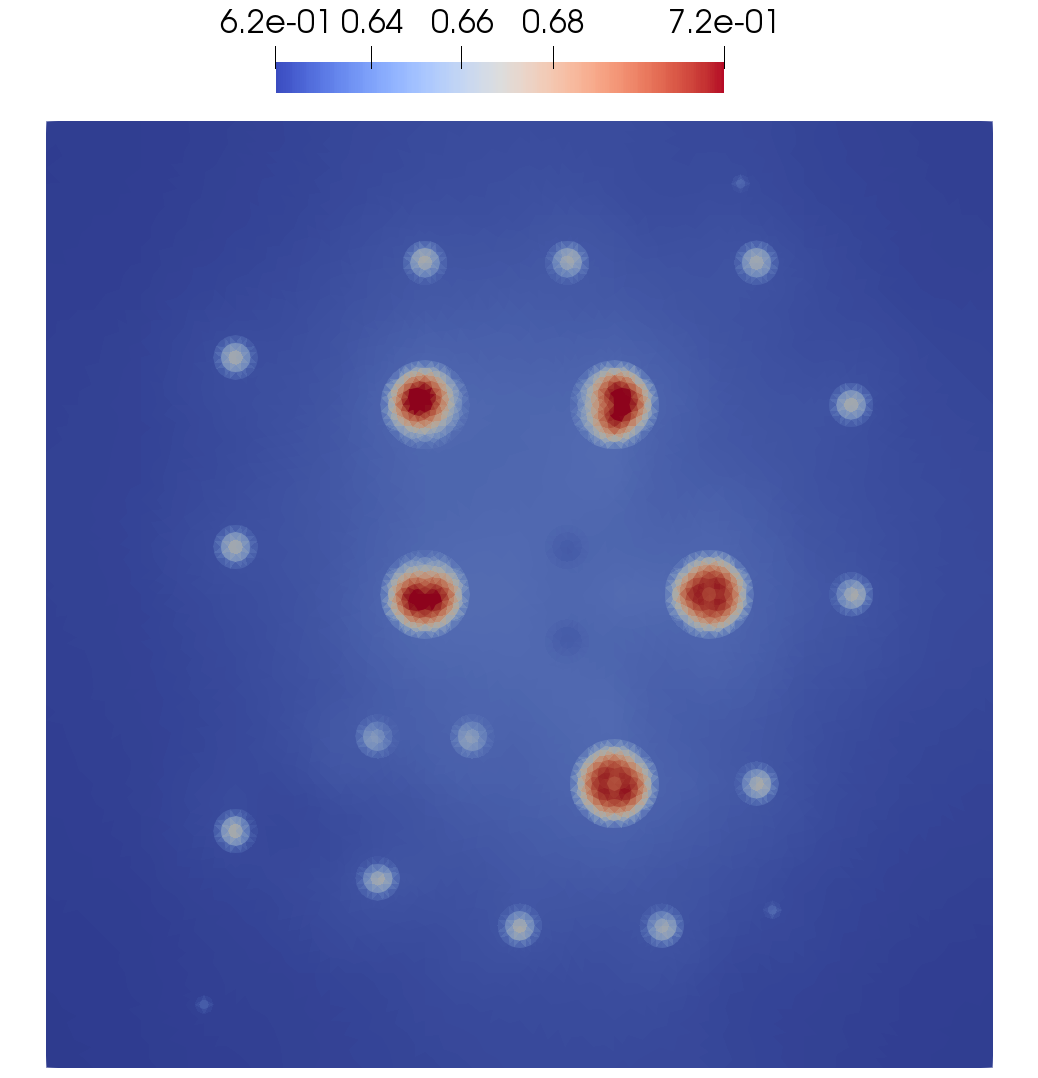}
\includegraphics[width=0.16\linewidth]{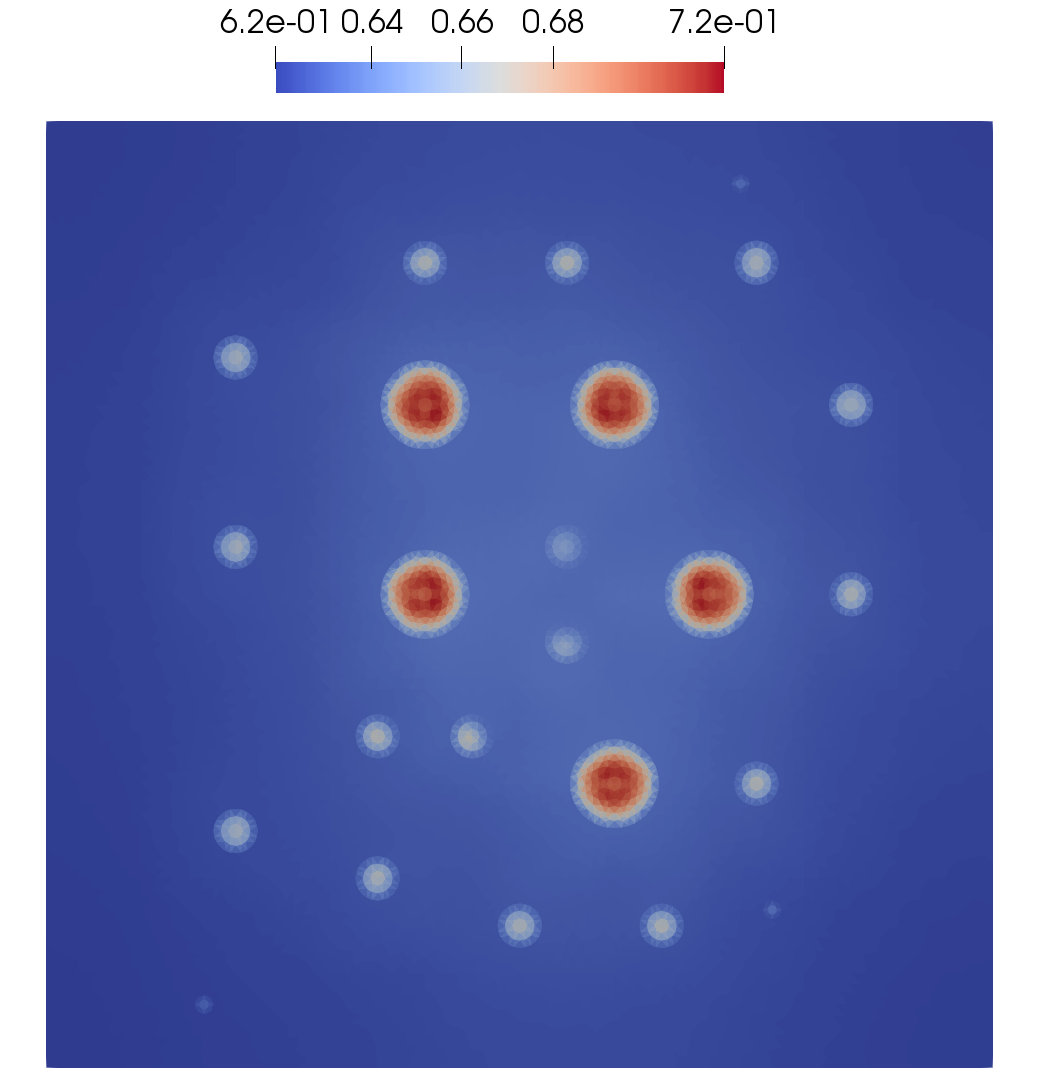}
\includegraphics[width=0.16\linewidth]{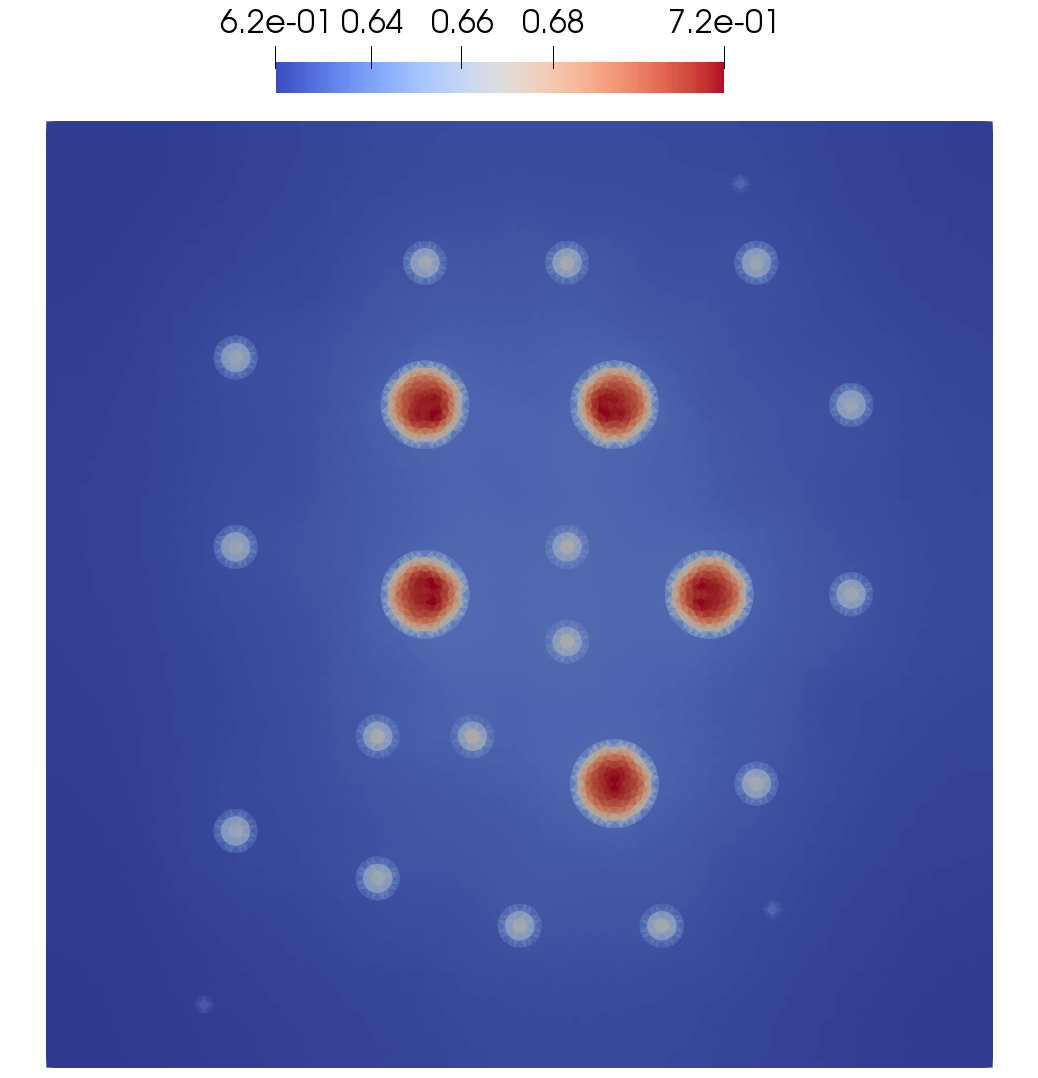}
\includegraphics[width=0.16\linewidth]{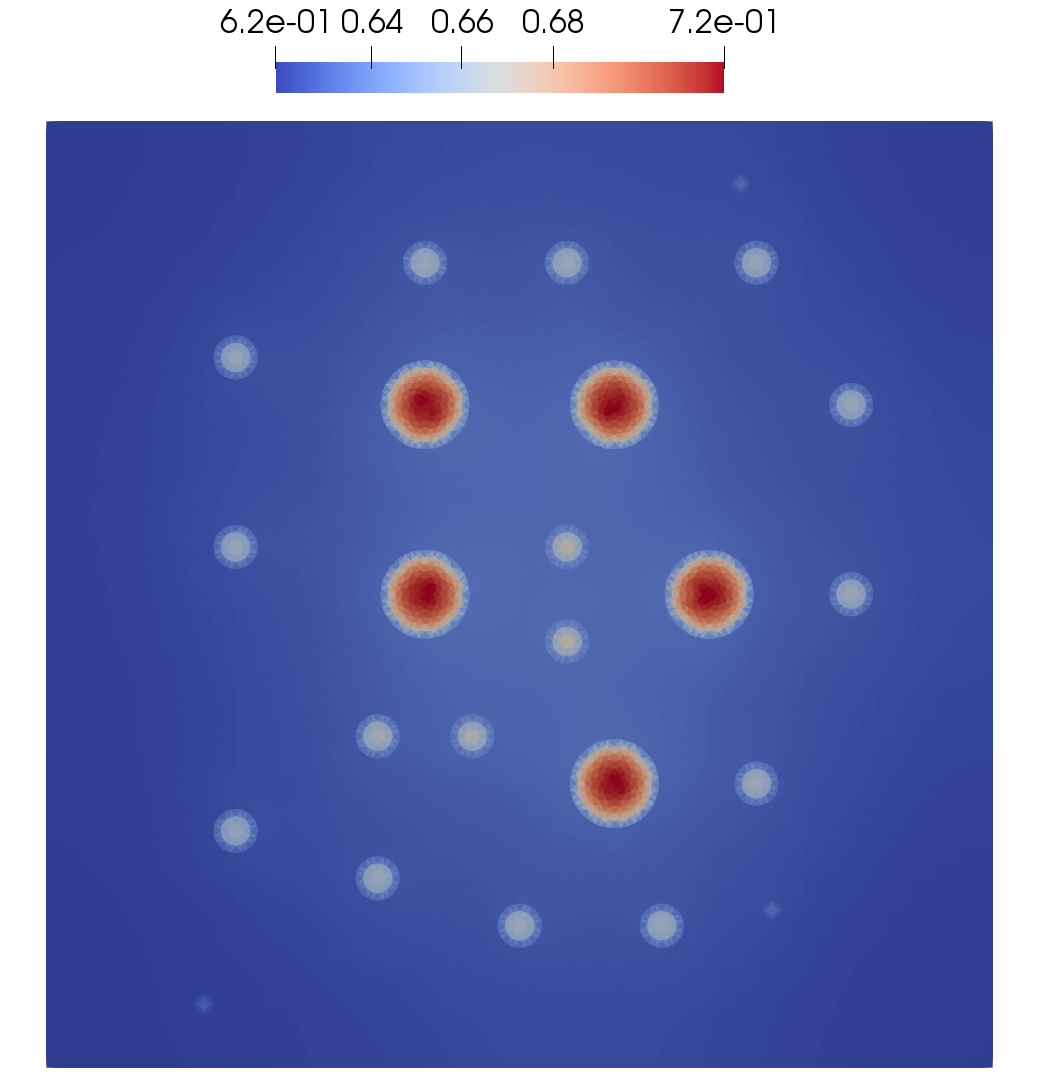}
\includegraphics[width=0.16\linewidth]{uf2}
\caption{Test 1a}
\end{subfigure}\\
\begin{subfigure}{1\textwidth}
\includegraphics[width=0.16\linewidth]{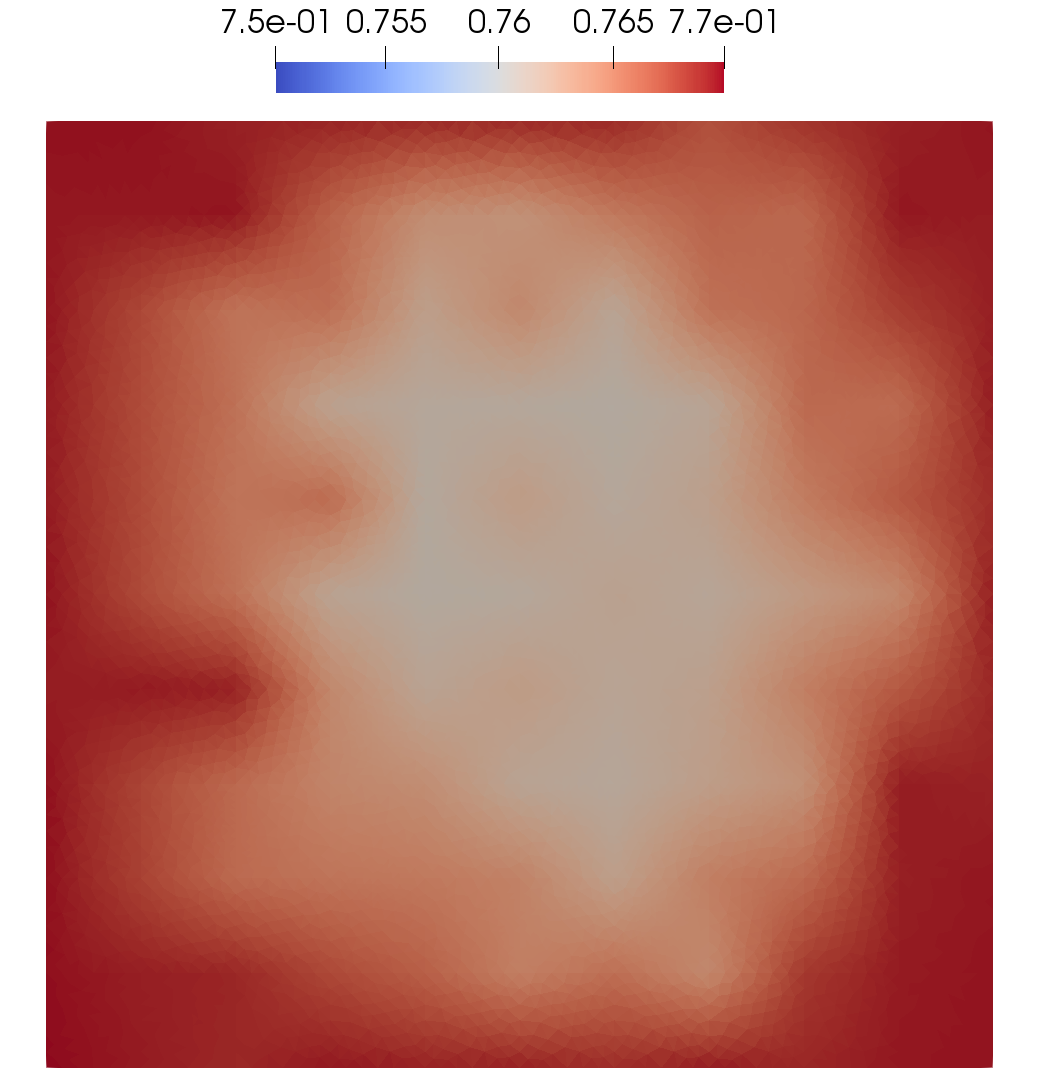}
\includegraphics[width=0.16\linewidth]{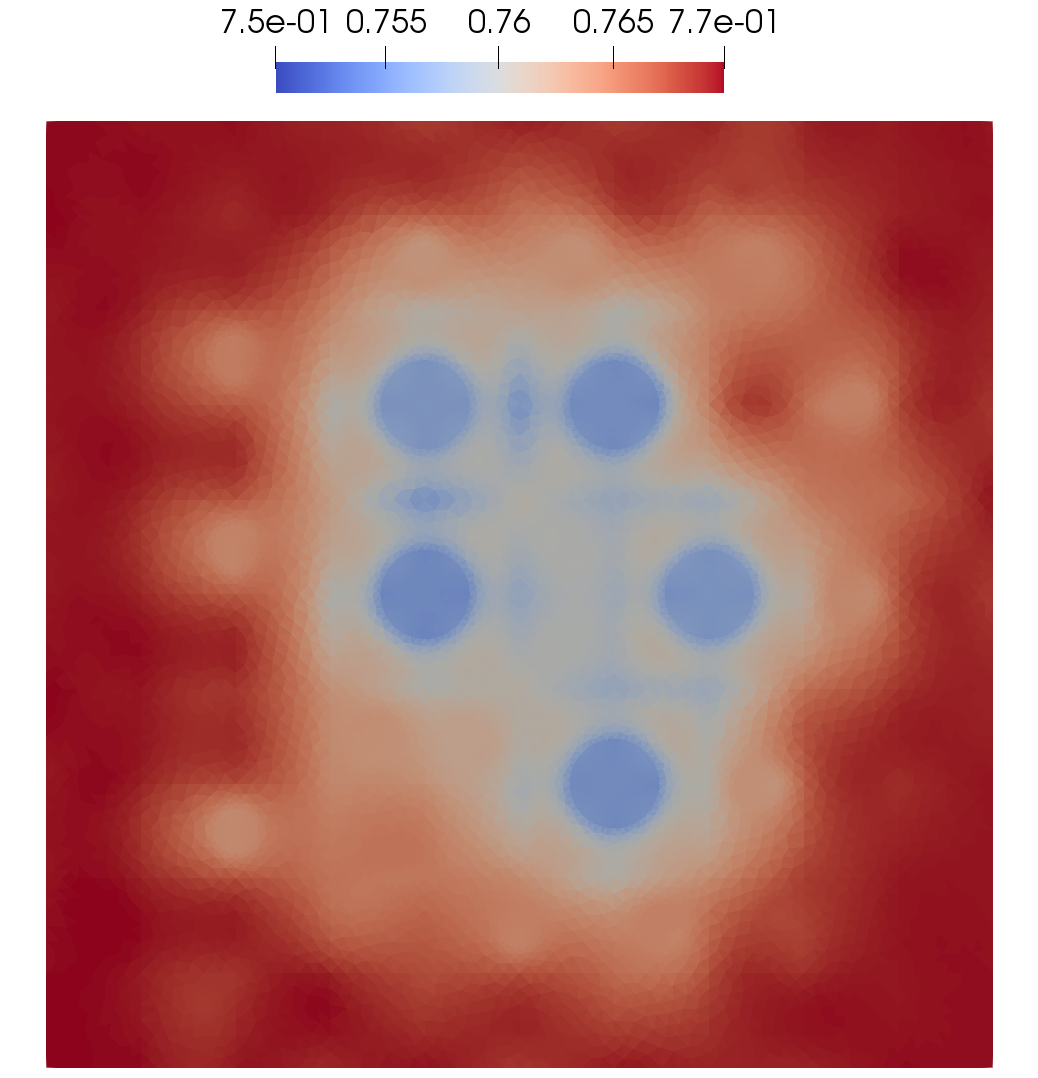}
\includegraphics[width=0.16\linewidth]{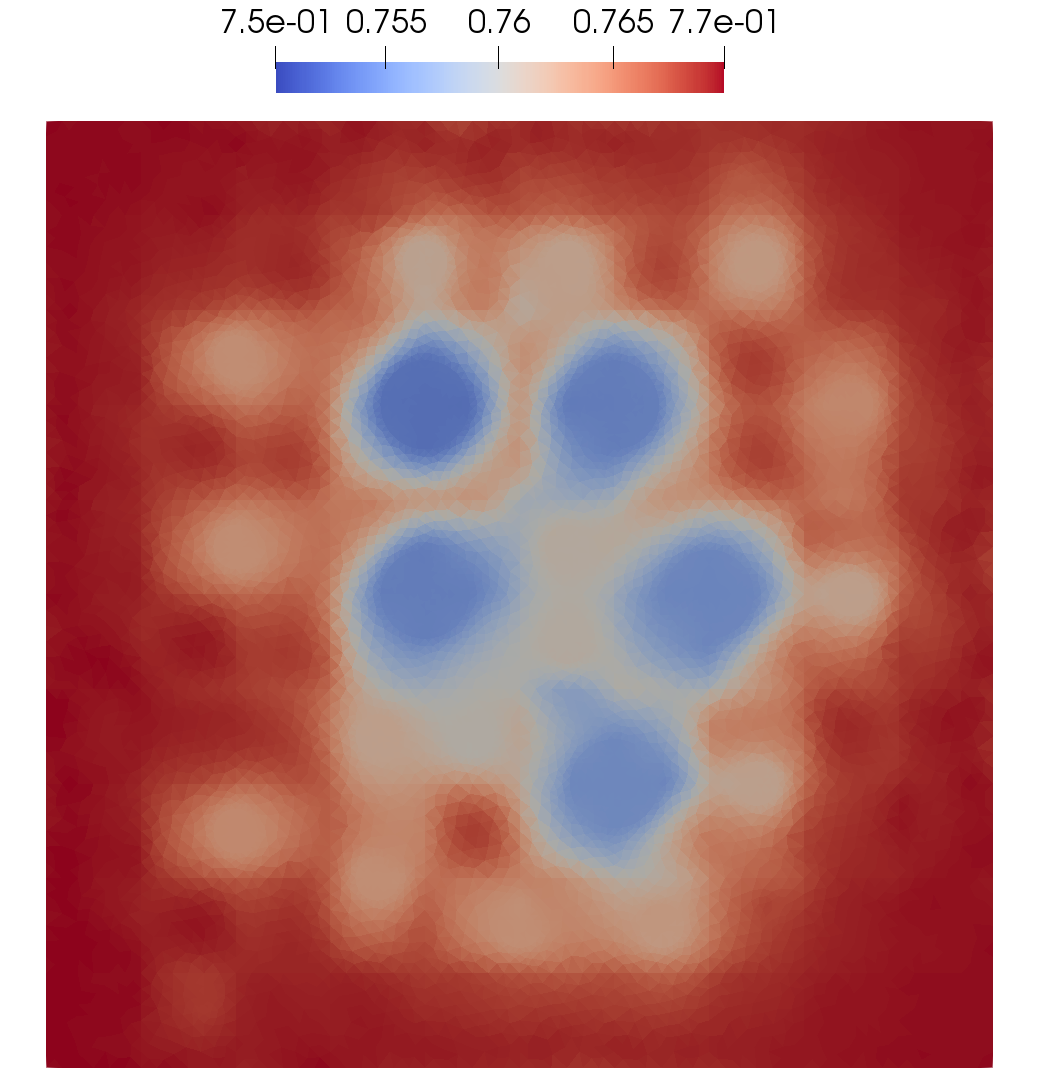}
\includegraphics[width=0.16\linewidth]{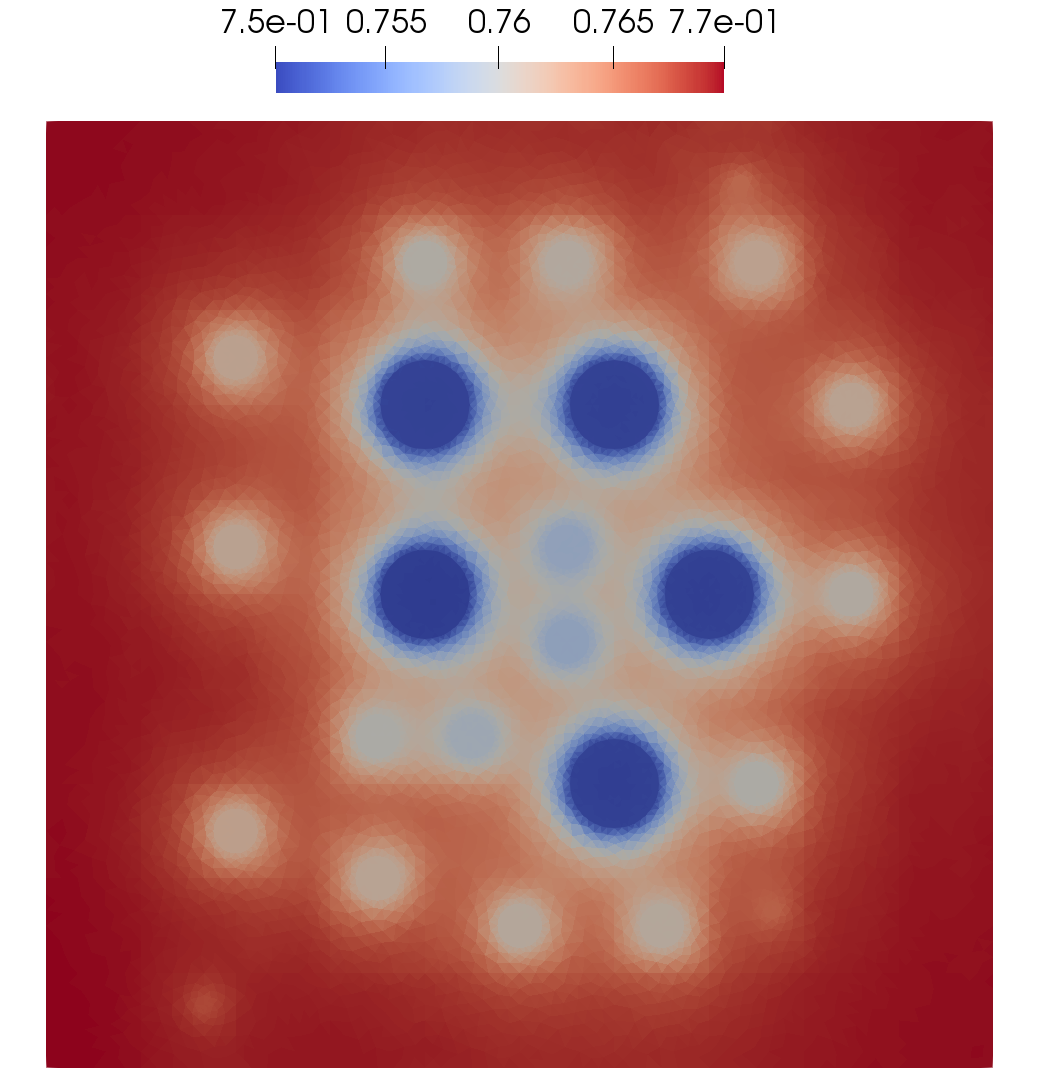}
\includegraphics[width=0.16\linewidth]{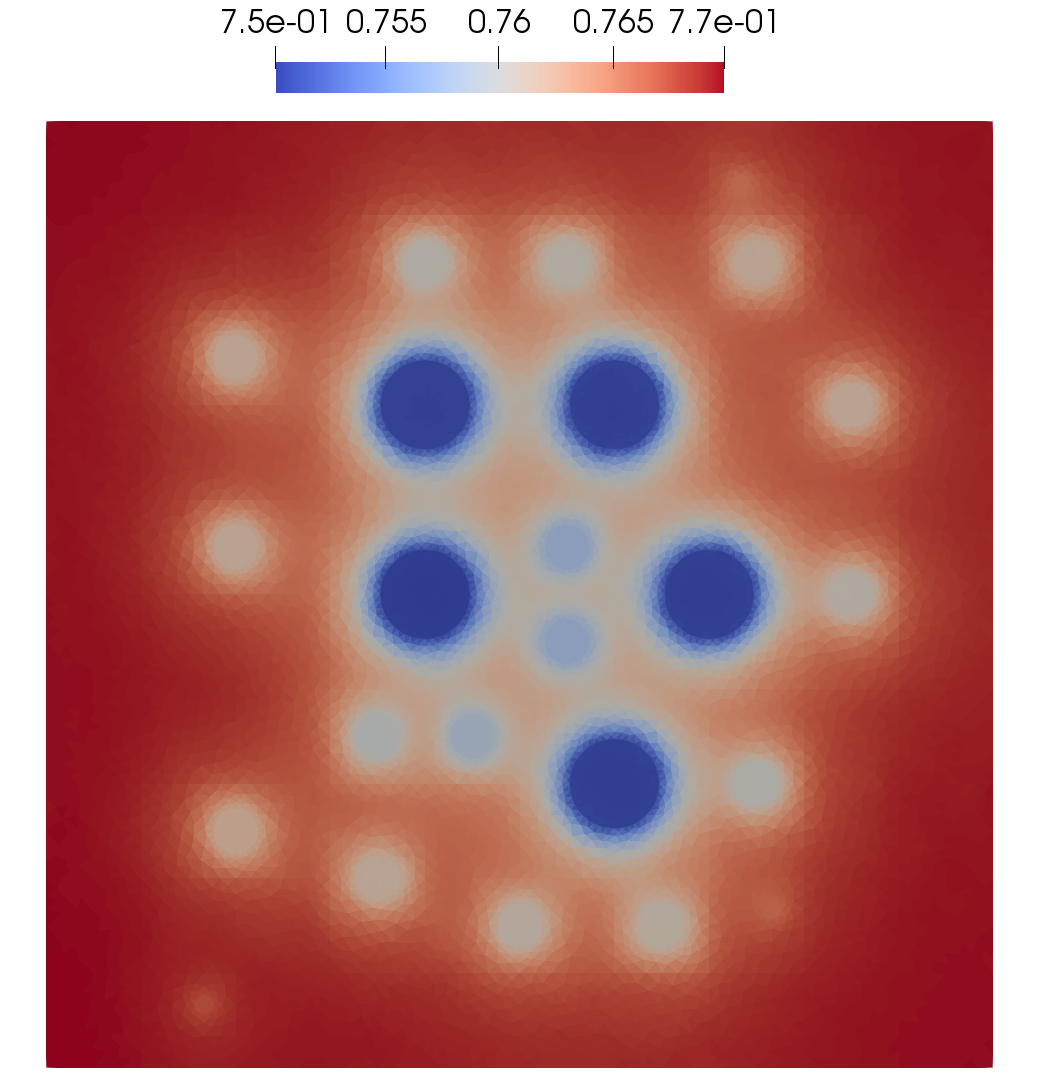}
\includegraphics[width=0.16\linewidth]{u2f1}\\
\includegraphics[width=0.16\linewidth]{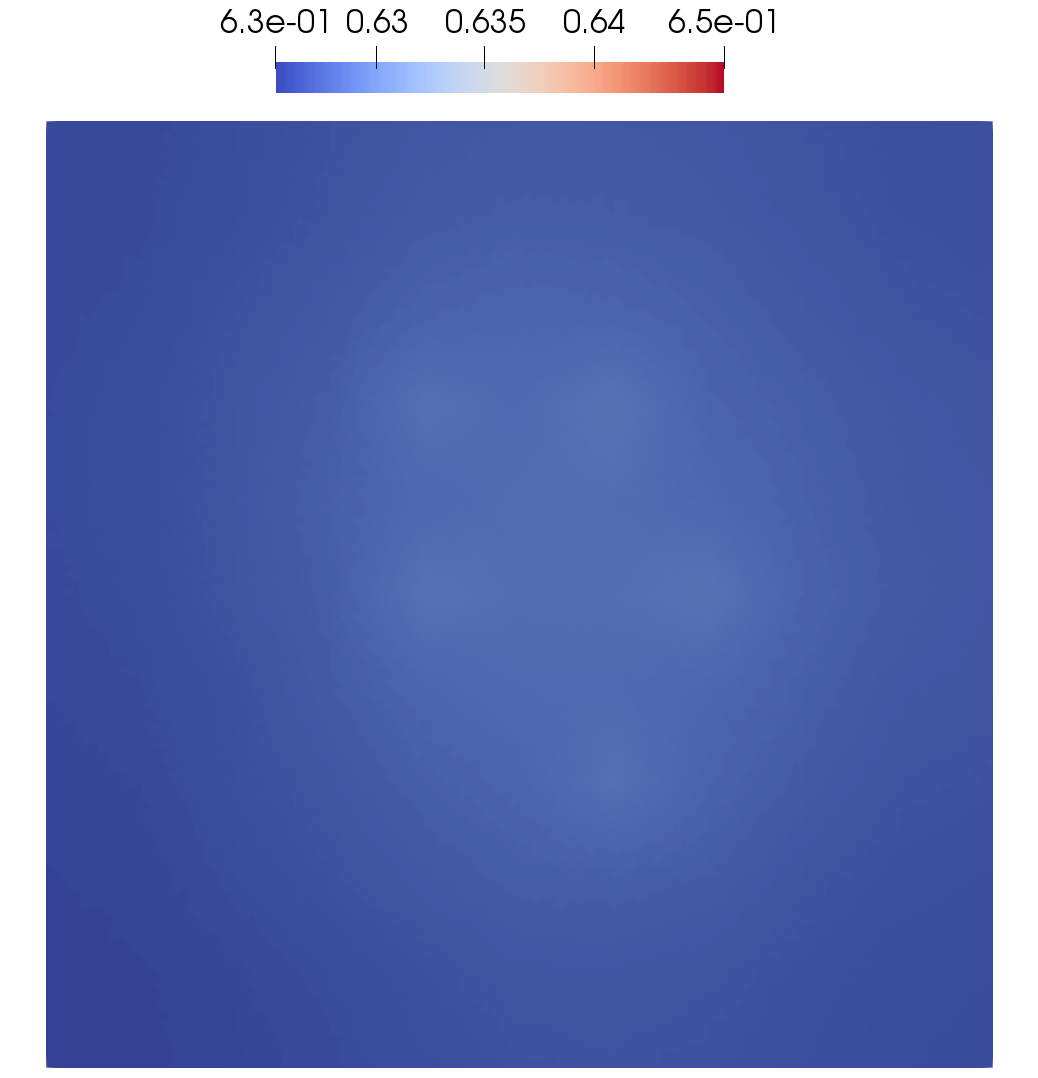}
\includegraphics[width=0.16\linewidth]{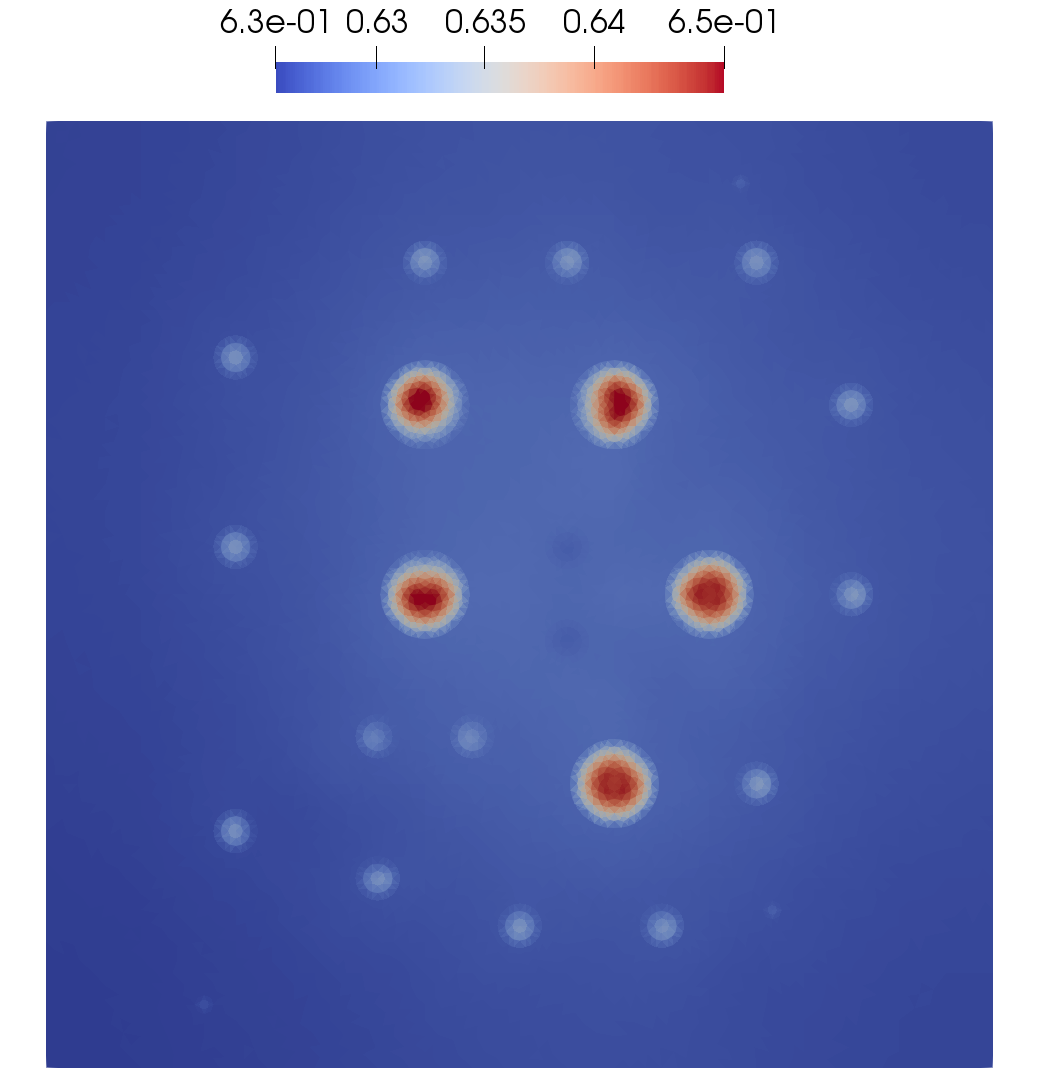}
\includegraphics[width=0.16\linewidth]{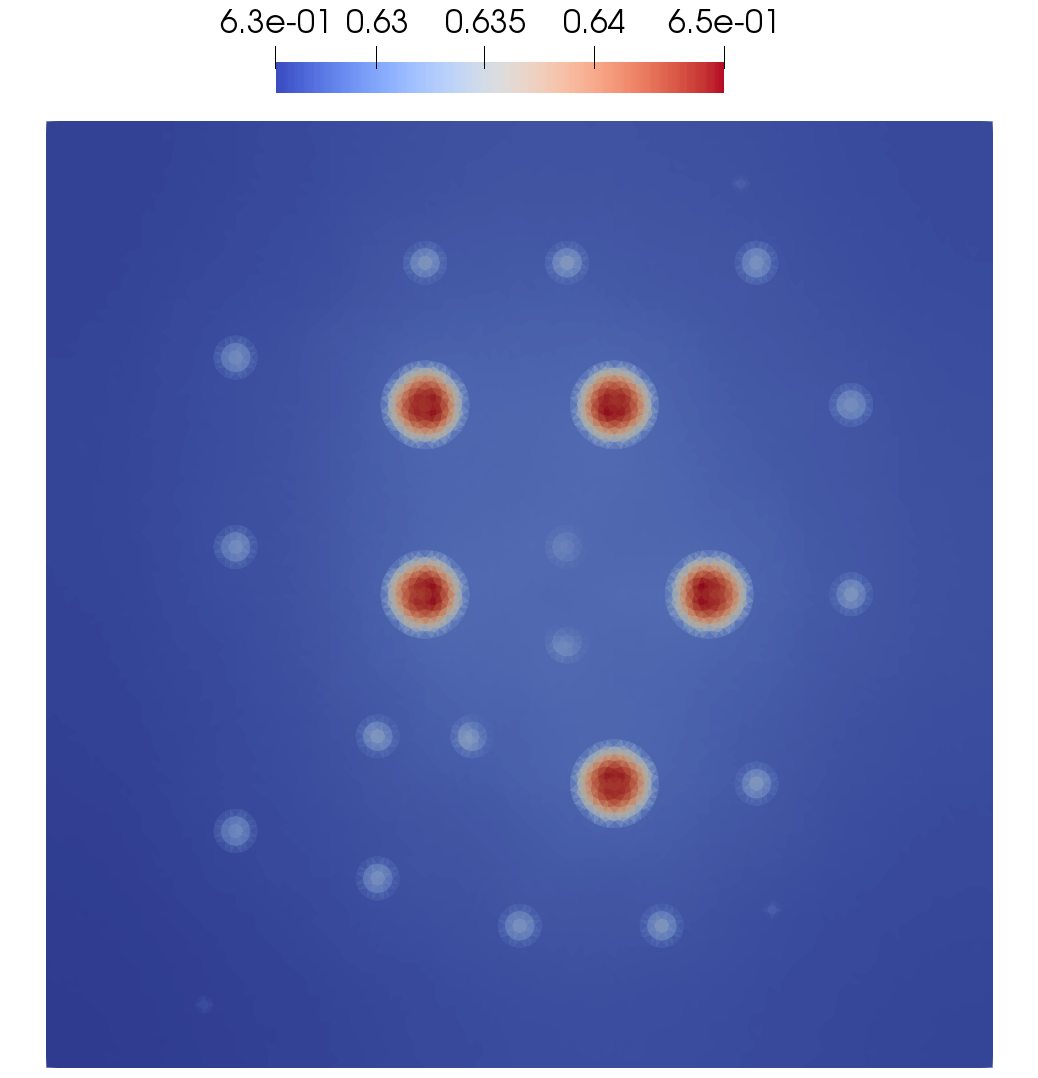}
\includegraphics[width=0.16\linewidth]{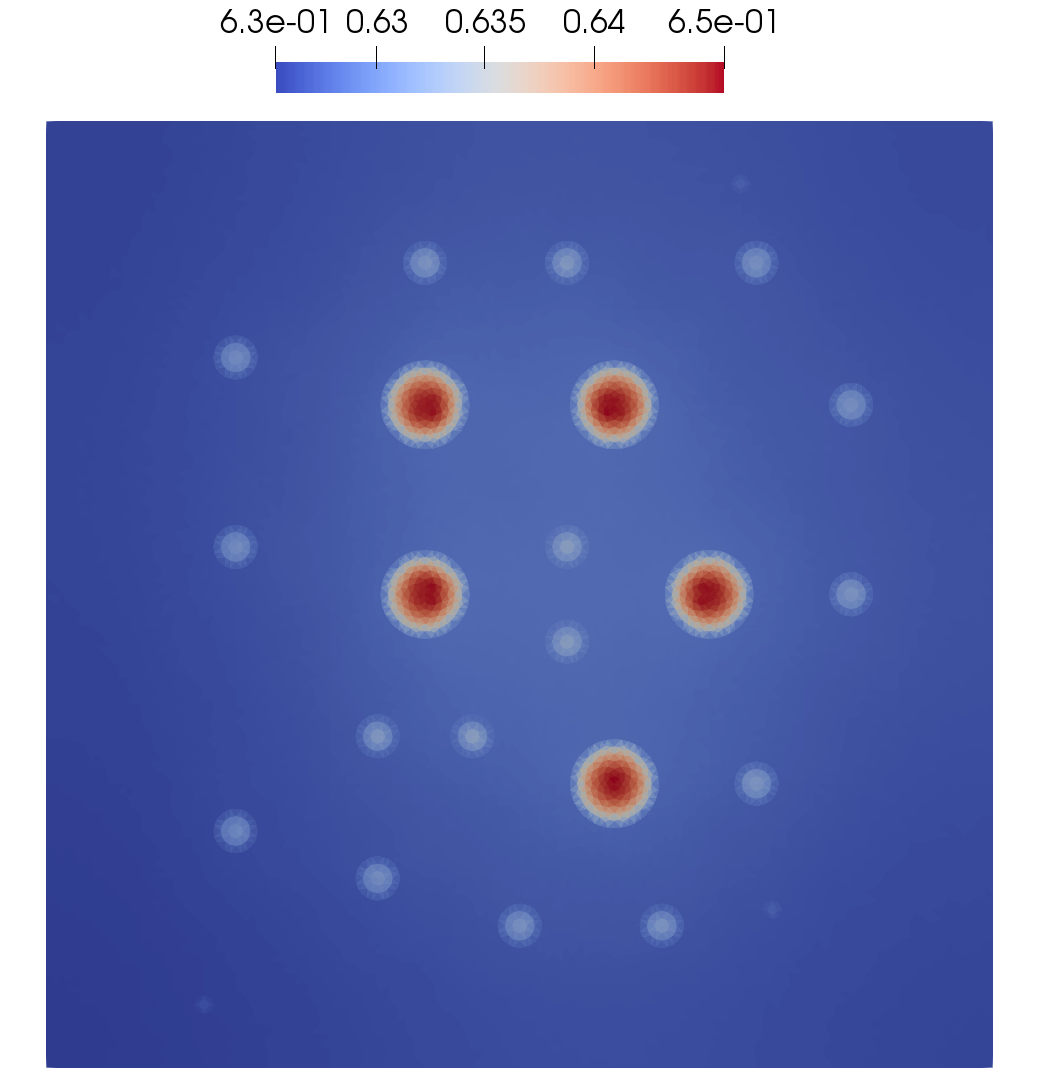}
\includegraphics[width=0.16\linewidth]{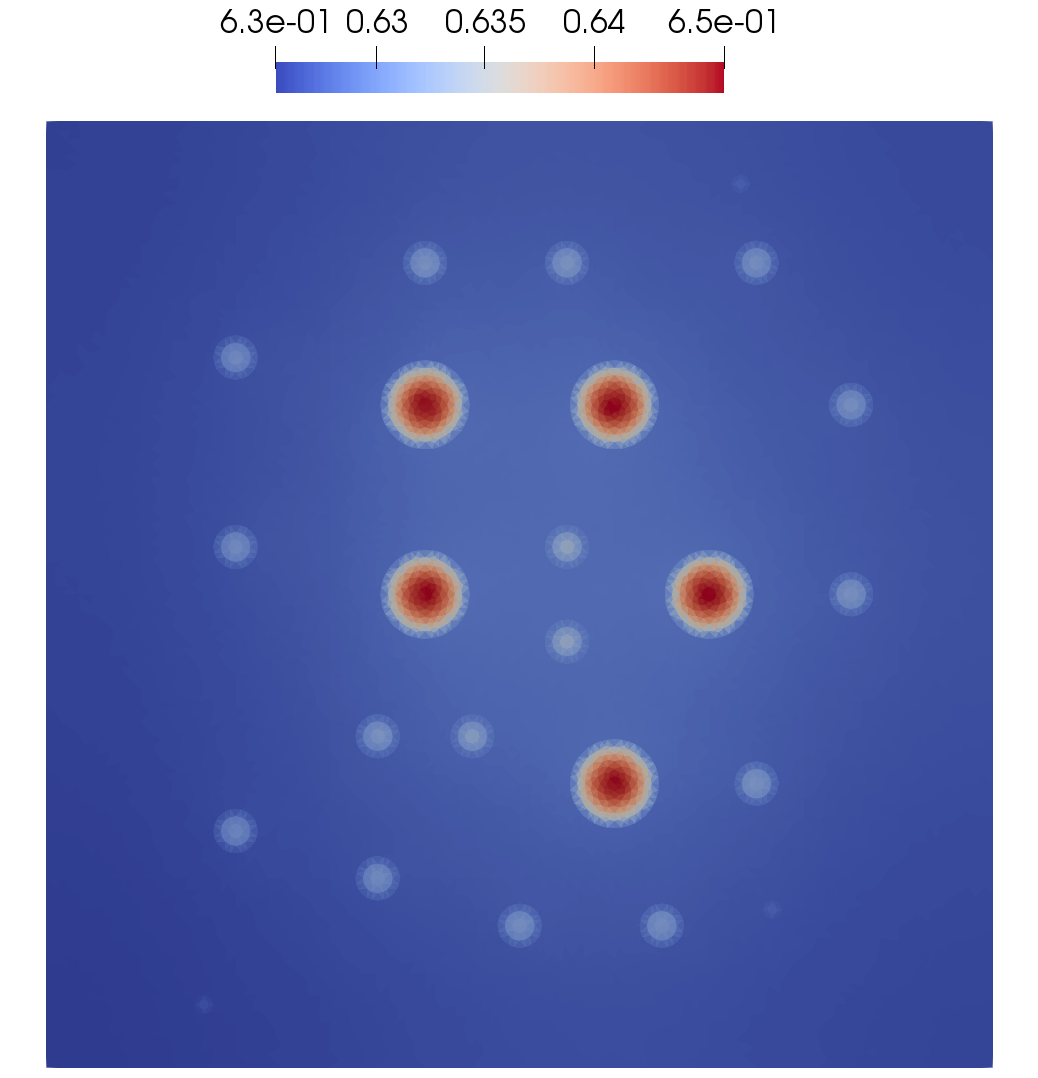}
\includegraphics[width=0.16\linewidth]{u2f2}
\caption{Test 1b}
\end{subfigure}
\caption{Test 1.  Solutions using multiscale solver for $M=1,2,4,6,8$ and fine grid (reference) solution (from left to right).  
First  and second rows: $u^1$ and $u^2$ for Test 1a (small diffusion) at final time. 
Third  and fourth rows: $u^1$ and $u^2$ for Test 1b (regular diffusion) at final time }
\label{sol1u}
\end{figure}

\begin{figure}[h!]
\centering
\begin{subfigure}{1\textwidth}
\includegraphics[width=0.16\linewidth]{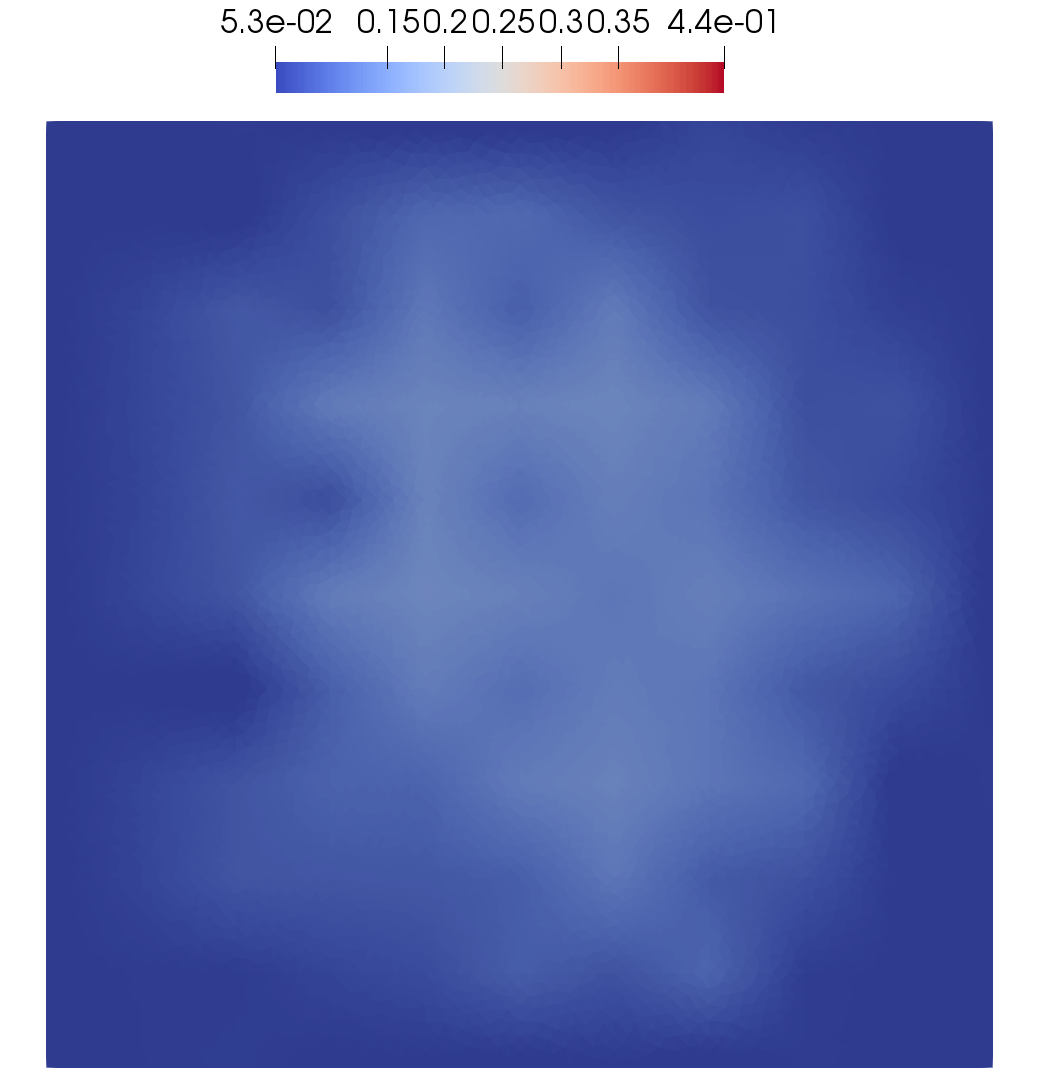}
\includegraphics[width=0.16\linewidth]{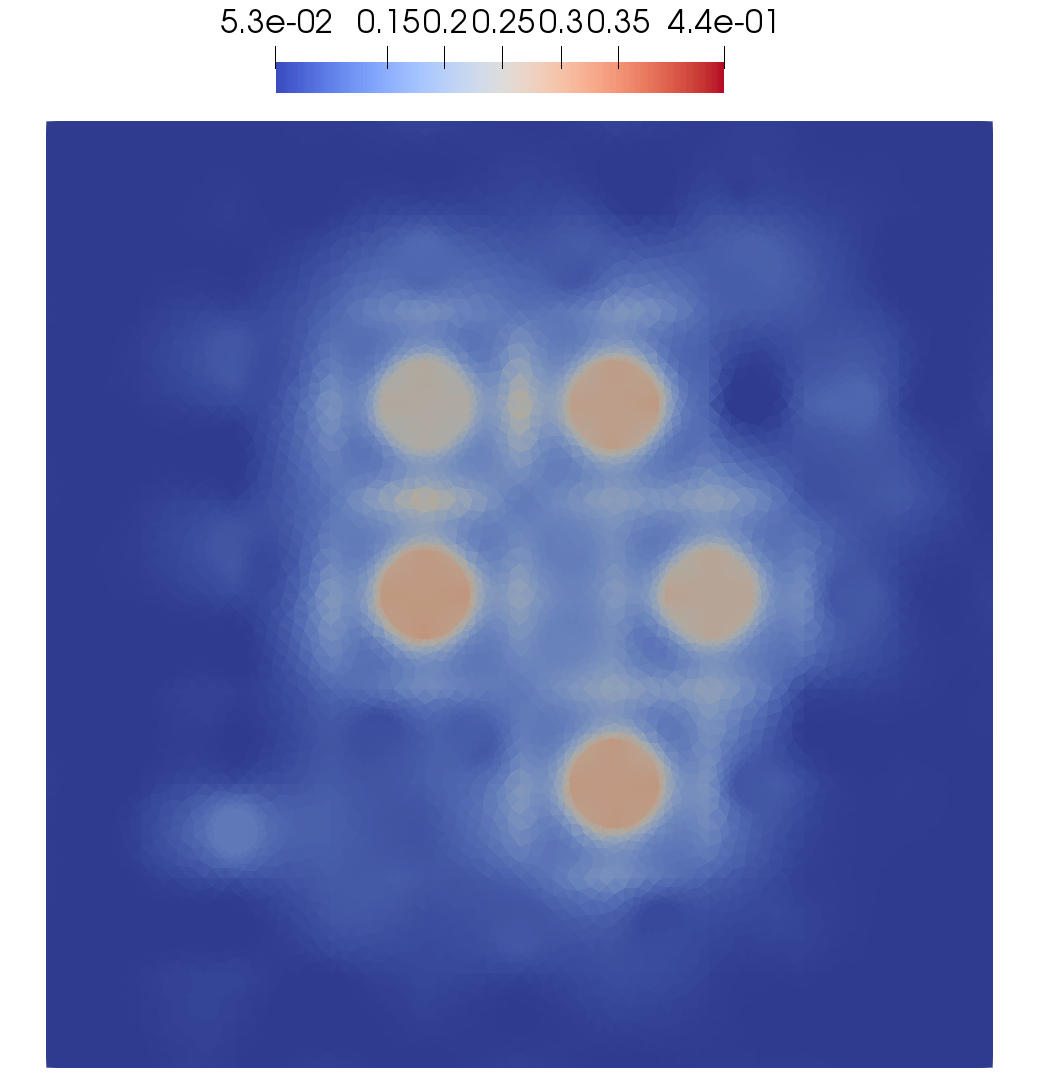}
\includegraphics[width=0.16\linewidth]{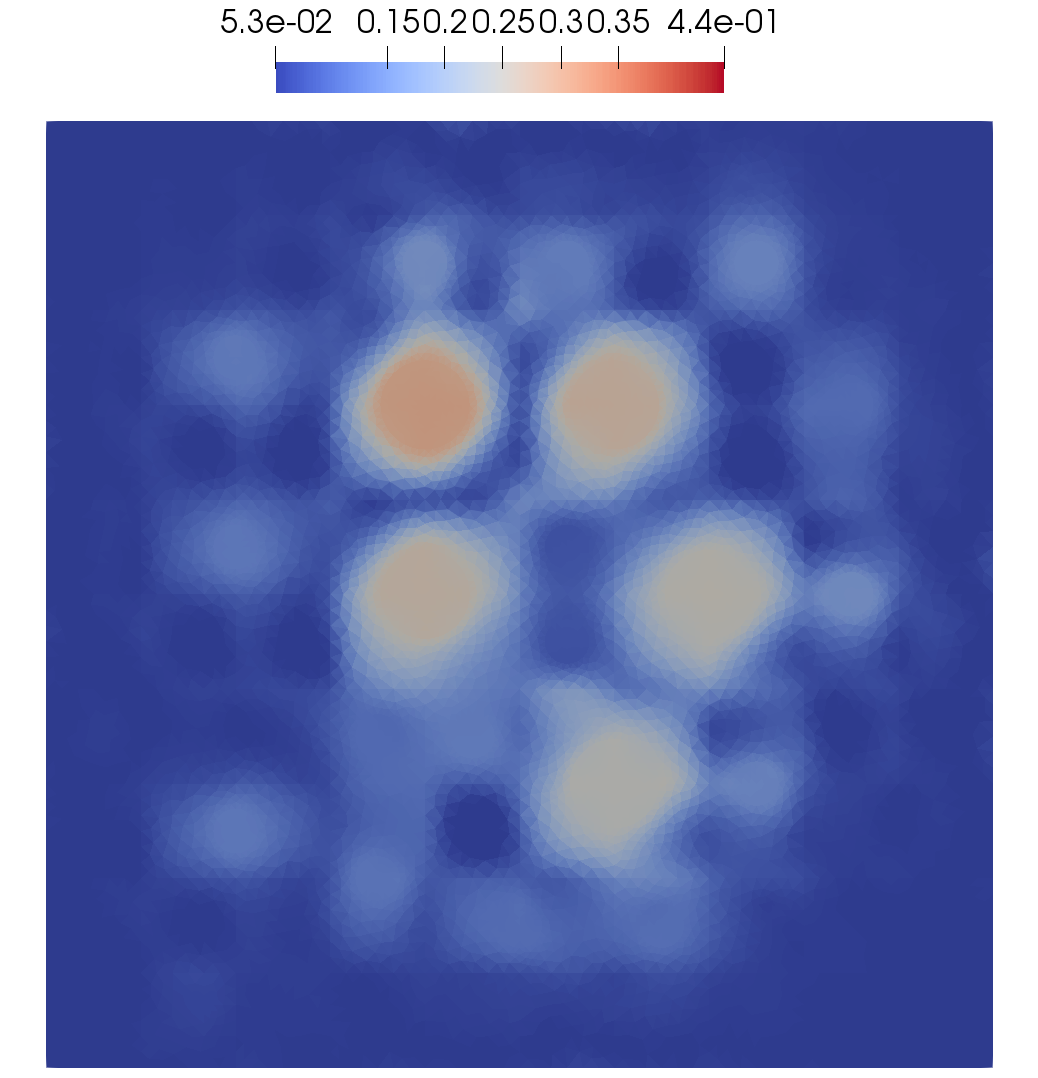}
\includegraphics[width=0.16\linewidth]{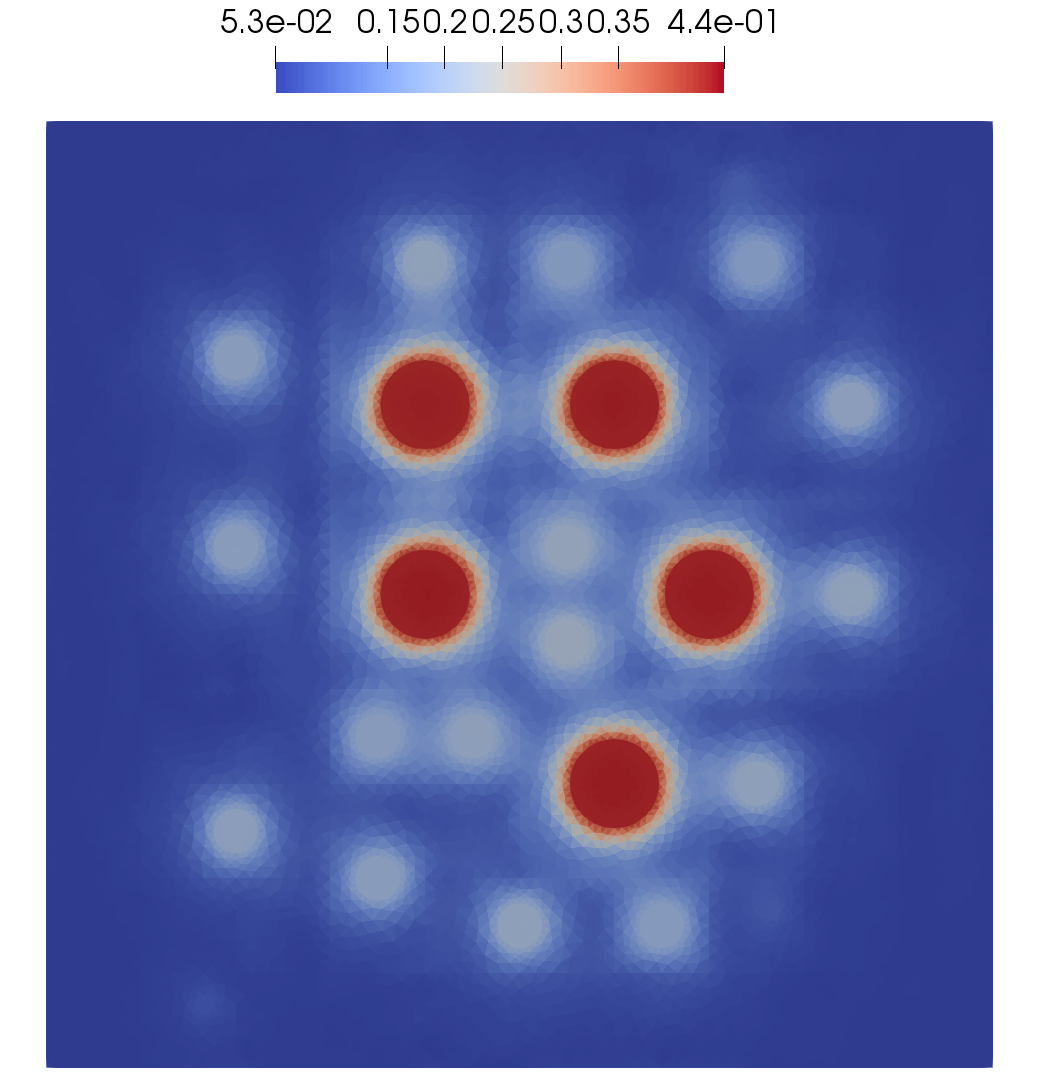}
\includegraphics[width=0.16\linewidth]{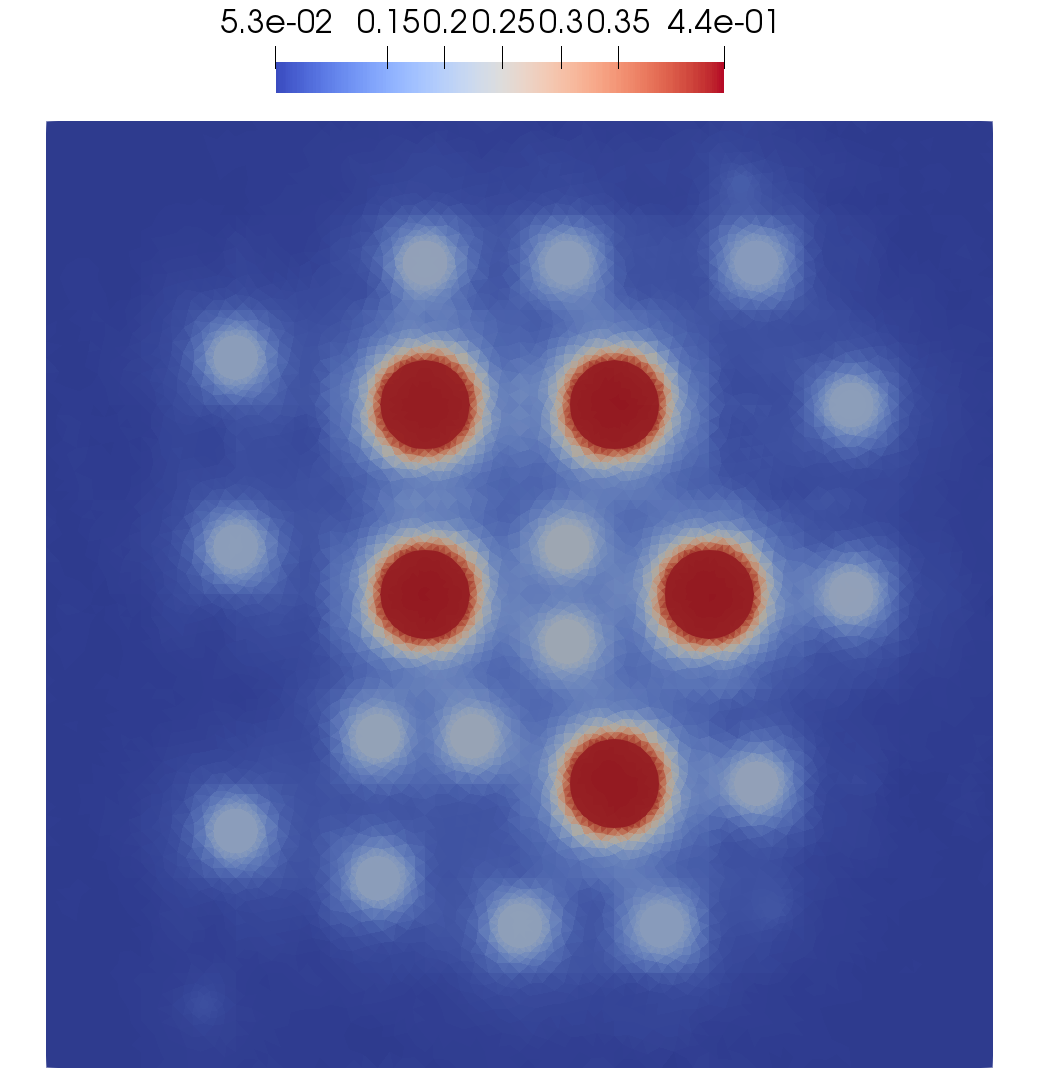}
\includegraphics[width=0.16\linewidth]{wf1}\\
\includegraphics[width=0.16\linewidth]{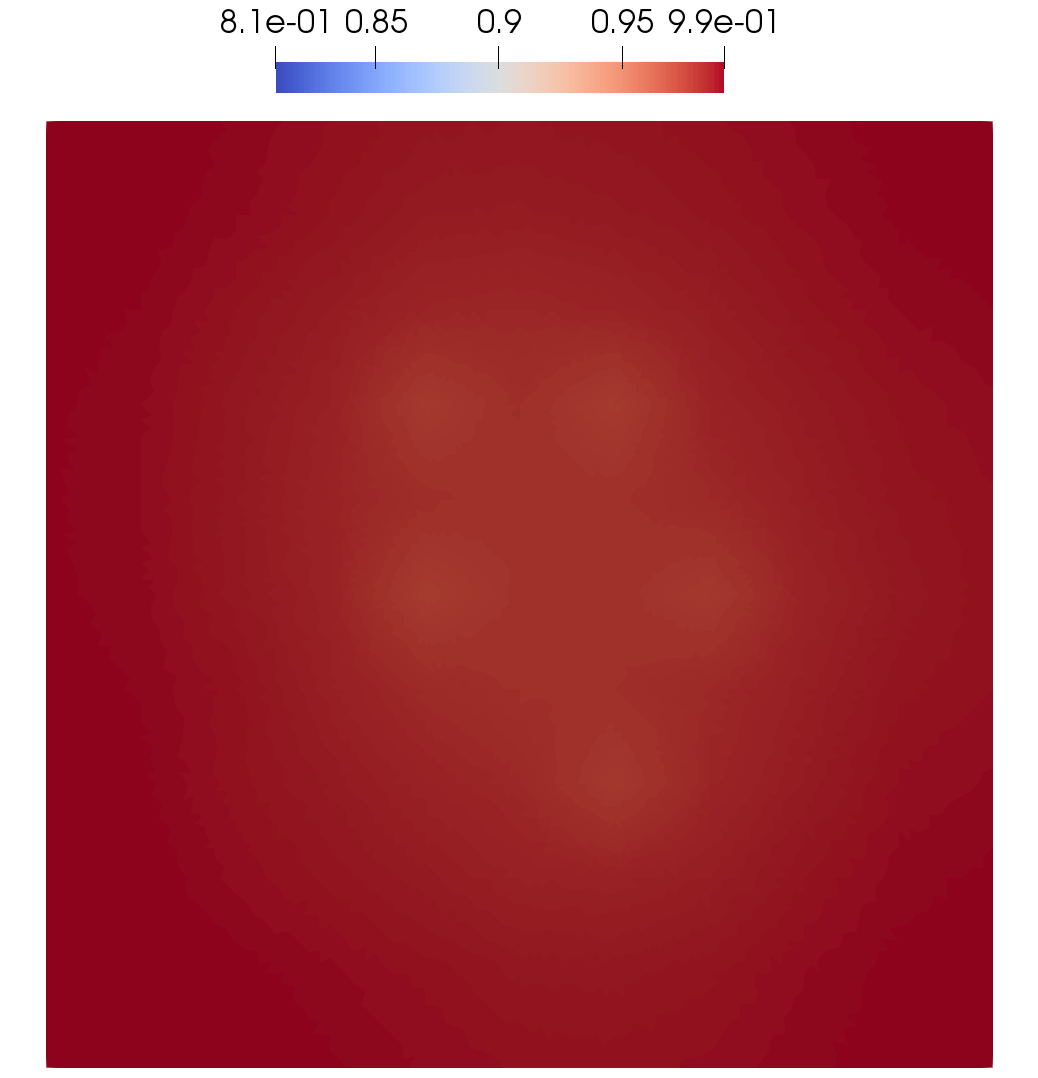}
\includegraphics[width=0.16\linewidth]{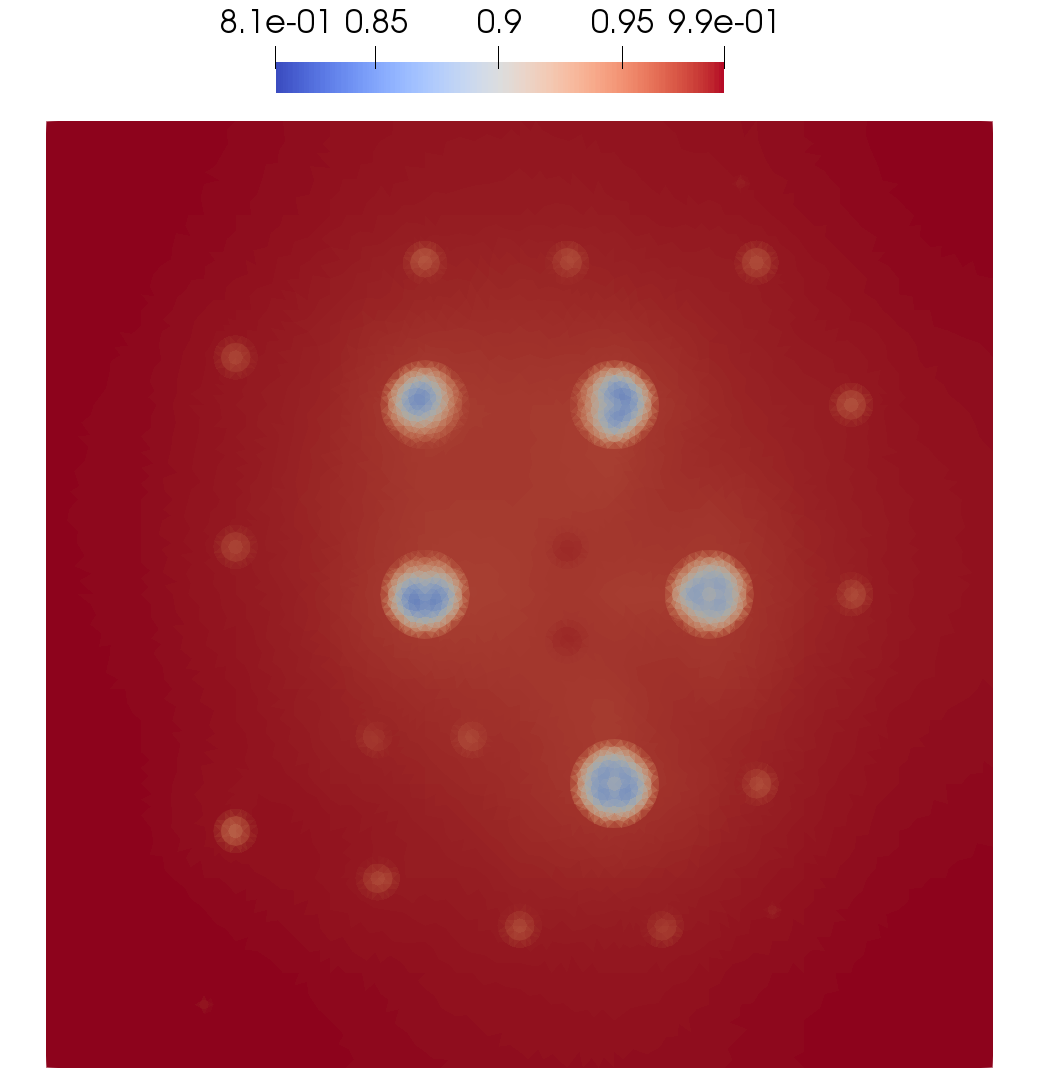}
\includegraphics[width=0.16\linewidth]{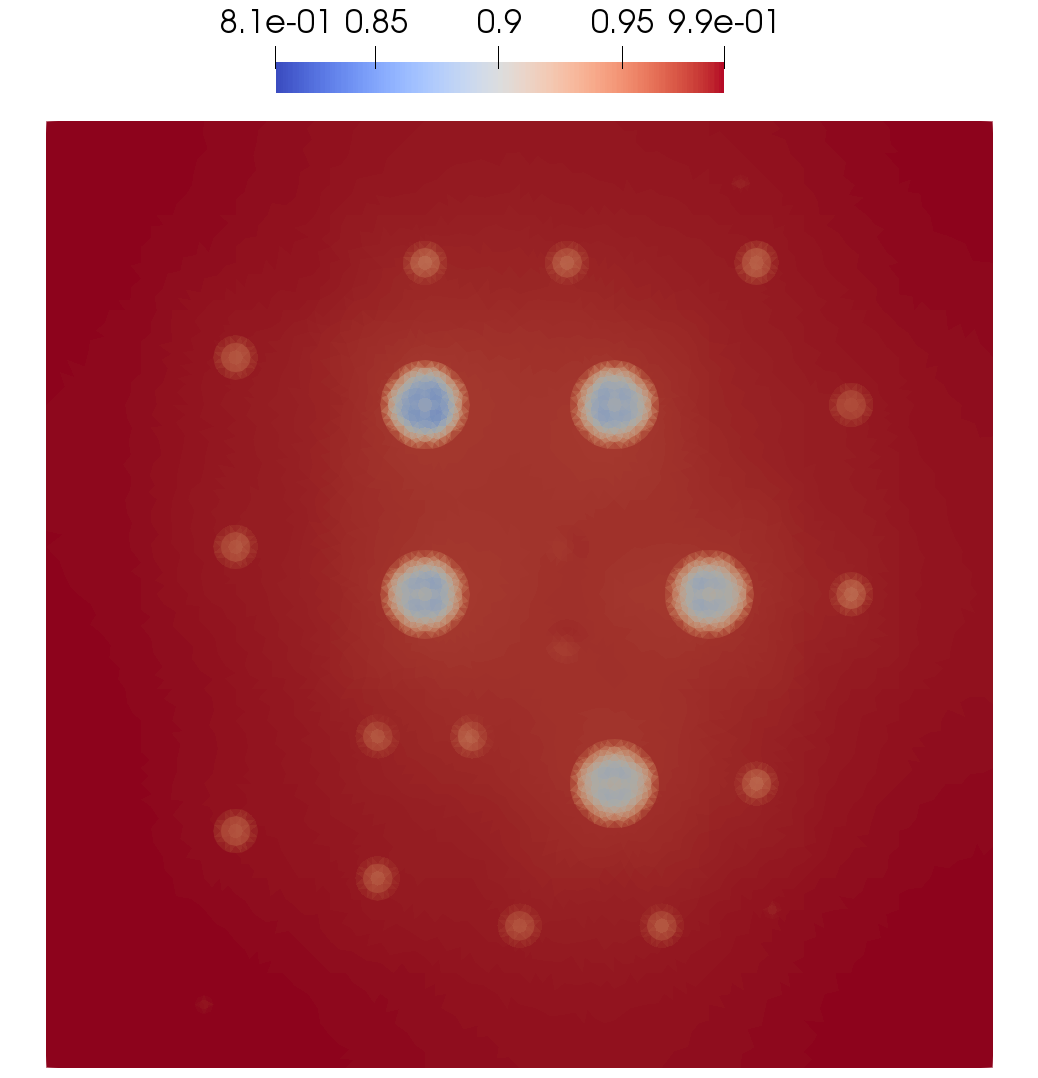}
\includegraphics[width=0.16\linewidth]{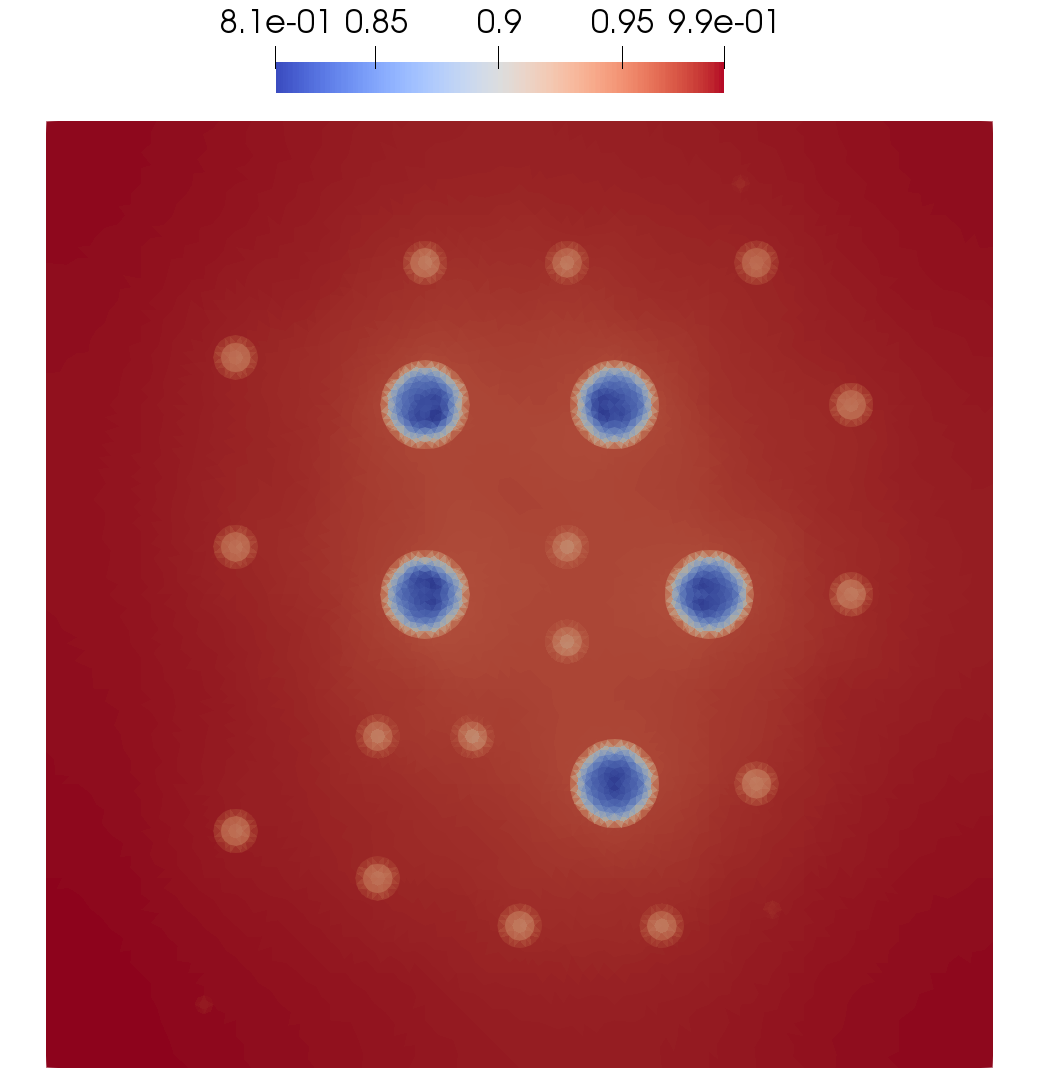}
\includegraphics[width=0.16\linewidth]{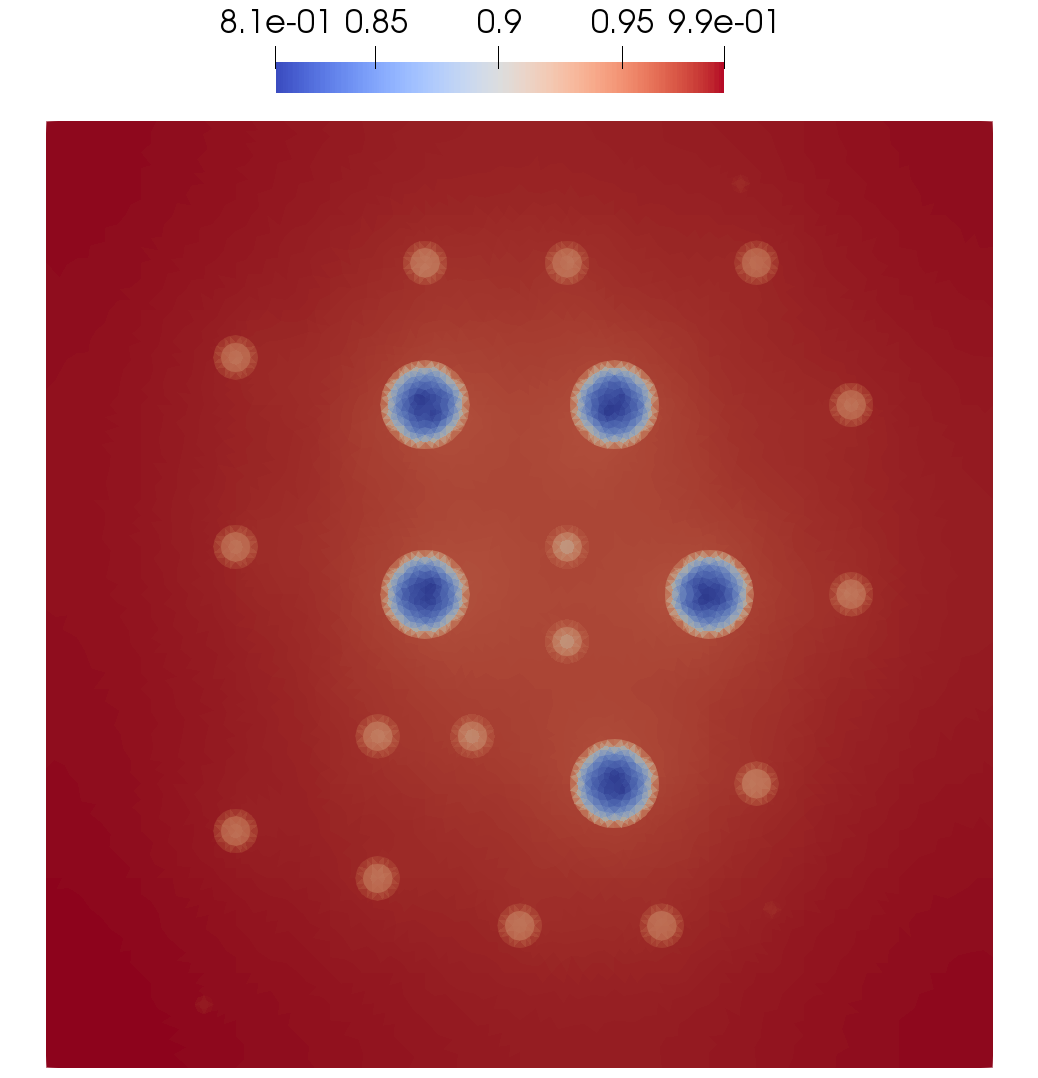}
\includegraphics[width=0.16\linewidth]{wf2}
\caption{Test 2a}
\end{subfigure}\\
\begin{subfigure}{1\textwidth}
\includegraphics[width=0.16\linewidth]{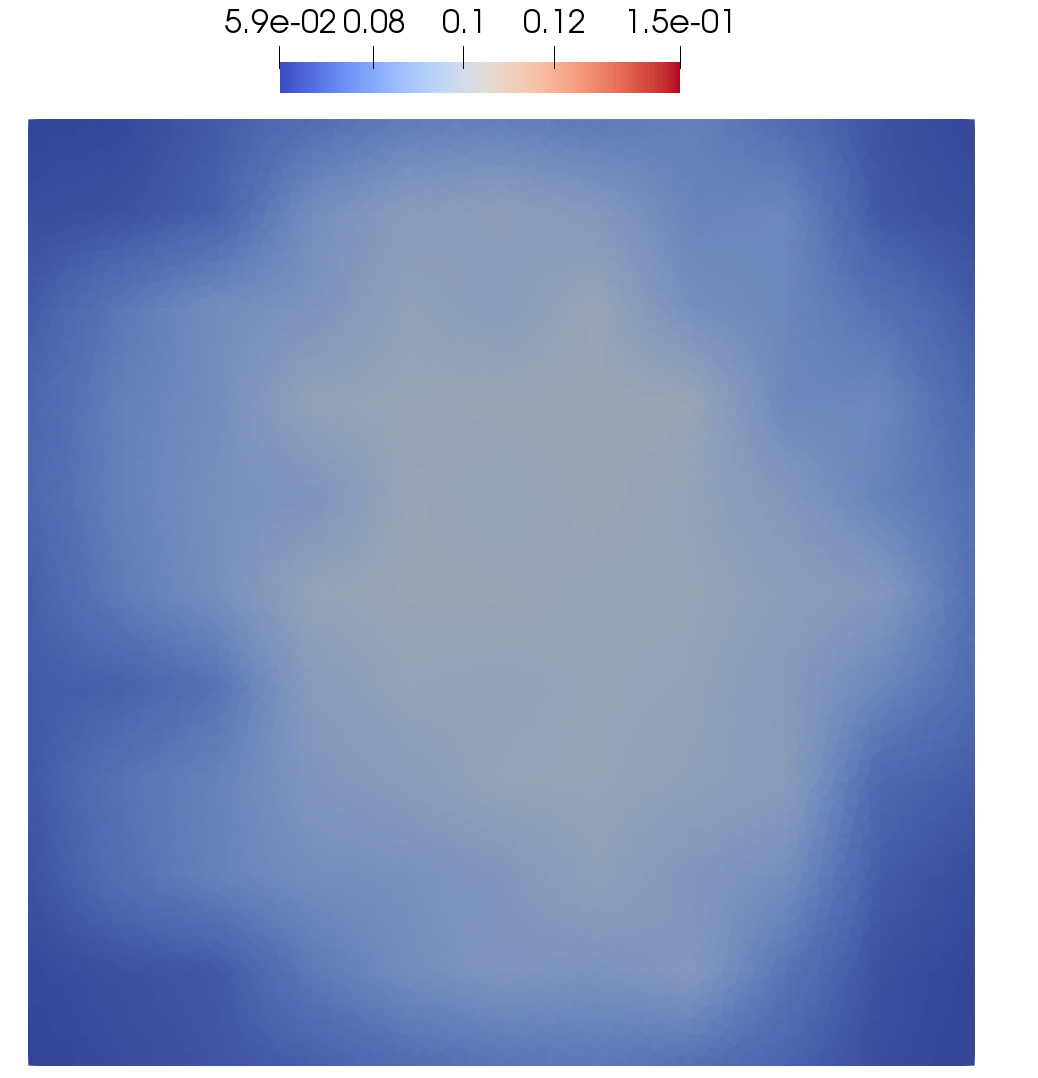}
\includegraphics[width=0.16\linewidth]{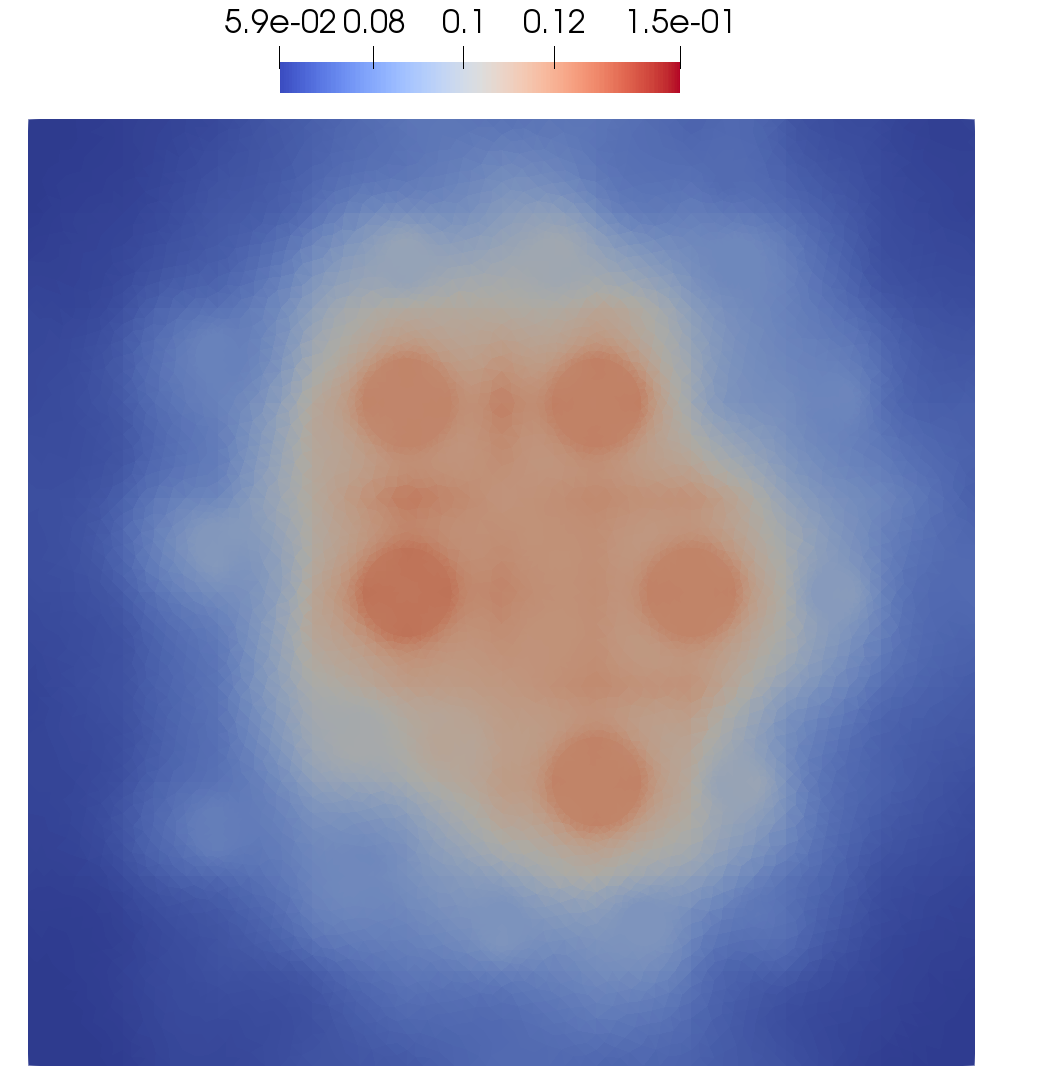}
\includegraphics[width=0.16\linewidth]{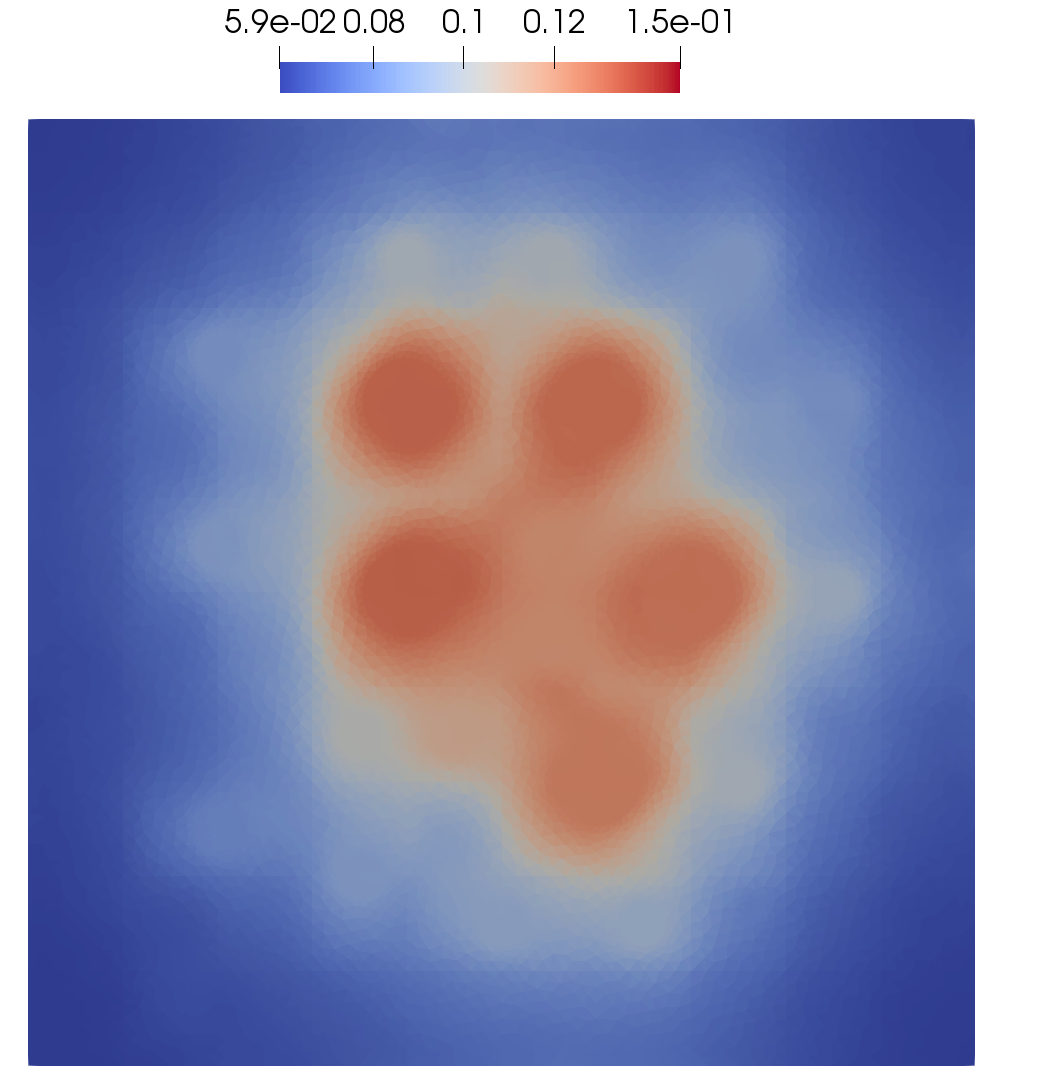}
\includegraphics[width=0.16\linewidth]{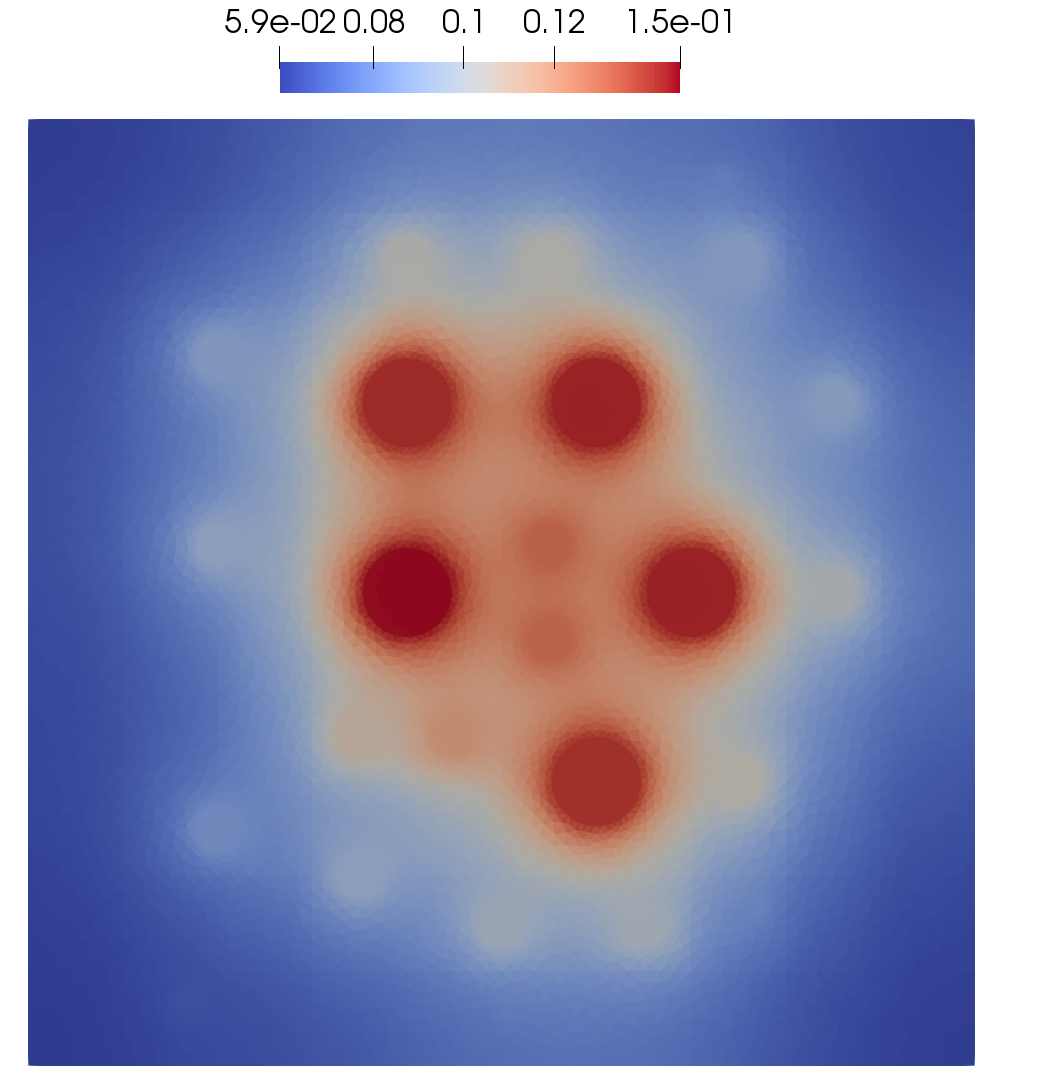}
\includegraphics[width=0.16\linewidth]{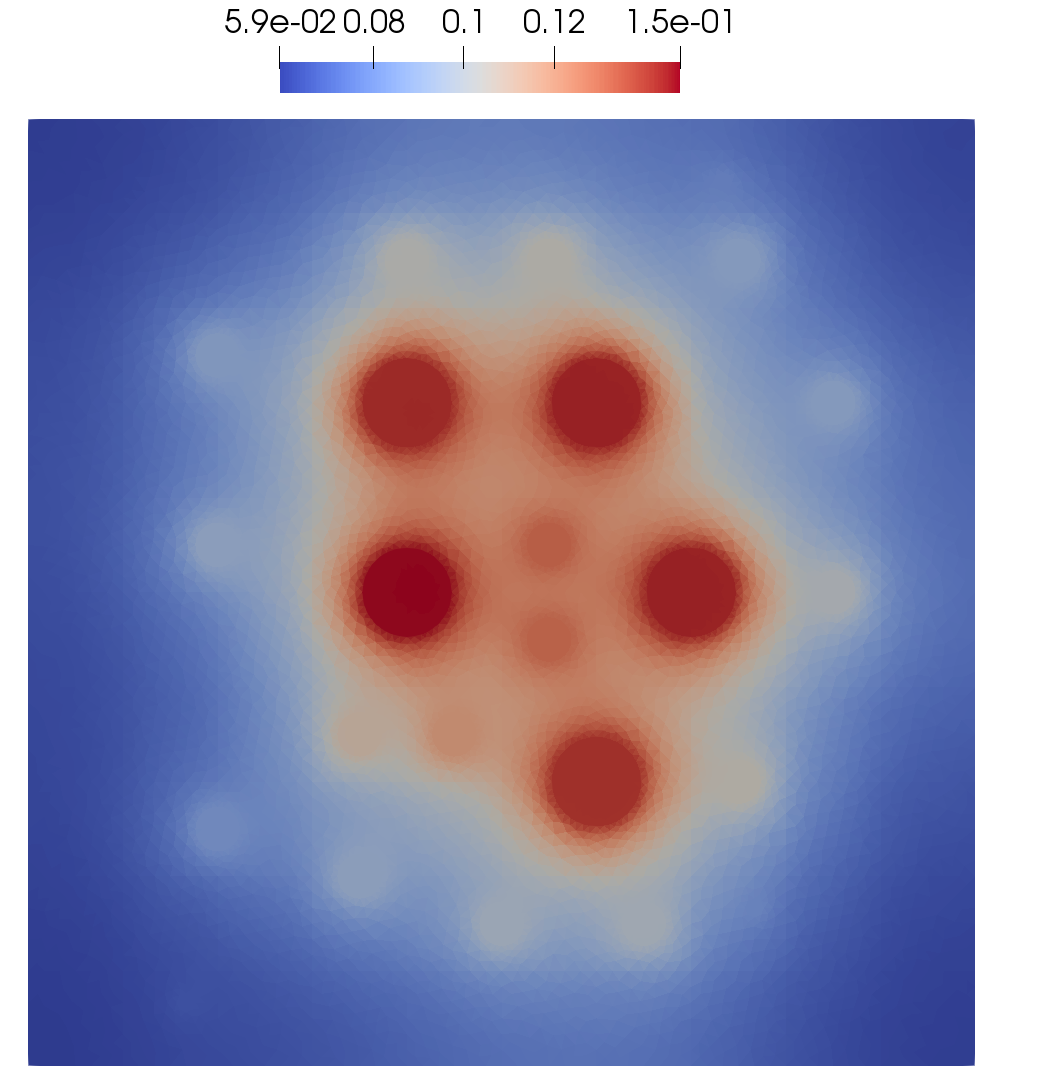}
\includegraphics[width=0.16\linewidth]{w2f1}\\
\includegraphics[width=0.16\linewidth]{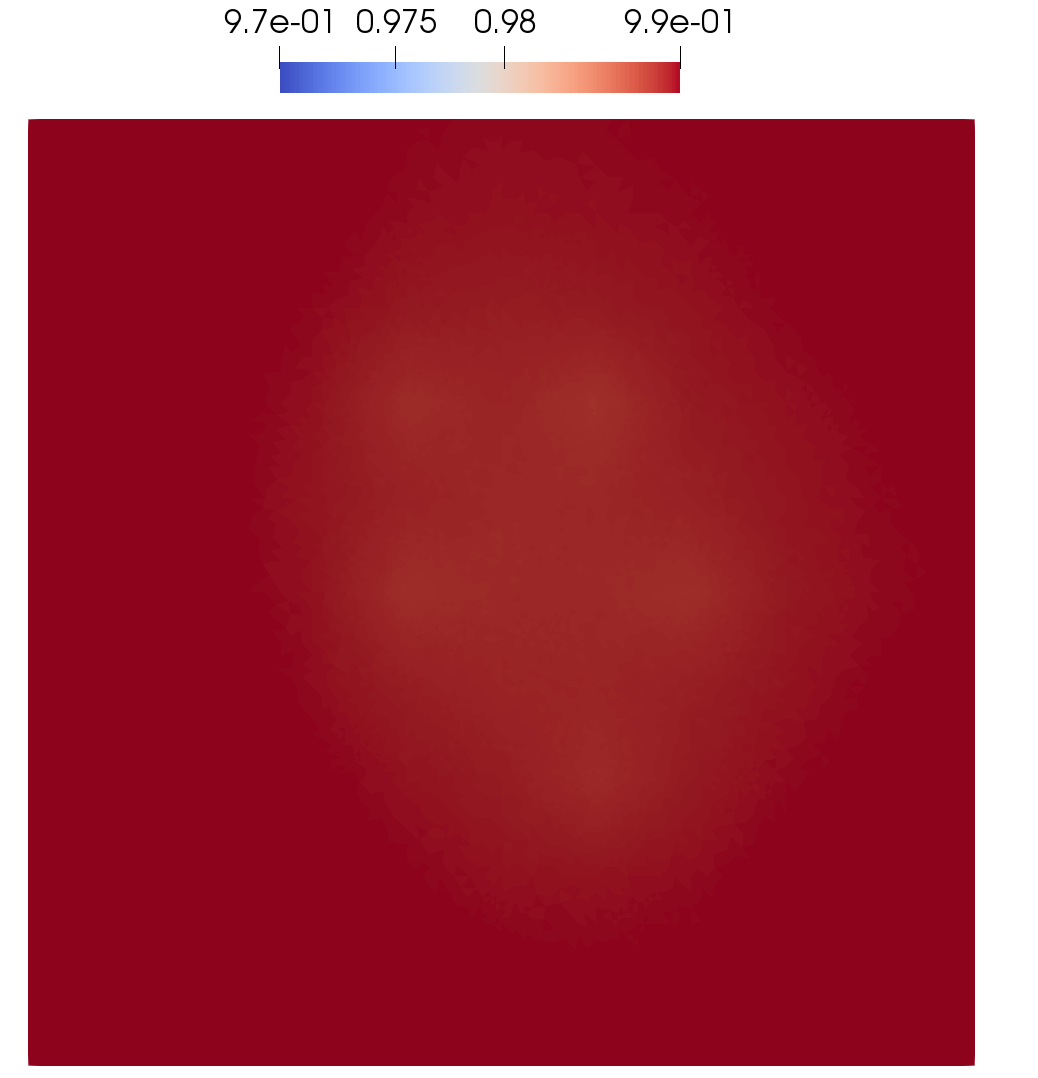}
\includegraphics[width=0.16\linewidth]{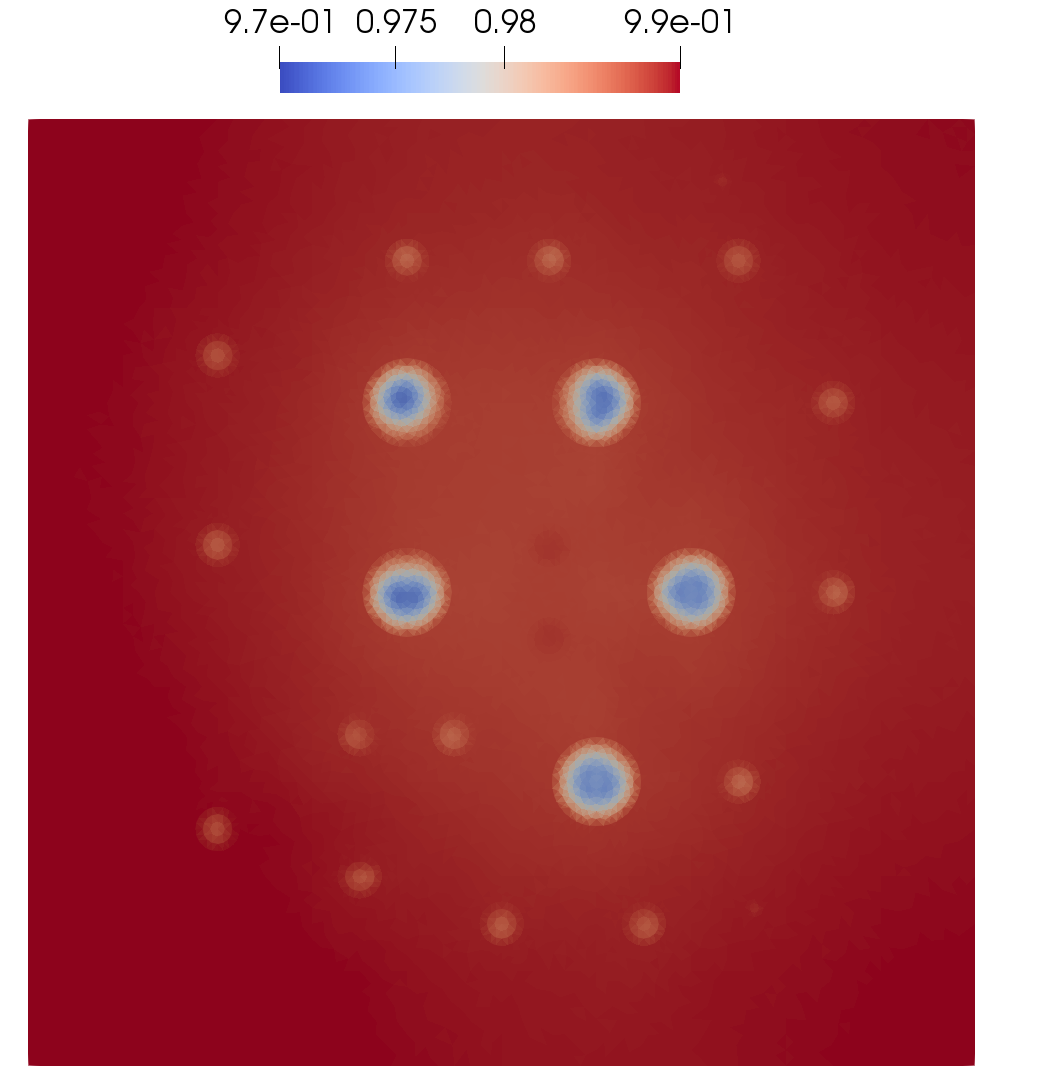}
\includegraphics[width=0.16\linewidth]{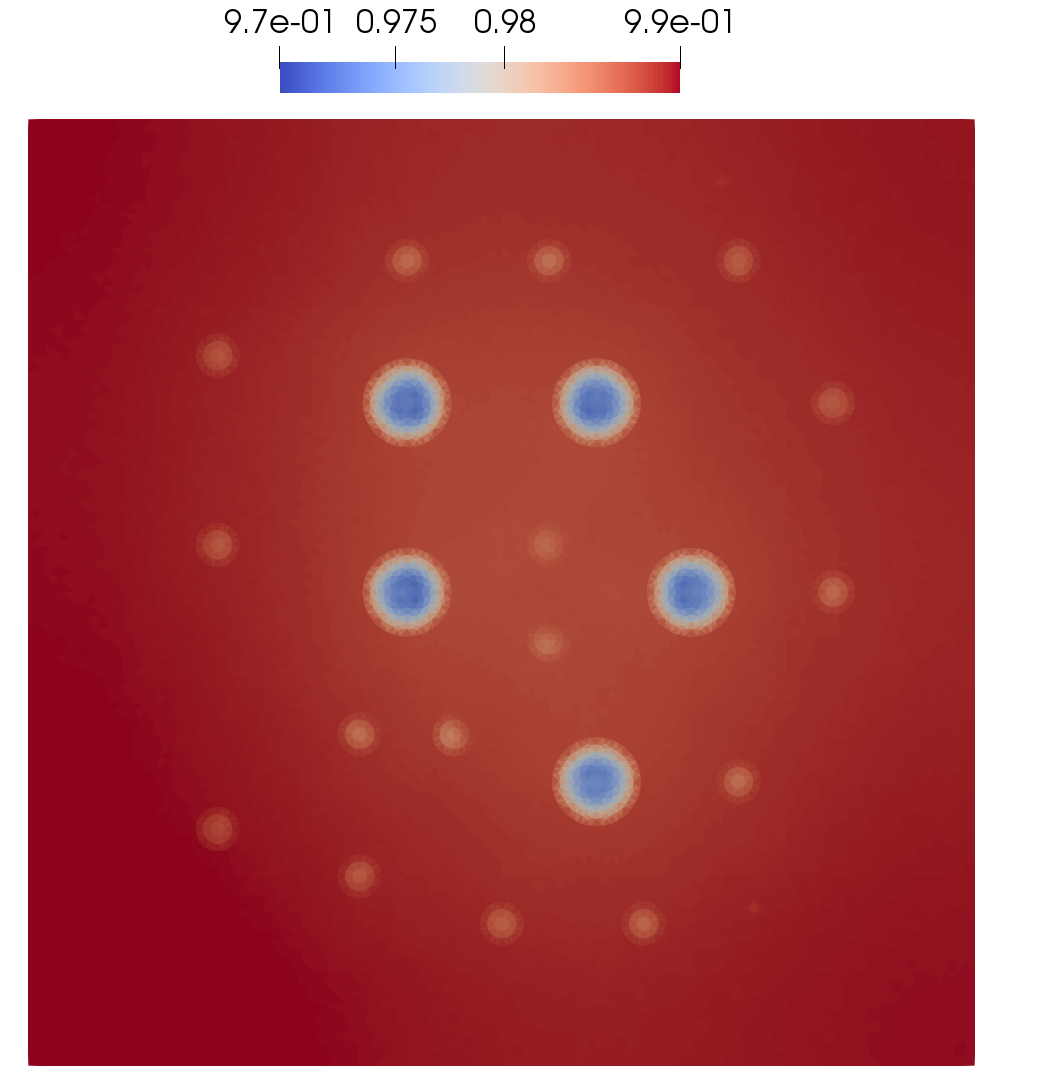}
\includegraphics[width=0.16\linewidth]{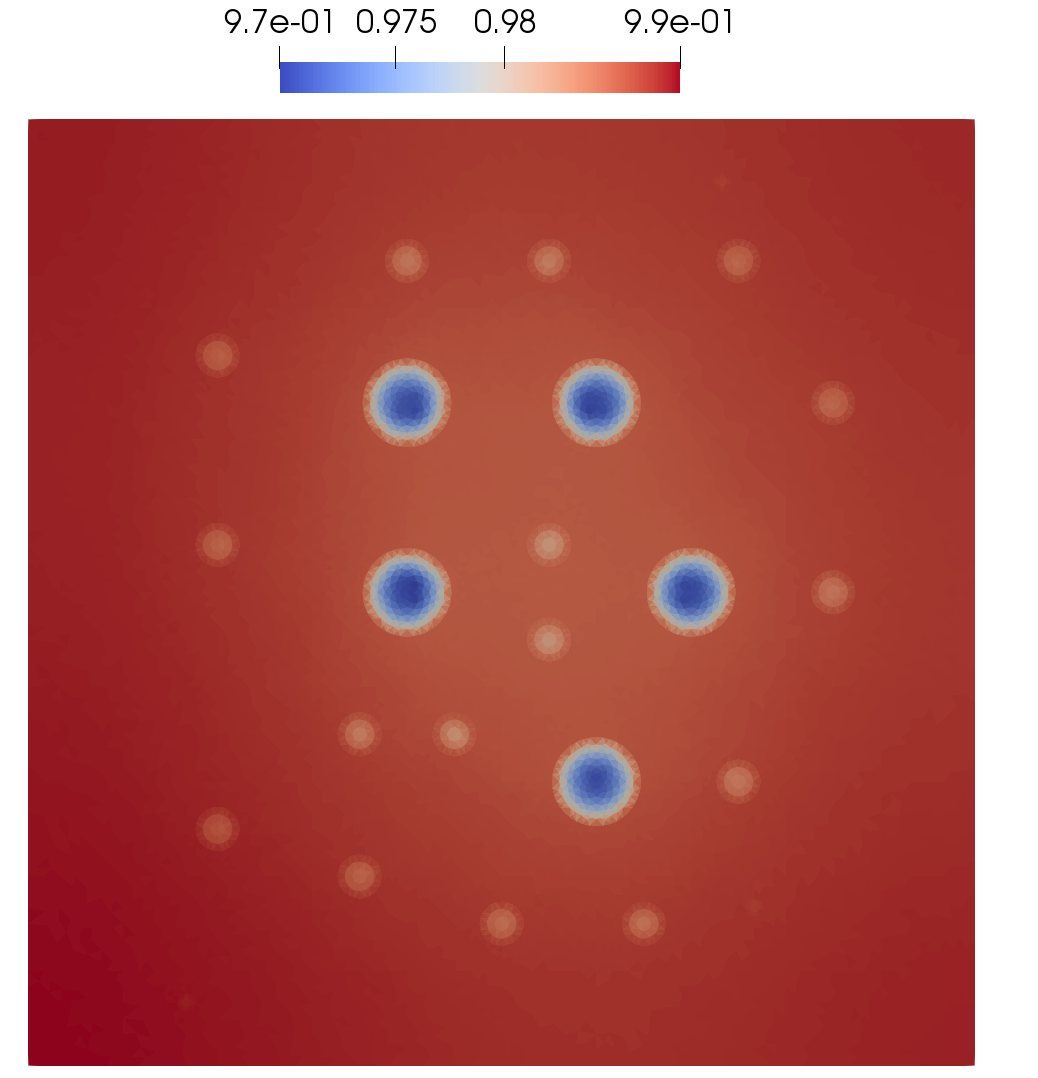}
\includegraphics[width=0.16\linewidth]{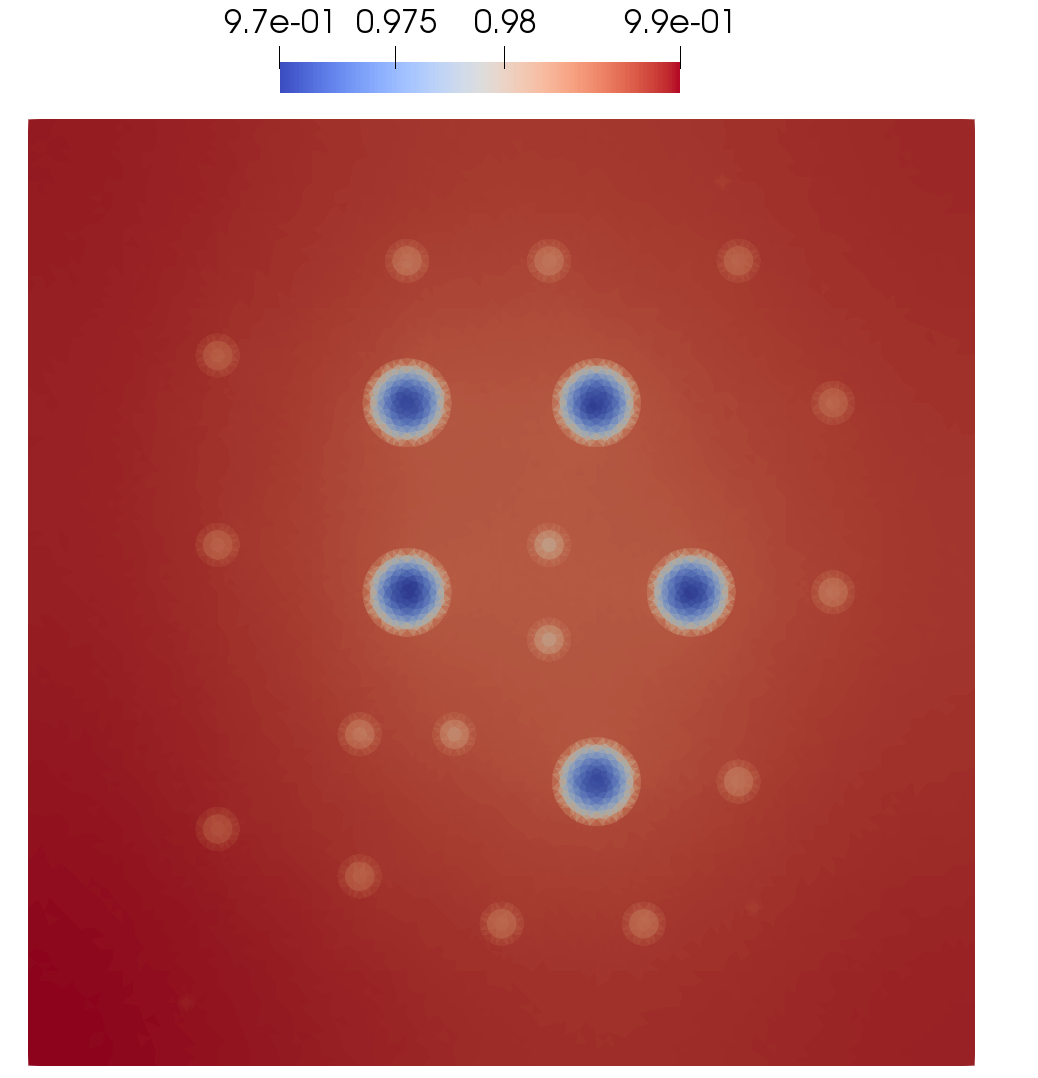}
\includegraphics[width=0.16\linewidth]{u2f2}
\caption{Test 2b}
\end{subfigure}
\caption{Test 2.  Solutions using multiscale solver for $M=1,2,4,6,8$ and fine grid (reference) solution (from left to right).  
First  and second rows: $u^1$ and $u^2$ for Test 2a (small diffusion) at final time. 
Third  and fourth rows: $u^1$ and $u^2$ for Test 2b (regular diffusion) at final time }
\label{sol2u}
\end{figure}

The diffusion coefficients for first and second species are assumed to be opposite. That is, we consider the case when the first species is less diffusive in the main subdomain (background) with $\varepsilon^1_m = \varepsilon^1_c/100$ while for the second species, we the background media to be more diffusive $\varepsilon^2_m = \varepsilon^2_c \cdot 100$. The corresponding first eight eigenvectors in local domain $\omega_i$ are given in Figure \ref{basis}, where we depicted local spectral solutions for first and second species. We observe that the basis functions are constructed in such a way that the variations of solutions in less diffusive subdomain are taken into account. Moreover, the first basis is constant. We note that, the multiscale basis functions are constructed by multiplication of the local eigenvectors to the corresponding linear partition of unity functions to preserve the continuity of the finite element approximation on the coarse grid. In our present method, the multiscale basis functions for small and regular diffusion are similar because we use similar contracts for both diffusion scenarios. The basis functions are based on diffusion operator without any effect on time and the reaction term. The main advantage of such approach is the ability to construct basis functions only once on offline stage and use it for small and regular diffusion terms with the same contrast and for any values of the reaction and time step size.  

To compare multiscale and fine-scale solutions, we calculate relative $L_2$ errors for first and second species population using the following formulas
\[
e_k = \left( \frac{\int_{\Omega} (u^k_{ref} - u^k_{ms})^2 \ dx}{\int_{\Omega} (u^k_{ref})^2 \ dx} \right)^{1/2}, 
\quad k = 1,2 
\]
where $u^k_{ref}$ refers to the reference solution and $u^k_{ms}$ is the solution using multiscale solver defined on the fine grid. 

\begin{figure}[h!]
\centering
\begin{subfigure}{0.48\textwidth}
\includegraphics[width=1\linewidth]{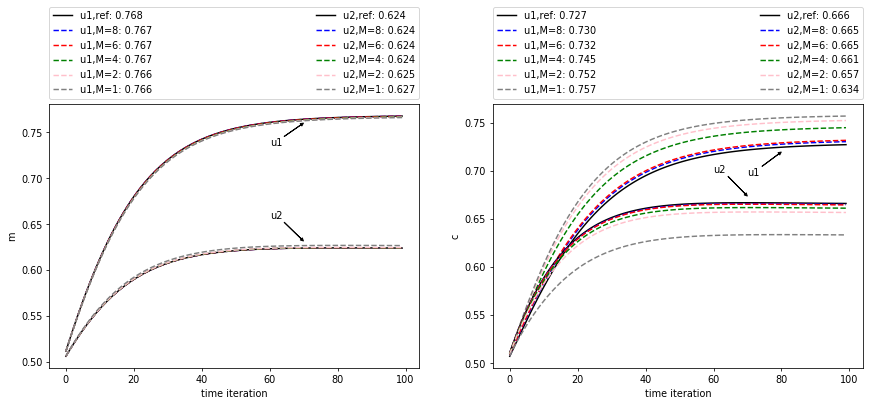}
\caption{Test 1a}
\end{subfigure}
\,\,\,\,
\begin{subfigure}{0.48\textwidth}
\includegraphics[width=1\linewidth]{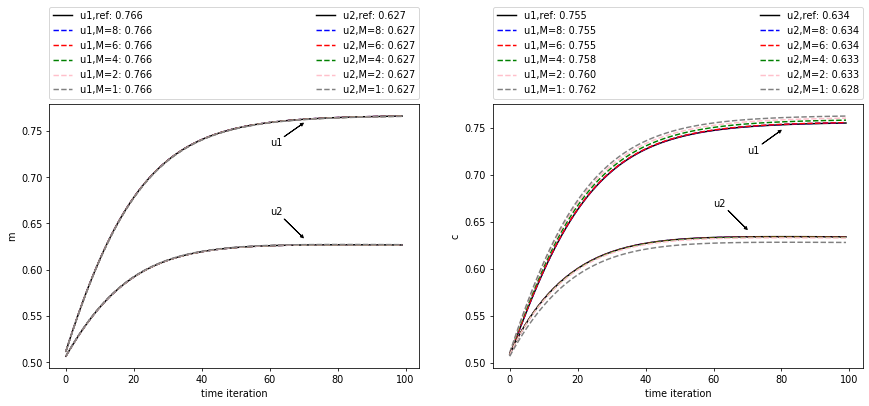}
\caption{Test 1b}
\end{subfigure}\\
\begin{subfigure}{0.48\textwidth}
\includegraphics[width=1\linewidth]{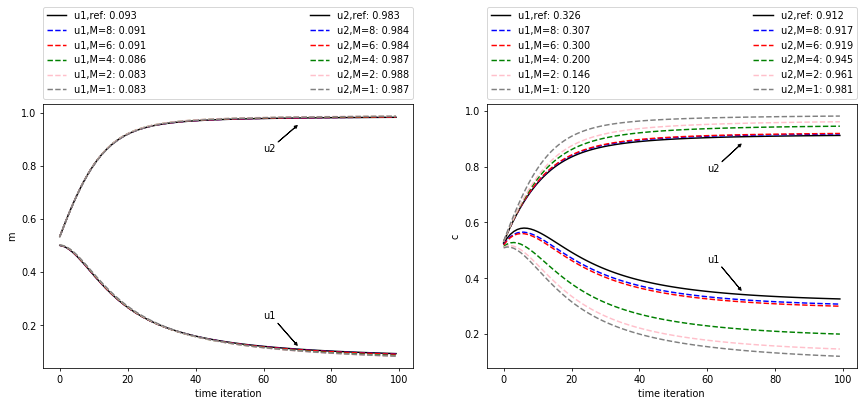}
\caption{Test 2a}
\end{subfigure}
\,\,\,\,
\begin{subfigure}{0.48\textwidth}
\includegraphics[width=1\linewidth]{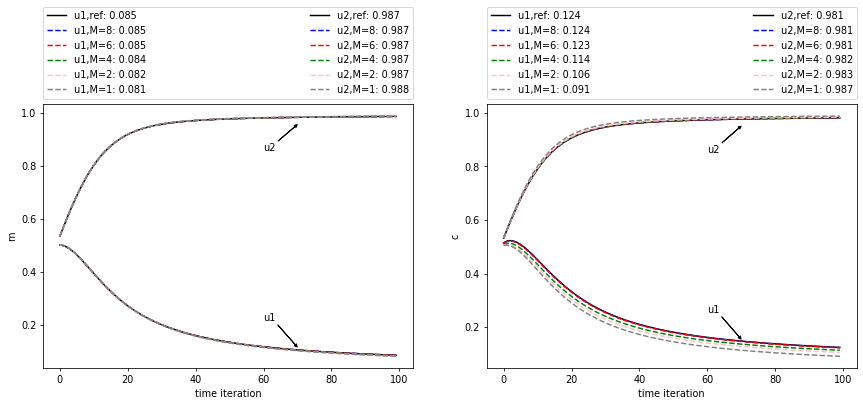}
\caption{Test 2b}
\end{subfigure}
\caption{Solution average for multiscale solver. $\bar{u}^1_m, \bar{u}^2_m$ in $\Omega_m$ (left) and $ \Omega_c$ (right) } 
\label{sol-ms-u}
\end{figure}

 Numerical solutions for the first and second species at the final time step for Test 1a/1b and Test 2a/2b are depicted in figures \ref{sol1u} and \ref{sol2u},, respectively.  Our multiscale method based solutions are represented in 1st, 2nd, 3rd, 4th, 5th columns,  for various number of multiscale basis functions ($M=1,2,4,6,8$). The fine grid (reference) solution is shown in the last column 6. We observe that the multiscale solution almost reproduces the reference solution, when we use a sufficient number of the multiscale basis functions ($M=6$) as seen from the figure. The values of the dynamic average solutions for $\bar{u}^1_m$, $\bar{u}^1_c$, $\bar{u}^2_m$ and $\bar{u}^2_c$(average for first and second species population in domain $\Omega_m$ and $\Omega_c$) are plotted in Figure \ref{sol-ms-u} at final time step.  Observe the convergence of the solutions with increasing number of the multiscale basis functions ($M$). It is clear that for $M=6$ one obtains multiscale solution that is almost indistinguishable to the reference solution for both test cases with small and regular diffusion. Moreover, we also see that the convergence is faster for larger diffusion scenarios.  

\begin{figure}[h!]
\centering
\begin{subfigure}{0.48\textwidth}
\includegraphics[width=1\linewidth]{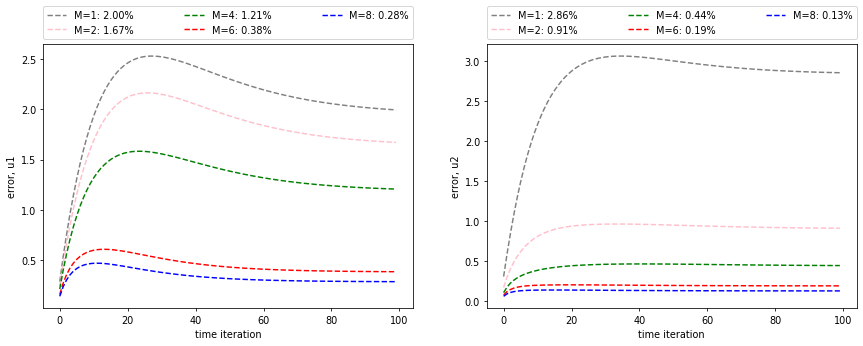}
\caption{Test 1a}
\end{subfigure}
\,\,\,\,
\begin{subfigure}{0.48\textwidth}
\includegraphics[width=1\linewidth]{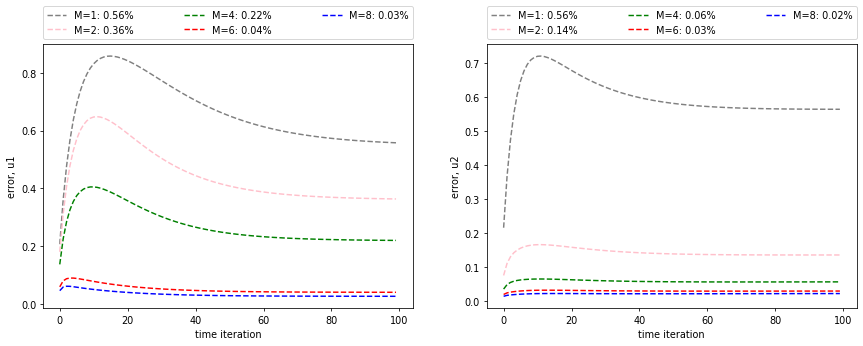}
\caption{Test 1b}
\end{subfigure}\\
\begin{subfigure}{0.48\textwidth}
\includegraphics[width=1\linewidth]{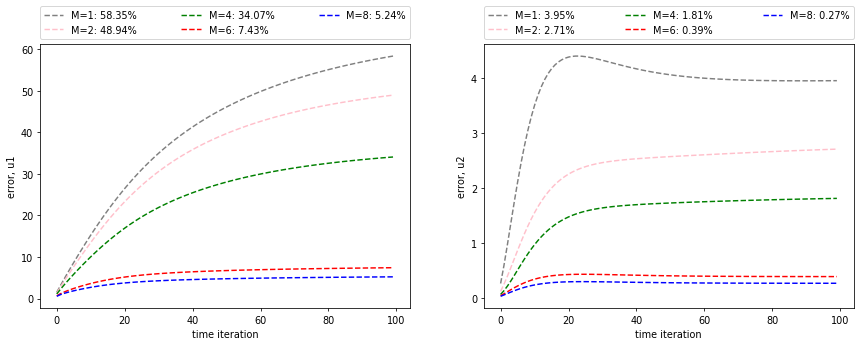}
\caption{Test 2a}
\end{subfigure}
\,\,\,\,
\begin{subfigure}{0.48\textwidth}
\includegraphics[width=1\linewidth]{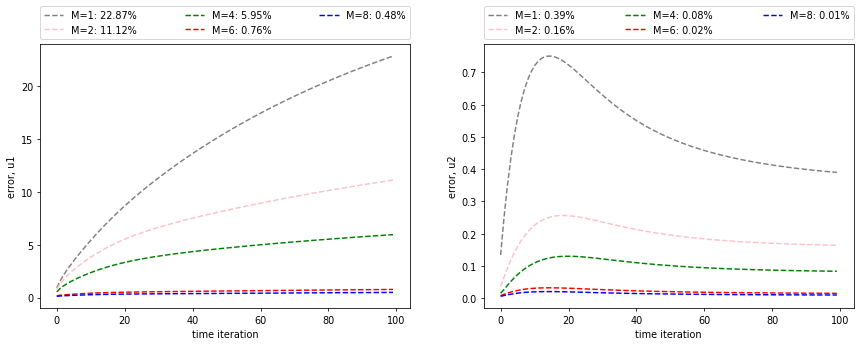}
\caption{Test 2b}
\end{subfigure}
\caption{Solution average errors for multisscale solver. $u^1$ (left) and $u^2$ (right)} 
\label{sol-ms-e}
\end{figure}

\begin{table}[h!]
\centering
\begin{tabular}{|cc|ccc|ccc|}
\hline
&
& $e_1$ (\%)  & $e_2$( \%) & time(sec) 
& $e_1$ (\%)  & $e_2$ (\%) & time(sec)  \\
 \raisebox{1.5ex}[0cm][0cm]{$M$ }  
&  \raisebox{1.5ex}[0cm][0cm]{$DOF$ }  &
 \multicolumn{3}{|c|}{Test 1a}&
 \multicolumn{3}{|c|}{Test 1b}\\ 
\hline
reference 	& 69948 	& - 	& - 	& 84.67 	& - 	& - 	& 229.95 \\
1 	& 121 	& 1.995 	& 2.857 	& 25.438 	& 0.557 	& 0.564 	& 24.834 \\
2 	& 242 	& 1.671 	& 0.911 	& 26.031 	& 0.363 	& 0.136 	& 28.272 \\
4 	& 484 	& 1.206 	& 0.444 	& 26.600 	& 0.220 	& 0.056 	& 27.043 \\
6 	& 726 	& 0.382 	& 0.190 	& 27.854 	& 0.040 	& 0.029 	& 28.030 \\
8 	& 968 	& 0.283 	& 0.127 	& 51.025 	& 0.026 	& 0.022 	& 30.520 \\
\hline
& &
 \multicolumn{3}{|c|}{Test 2a}&
 \multicolumn{3}{|c|}{Test 2b}\\ 
\hline
reference 	& 69948 	& - 	& - 	& 121.18 	& - 	& - 	& 405.20 \\
1 	& 121 	& 58.346 	& 3.951 	& 25.833 	& 22.866 	& 0.390 	& 25.591 \\
2 	& 242 	& 48.945 	& 2.706 	& 26.059 	& 11.115 	& 0.164 	& 25.982 \\
4 	& 484 	& 34.067 	& 1.811 	& 26.725 	& 5.953 	& 0.084 	& 26.755 \\
6 	& 726 	& 7.435 	& 0.389 	& 27.531 	& 0.764 	& 0.015 	& 27.737 \\
8 	& 968 	& 5.237 	& 0.267 	& 30.426 	& 0.482 	& 0.010 	& 29.409 \\
\hline
\end{tabular}
\caption{Errors at final time and time for the solution}
\label{table-ms}
\end{table}

\begin{figure}[h!]
\centering
\begin{subfigure}{0.48\textwidth}
\includegraphics[width=1\linewidth]{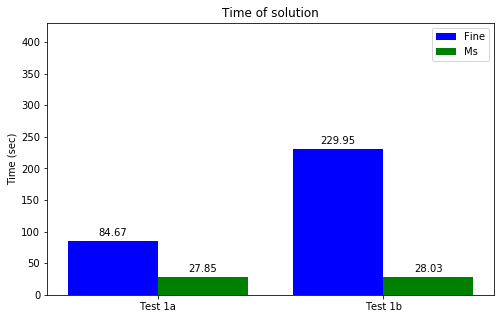}
\caption{Test 1}
\end{subfigure}
\,\,\,\,
\begin{subfigure}{0.48\textwidth}
\includegraphics[width=1\linewidth]{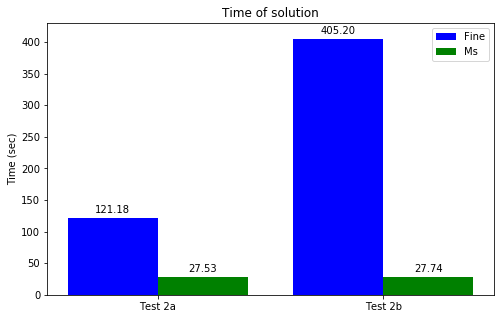}
\caption{Test 2}
\end{subfigure}
\caption{Simulation times for fine-scale and multiscale models} 
\label{sol-ms-t}
\end{figure}

Figure \ref{sol-ms-u} depicts the solution average errors for first and second species. It is evident from the figure that the errors reduce as we increase number of multiscale basis functions in a similar fashion. 
Table \ref{table-ms}, lists the errors at the final time for Test 1a/1b and Test 2a/2b.  We observe that one can obtain an accurate numerical solutions using 6 multiscale basis functions for all cases in a given heterogeneity distribution on a $10 \times 10$ coarse grid. The solution time for both reference solution (on the fine grid) as well as the multiscale solution (using the multiscale solver) with the given number of basis functions is also provided in the tables. In general, sufficient number of multiscale basis functions depends on the size of the coarse grid and heterogeneity distribution in each local domain. Additionally, the number of multiscale basis functions can be calculated adaptively based on the eigenvalues.

The reference solution is obtained using the finite volume method in which case the number of unknowns equals the number of fine grid cells, that is, $DOF_h = 69948$. Implementation of the solver (online stage) is performed on python programming language and PETSc solver (petsc4py). For multiscale and fine scale solutions, we use a GMRES iterative solver with ILU preconditioner. We remark that in the present work we did not develop and explore iterative solvers for fine-scale and multiscale systems. Offline calculations (multiscale basis construction and generation of the projection matrix) are implemented using C++ and local calculations in $\omega_i$ are performed in parallel, because basis generations in each local domain can be done independently. Solution time data presented in Table \ref{table-ms} is a measurement of the solution time exclusively, without any additional time taken for saving the solution for visualization and error calculations. From the values recorded in the tables, we notice that the computation is faster using multiscale solver due to the reduction of the system size in all test cases. Moreover, we observe that the solution time using multiscale solver doesn't depend on the diffusion values and test parameters where as simulations on the fine grid take longer time for larger diffusion. This is because the latter (fine grid) method affects the matrix arising from the linear system and also depends on test parameters (right hand side in the linear system). For instance, the simulation time on fine grid is $405$ seconds (with $DOF_h = 69948$) while the simulation time using multiscale solver with $M=6$ ($DOF_H = 726$) is just 27 seconds with an error less than 1 \% for Test 2b. The computational efficiency of the multiscale solver with $M=6$ is further illustrated in Figure \ref{sol-ms-t}.

\section{Conclusion}

A constructive reduced order model for computing numerical solutions of reaction-diffusion systems in heterogeneous domains is presented herein. The fine-scale solver has been built using finite volume method. 
Fully implicit (FI) coupled and semi-implicit (SI) uncoupled schemes have been utilized for time approximations. The numerical results show that the simulation running time using an uncoupled scheme is significantly shorter, with good accuracy, using small time stepping size in the reaction-diffusion equations. Further reduction of the simulation time has been achieved by reducing the size of the discrete system via an efficient multiscale solver developed based on the Generalized Multiscale Finite Element Method (GMsFEM). We have also presented the construction of multiscale basis functions, in each local domain, based on the diffusion operators for each species separately. The basis functions are in turn used to build (using projection) coarse scale system with smaller number of unknowns ($DOF_H << DOF_h$). The coarse grid reduced order model is then solved leading to the reconstruction of the fine-scale solution for accurate reaction term approximation. It is therefore imperative to update the right-hand side (reaction term) based on the fine grid solution information and project it back to the coarse grid. Our approach reveals that the multiscale basis functions do not depend on the reaction term and can be used for multi-parameter simulations with various values of reaction coefficients and different time step sizes. 

The multiscale solver has been successfully applied to a two-species competition reaction-diffusion model in a heterogeneous domain. Numerical results are presented for a variety of scenarios by assigning values of the parameters. Specifically, two test problems with small and regular diffusion coefficients in two-dimensional formulation are discussed in detail. The solution average, solutions, solution times and errors are presented numerically. One striking feature that the influence of number of multiscale basis functions on the accuracy and ability to work with different values of the diffusion coefficients is demonstrated in our work. Very good results are obtained with small errors by employing a sufficient number of the multiscale basis functions. One of the implications of our numerical results is that the dominance of a particular species over the other may be dictated by diffusion. Finally, the proposed multiscale solver is accurate, faster, and computationally effective and can handle many complex multi-component reaction-diffusion systems.

\color{black}

\bibliographystyle{plain}
\bibliography{lit}

\end{document}